\documentclass{IEEE-Con-Sys-mag}

\jvol{XX}
\jnum{XX}
\paper{8}
\jmonth{July}
\jname{IEEE CONTROL SYSTEMS}
\pubyear{2022}

\title{Towards a Theory of Control Architecture\stitle{A quantitative framework for layered multi-rate control}}
\author{Nikolai Matni, Aaron D. Ames, and John C. Doyle}
\affil{N. Matni (\href{mailto:nmatni@seas.upenn.edu}{nmatni@seas.upenn.edu}), Dept. of Electrical and Systems Engineering, University of Pennsylvania, Philadelphia, USA.\\
A. D. Ames (\href{mailto:ames@caltech.edu}{ames@caltech.edu}) and J. C. Doyle (\href{mailto:doyle@caltech.edu}{doyle@caltech.edu}), Control and Dynamical Systems, California Institute of Technology, Pasadena, USA. }

\usepackage{booktabs}
\usepackage{array}
\newcommand{\PreserveBackslash}[1]{\let\temp=\\#1\let\\=\temp}
\newcolumntype{C}[1]{>{\PreserveBackslash\centering}p{#1}}

\usepackage{caption}
\usepackage{amsthm}
\usepackage{amssymb}
\usepackage{amsmath}
\usepackage{mathrsfs}
\usepackage{hyperref}
\usepackage{makecell}
\usepackage[most]{tcolorbox}
\hypersetup{
    unicode=false,          
    pdftoolbar=true,        
    pdfmenubar=true,        
    pdffitwindow=false,     
    pdfstartview={FitH},    
    pdfnewwindow=true,      
    colorlinks=true,       
    linkcolor=blue,          
    citecolor=blue,        
    filecolor=magenta,      
    urlcolor=red,           
    breaklinks=true
}
\usepackage{cleveref}
\usepackage{xcolor}
\usepackage[sort, numbers, compress]{natbib}
\usepackage{here}
\usepackage{caption}
\usepackage{subcaption}
\usepackage{graphicx}
\usepackage{graphbox}

\newcommand{\Projx}{\pi_x}

\newcommand{\Projv}{\pi^v}

\newcommand{\floor}[1]{\left\lfloor #1 \right\rfloor}

\newcommand{\R}{\mathbb{R}}

\newcommand{\edit}[1]{{#1}}

\DeclareMathOperator{\argmin}{argmin}

\newcommand{\norm}[1]{\lVert #1 \rVert}

\renewcommand{\argmin}{\operatornamewithlimits{argmin}} 

\DeclareMathOperator*{\esssup}{esssup}

\newtheorem{my_theorem}{Theorem}
\newenvironment{my_thm}
  { \begin{tcolorbox}[
 colframe=yellow!70!white,
 colback=yellow!17!white,
 arc=8pt,
 breakable,
 left=1pt,right=1pt,top=1pt,bottom=1pt,
 boxrule=0.3pt,
 ]
\begin{my_theorem}}
  {\end{my_theorem}\end{tcolorbox}}

\theoremstyle{definition}
\newtheorem{my_example}{Example}
\newenvironment{my_exam}
  { \begin{tcolorbox}[
 colframe=green!70!white,
 colback=green!17!white,
 arc=8pt,
 breakable,
 left=1pt,right=1pt,top=1pt,bottom=1pt,
 boxrule=0.3pt,
 ]
\begin{my_example}}
  {\end{my_example}\end{tcolorbox}}

\theoremstyle{theorem}\newtheorem*{no-theorem}{Theorem}

\renewcommand{\vec}[1]{\boldsymbol{#1}}

\begin{document}

\maketitle

\begin{summary}
This paper focuses on the need for a rigorous theory of \emph{layered control architectures (LCAs)} for complex engineered and natural systems, such as power systems, communication networks, autonomous robotics, bacteria, and human sensorimotor control. All deliver extraordinary capabilities, but they lack a coherent theory of analysis and design, partly due to the diverse domains across which LCAs can be found.  In contrast, there is a core universal set of control concepts and theory that applies very broadly and accommodates necessary domain-specific specializations.  However, control methods are typically used only to design algorithms in components within a larger system designed by others, typically with minimal or no theory.  This points towards a need for natural but large extensions of robust performance from control to the full decision and control stack.  It is encouraging that the successes of extant architectures from bacteria to the Internet are due to strikingly universal mechanisms and design patterns.  This is largely due to convergent evolution by natural selection and not intelligent design, particularly when compared with the sophisticated design of components. Our aim here is to describe the universals of architecture and sketch tentative paths towards a useful design theory.
\end{summary}

\section{Introduction}\label{sec:intro}
Complex engineered and natural control systems, such as those used in robotics, the power grid, human sensorimotor control, and the internet, are characterized by needing to operate robustly and reliably across many spatiotemporal scales, despite being implemented using highly constrained hardware and software. Remarkably, a universal design pattern centered around \emph{layered control architectures (LCAs)} has emerged to address these challenges across vastly different domains.  These LCAs are the central object of study of this paper.

\begin{figure}
    \centering
    \includegraphics[width=\columnwidth]{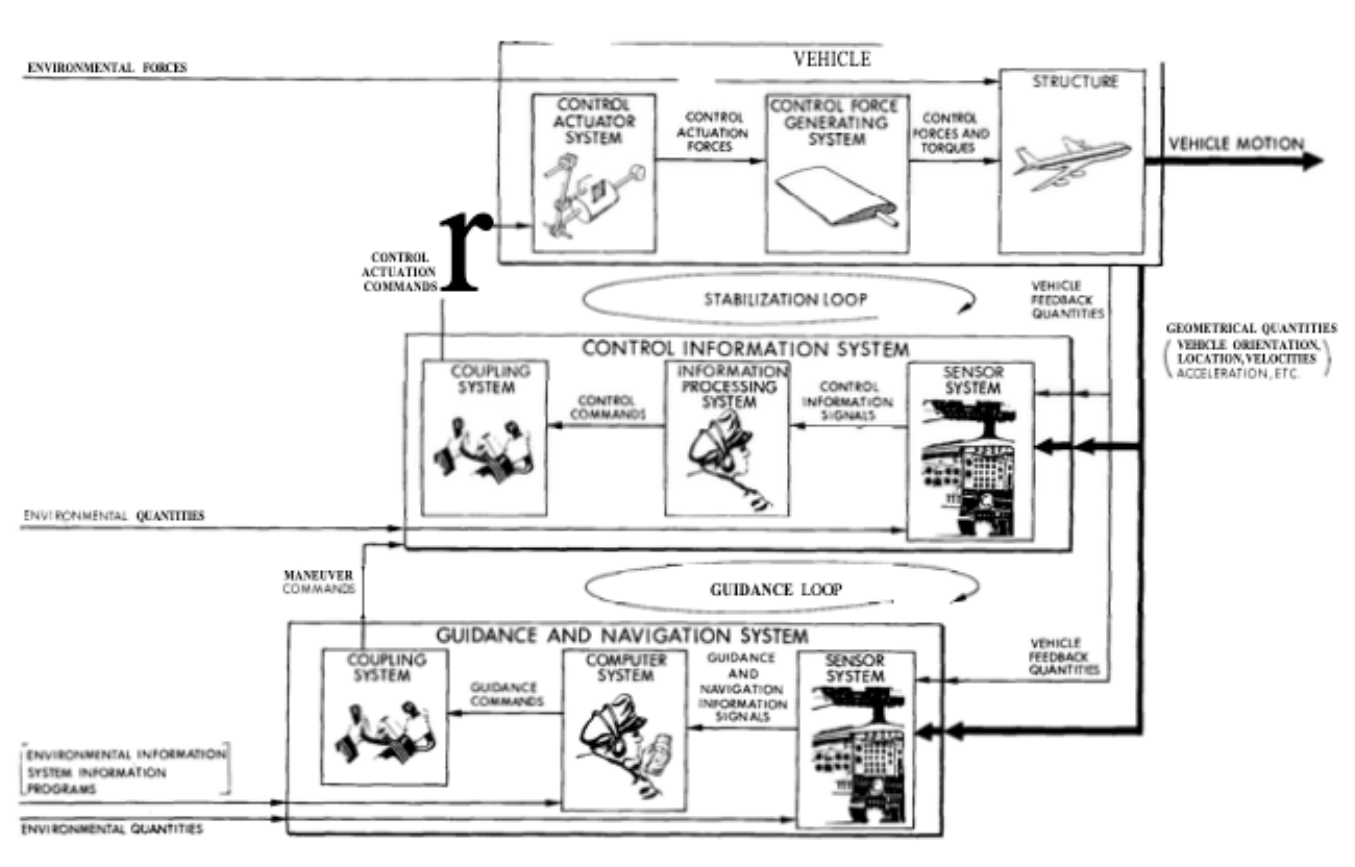}
    \caption{Figure taken from~\cite{apollo}, showing the GNC architecture used for Apollo missions.}
     \label{fig:apollo}
\end{figure}

Before proposing a broad definition of LCAs, we consider a familiar representative example from aerospace engineering, namely the widely used Guidance, Navigation, and Control (GNC) approach to aircraft control.  Here, the overall task of flying an aircraft from an initial location to a goal location is decomposed into tractable subproblems, as illustrated in Fig.~\ref{fig:apollo}, taken from the Apollo mission documentation~\cite{apollo}:
\begin{itemize}
    \item \textbf{Guidance}  determines a desired trajectory from the aircraft's current location to a goal location, in addition to nominal control actions, e.g., changes in forward and rotational velocity, for following the desired trajectory.
    \item \textbf{Navigation} is tasked with estimating the aircraft's state from onboard sensors, such as accelerometers and gyroscopes, and external signals, such as GPS.
    \item \textbf{Control} applies forces directly to the aircraft via actuators, e.g., steering, thrust, aileron deflection, in order to execute the trajectory planned by Guidance, all while maintaining aircraft stability.
\end{itemize}

We highlight some salient features of the GNC approach that we aim to capture in our broader theory of LCAs. The first, and most important aspect, is that an overall complex task (aircraft control) is decomposed into modular subtasks (Guidance, Navigation, and Control) of different complexity that operate at different frequencies over different spatiotemporal resolutions.  These control modules, or as we will call them, \emph{layers}, are allowed to interact, but only via well defined interfaces.  For example, the Control layer must operate at a high-frequency as it is tasked with stabilizing the unstable aircraft dynamics about a nominal trajectory in the face of an uncertain and dynamic environment, and hence is limited to simple feedback laws that can be implemented in real-time (e.g., PD control or LQR).  In contrast, the Guidance layer, which must contend with vast spatiotemporal scales in planning an aircraft's route, typically issues commands at a much slower frequency than the Control layer, as it must solve a longer horizon trajectory planning problem. Despite this modularization, the layers are nevertheless coupled via the exchange of a reference trajectory from Guidance to Control, and a tracking error from Control to Guidance.  Enabling both Guidance and Control is the Navigation layer, which is responsible for aircraft state estimation. 

\begin{figure}
    \centering
\includegraphics[width=0.95\columnwidth]{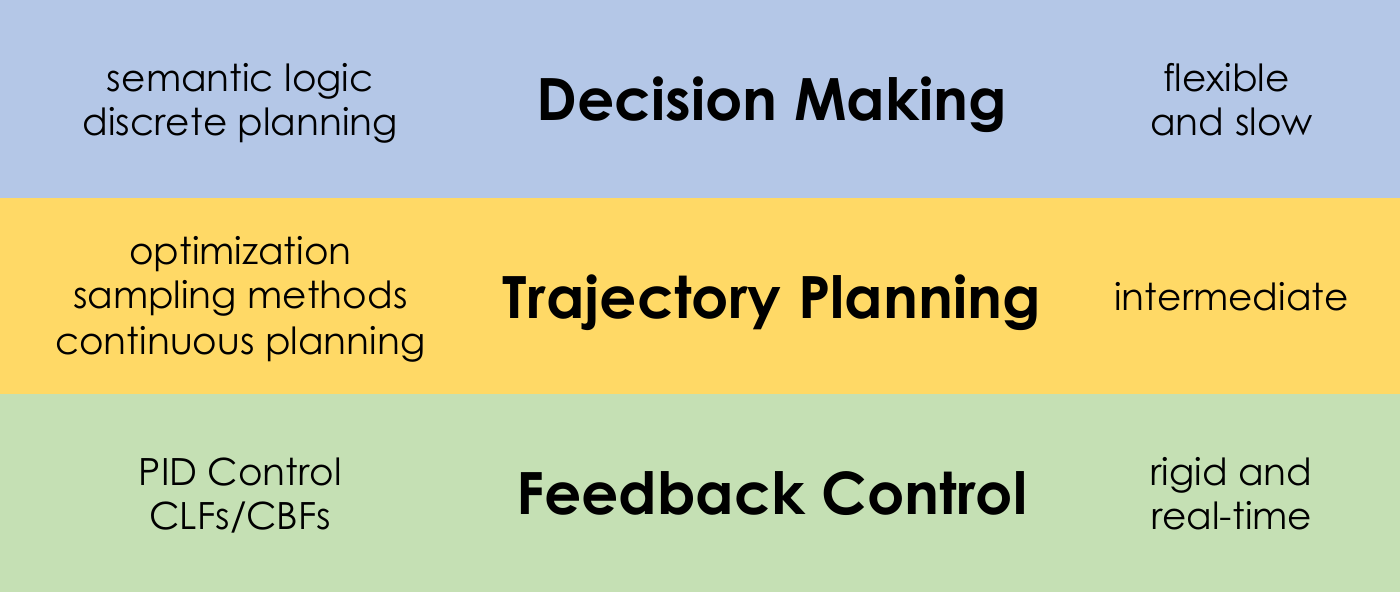}
    \caption{\footnotesize Layered control architectures are ubiquitous across natural and engineered systems.  We seek to initiate a quantitative study of layered control architectures based on the illustrated three-layer abstraction.}
    \label{fig:three_layers}
\end{figure}

\paragraph{A model LCA}
This paper seeks to initiate a \emph{quantitative study of LCAs} such as the one described above. To ground our discussion of LCAs, we begin with the ``model LCA'' shown in Fig.~\ref{fig:three_layers}, which is composed of three layers which broadly decompose across time-scales and complexity/flexibility:\footnote{By convention, we place slower more complex layers ``higher'' in the stack, and faster more rigid layers ``lower.''}
\begin{itemize}
    \item \emph{Decision making}: the top layer operates at the slowest frequency of the architecture, but is tasked with making complex logical decisions.  Primarily the domain of semantic logic and other discrete decision making tools, the decision making layer establishes mission objectives (e.g., which locations to visits via goal waypoints) and other system actions (e.g., delivering a payload, exchanging or collecting information).  In the context of GNC, this would be a higher ``mission layer'' that specifies the goal location.
    \item \emph{Trajectory planning}: the intermediate layer, sitting between decision making and feedback control, operates at a moderate frequency to generate trajectories that accomplish the mission objectives specified by the decision making layer.  Typical techniques employed at this layer include optimization-based (e.g., model predictive control, mixed integer programming) and sampling-based (e.g., rapidly-exploring random tree search) methods.  The generated trajectories, which are constrained to satisfy the mission objectives specified by the decision making layer, are transmitted to the feedback control layer.  This is precisely the Guidance layer in GNC.
    \item \emph{Feedback control}: the bottom layer operates at the fastest frequency of the architecture, and is tasked with tracking the trajectories generated by the planning layer.  This layer is the home of feedback control, and while offline computation to synthesize control gains may be sophisticated and expensive, online evaluation is typically constrained to be simple, fast, and rigid. In addition to ensuring that the system tracks the desired trajectory, feedback control also provides robustness to high-frequency and dynamic disturbance processes.  This is the Control layer in GNC.
\end{itemize}
\edit{A small note on terminology is in order before proceeding: although the word \emph{feedback} only appears in the bottom layer, it should be understood that some degree of feedback, either implicit or explicit, is present at all layers. For example, if trajectory planning is implemented using model predictive control, implicit feedback is provided by measuring the current system state.  As such, we ask the reader to interpret the use of the word feedback as indicating \emph{real-time explicit feedback control}, unless described otherwise.}

This architectural pattern or similar ones, which should be familiar to control theorists, appears consistently and broadly across domains despite both extreme diversity in the systems on which it is deployed, and the remarkable advances in sensing, actuation, and computation that have occurred over the past decades.  Despite the surprising unviersality of LCAs, they have yet to be a central object of study within the systems and controls community.  This paper is motivated by this current gap in the literature.  
 %

\edit{\paragraph{A brief overview of control architecture research}  While we defer more detailed literature reviews to appropriate sections, we pause to highlight that this manuscript builds upon and is inpsired by a rich literature, both academic and industrial, on process control and automation architecture.  Work providing a qualitative perspective about control architectures can be found in~\cite{bauer2017changes,samad2007system}.  While these works place a heavier emphasis on industrial applications and implementations, e.g., the use of Programmable Logic controllers to implement distributed control systems, they also touch upon topics core to this paper.  An interesting observation is that both papers acknowledge the importance of control architecture, while also recognizing its mercurial and difficult to define nature.  Although different terminology is used, a layered and multi-rate perspective is provided in both---indeed using model predictive control (MPC) for planning, and simple feedback control (PD control) for tracking is identified as a common design pattern, and is one that we revisit in great detail in the sequel.  We view these important qualitative perspectives as complementary to the frameworks we propose, and as further supporting the need for a more rigorous quantitative framework for reasoning about LCAs.  Domain-specific work centered around control architecture can be found in~\cite{samad2017controls} for smart-grid applications, in~\cite{lee2015cyber} for cyber-physical-system (CPS) applications, and in~\cite{chiang2007layering,palomar2006tutorial} for internet congestion control.  Once again, we see the key themes of this manuscript, such as layered multi-rate control implemented using diverse components, discussed.  For example, the templates proposed in~\cite{samad2017controls} can be directly mapped to the proposed layered strategies in \nameref{sec:layering}, ~\citet{lee2015cyber} propose a five layer architecture (called the 5C architecture), with each layer having different complexity and spatiotemporal scope, and initial quantitative methodologies for layering as optimization decomposition can be found in~\cite{chiang2007layering,palomar2006tutorial}.  This latter perspective serves as a key starting point for the framework proposed in this paper.  Finally, we note that although not the subject of this manuscript, an important enabling technology for control architecture design will inevitably be appropriate \emph{modeling languages and frameworks}, which may for instance be inspired by or build upon SysML~\cite{friedenthal2014practical}.}

\paragraph{Paper organization}  The rest of the paper is broadly organized into three parts.  In Part 1, composed of the next two sections, we first propose a framework for \emph{deriving} \nameref{sec:layering}.  We then provide concrete instantantiations of \nameref{sec:robotics} to illustrate the already impressive practical impact of layered control system design.  In Part 2, composed of the subsequent two sections, we propose an alternative perspective, and frame \nameref{sec:sweet-spots}.  A key takeaway of this section is that matching diversity across layers with diversity in control tasks can lead to LCAs that perform better than any individual layer could on its own.  We illustrate these concepts with \nameref{sec:neuro}.  Finally, in Part 3, composed of the penultimate section of the paper, we indulge in a more speculative discussion, and introduce qualitative definitions of what we believe to be other \nameref{sec:LLL}.  Finally, we end with \nameref{sec:conclusion}.




\section[Layered Control Architectures via Optimal Control Decomposition]{Part 1.1: Layered Control Architectures via Optimal Control Decomposition}
\label{sec:layering}
We propose a minimal quantitative framework for deriving and reasoning about LCAs such as those illustrated in Fig.~\ref{fig:three_layers}. Our starting point is a control policy synthesis problem which captures the key ingredients of modern complex systems that LCAs have evolved to address, namely: (i) the mix of discrete/logical decision making with continuous dynamics and control, and (ii) the diversity in time-scales at which different layers of a system (and its environment) evolve. Inspired by the Layering as Optimization Decomposition~\cite{chiang2007layering,palomar2006tutorial} approach to layered architectures, originally applied to network congestion control, our strategy is to systematically decompose the overall synthesis problem into tractable subproblems, each associated with a specific layer. 

We consider \emph{specifications} that the system must meet and \emph{safety constraints} that the system must obey.  We restrict ourselves to specifications expressed using formal logic, although alternative formulations are certainly possible.  These specifications describe system goals: for example, in robotic applications, such a goal might be navigate to a target, or to perform a household task.  Safety constraints, in contrast, are often expressed in terms of set membership constraints on a system's physical state: for example, in aerospace applications, such safety constraints may be expressed in terms of state/input inequality constraint enforcing the aerodynamic flight envelope.  Design problems also typically include auxiliary performance objectives (e.g., fuel efficiency, speed, robustness), which are optimized subject to the specification and safety constraints.  The goal then is to find the most efficient system design, as measured by the auxiliary performance objectives, that meets the system specifications and safety constraints.

Our framework thus centers around an \emph{overall synthesis problem} that seeks to find a control policy that satisfies high-level specifications, subject to system dynamics as well as state and input constraints.  In the interest of clarity, we omit auxiliary performance objectives in the initial formulation, but highlight natural ways in which they can be incorporated throughout.  
After defining the overall synthesis problem, we show that through suitable \emph{relaxations and decompositions}, (i) the three layer architecture described above can be derived, and (ii) familiar optimization-based trajectory planning and feedback control algorithms emerge naturally in an attempt to minimize the errors induced by these relaxations and decompositions.




\paragraph{Notation}
We use subscripts to denote continuous time, e.g., $x_t$ is the state $x$ at time $t\in\R$, and parenthesized numbers to denote discrete time, e.g., $x(k)$ is the state $x$ at discrete time step $k\in\mathbb{N}$. 
We use boldface font to denote infinite-horizon signals in both continuous and discrete time, i.e., $\vec x=(x_t)_{t\geq 0}$ or $\vec x = (x(0),x(1),\dots)$, depending on context.  We use the notation $x(T_1:T_2)$ to denote finite-horizon discrete-time signals, i.e., $x(T_1:T_2)=(x(T_1), x(T_1+1),\dots, x(T_2))$.  Finally, we use parentheses to concatenate two vectors, i.e., if $a\in\R^n$ and $b\in\R^m$, then $c=(a,b)\in\R^{n+m}$.

\subsection{The overall synthesis problem}
We assume that the system specifications are defined using a form of temporal logic, with Linear Temporal Logic (LTL) and Signal Temporal Logic (STL) being the most commonly used in the controls community---in the sequel, we use *TL to denote such general temporal logics.  In particular, we let the system goals be specified by a given *TL formula $\varphi$ defined over a finite set $\mathcal{AP}$ of atomic propositions.

The system state evolves according to the continuous-time dynamics
\begin{eqnarray}
\label{eqn:continuoustime}
\Sigma: \quad \dot x_t = f(x_t,u_t),
\end{eqnarray}
where $x_t \in \R^n$ is the system state, and $u_t \in \R^m$ is the control input.  In order to simplify exposition, we assume nominal dynamics without any model uncertainty or process noise.  Of course, real systems are subject to both, and how to systematically account for such uncertainty in LCAs remains an important open problem (see \nameref{sec:robust-lca}).  The state and control inputs are subject to safety constraints of the form $x_t \in \mathcal{X}\subseteq \R^n$ and $u \in \mathcal{U}\subseteq \R^m.$ 

We also define, for a suitable sampling time $\tau$, the corresponding discrete time model
\begin{eqnarray}
\label{eqn:discretetime}
\Sigma_d: \quad x(k+1) = f_d(x(k), u(k)),
\end{eqnarray}
where $x(k)=x_{k\tau}$, $u(k)=u_{k\tau}$, and $f_d$ is a discretization of the continuous time dynamics $f$.  Albeit somewhat cumbersome, we introduce both continuous and discrete time dynamics to highlight the \emph{multi-rate} nature of typical LCAs, wherein the decision making, trajectory planning, and feedback control layers all operate at different \emph{loop rates}, i.e., each layer recomputes or updates its action at a different frequency.  Where the switch from continuous to discrete time models is made in the LCA is often subject to computational constraints and loop rate requirements, which are in turn dictated by system specifications, safety constraints, and dynamics, see \nameref{sec:multirate} and \nameref{sec:cts_multilayer}.  Nevertheless, a common design pattern, which we adopt here, is to use continuous time models for real-time feedback control (to emphasize fast loop-rates), and discrete time models for both trajectory planning and decision making.


To verify the satisfaction of the *TL specifications $\varphi$, we assume the existence of a \emph{labeling} function $L:\mathcal{X} \to 2^{\mathcal{AP}}$ that associates a \emph{label} encoding the  \texttt{TRUE} atomic propositions at every state $x \in \mathcal{X}$.  
Finally, we define the trace, or run, of system $\Sigma$ under a control input sequence $\vec{u}$ to be the sequence:
$$
\xi(\vec x, \vec{u}):= (L(x(0)),u(0))((L(x(1)),u(1))\cdots,
$$
where here the signals $\vec x = (x(0), x(1), \dots)$, $\vec u = (u(0), u(1), \dots)$, are discrete time trajectories from the sampled system $\Sigma_d.$  If such a trace satisfies the specification $\varphi$, we write $\xi(\vec{x}, \vec{u})\models \varphi.$

We can now finally pose the overall synthesis problem as finding a possibly time-varying state-feedback policy $u_t:\mathcal X \to \mathcal U$ such that:
\begin{equation}\label{eq:master}
\begin{array}{rcl}
    \xi(\vec x, \vec u) &\models& \varphi, \\
    \dot x_t &=& f(x_t,u(x_t)), \ x_0 \text{ given},\\
    x_t &\in& \mathcal X, \,  u_t(x_t) \in \mathcal U\quad \forall t\geq 0.
\end{array}
\end{equation}

\begin{figure}[t]
    \centering
    \includegraphics[width=0.5\columnwidth]{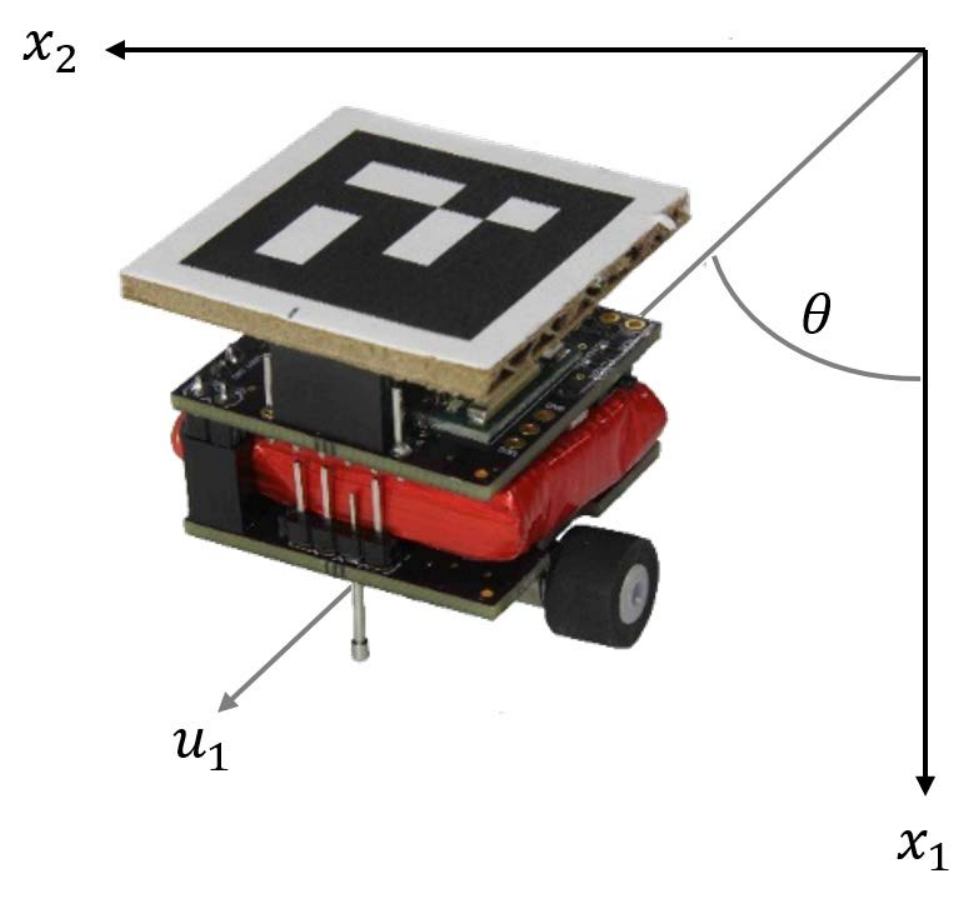}
    \caption{A mobile robot (the Gritbot) modeled as a Dubins' car.  In this case, the Gritbot that is deployed in the Robotarium which utilizes a LCA to allow for the implementation of user algorithms in a safe fashion \cite{pickem2017robotarium}. }
    \label{fig:gritbot}
\end{figure}

\begin{figure*}[t]
    \centering
    \includegraphics[width=.95\textwidth]{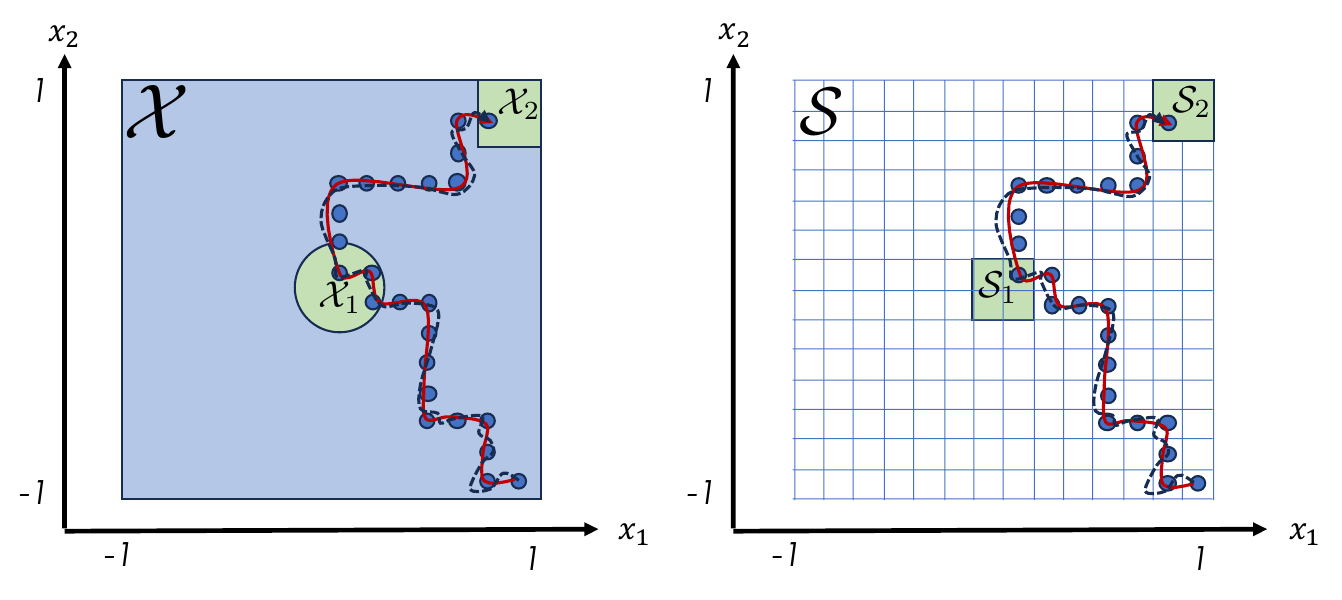}
    \caption{Moving from continuous (left) to discrete (right) state space in the robot navigation running example.  We see that due to the resolution of the partition, the set of discrete states $\mathcal{S}_1$ is an under-approximation of the corresponding first objective set $\mathcal{X}_1$.  The blue circles illustrate the discrete state space trace defined by the decision making layer; the red arrow shows the reference trajectory generated by the planning layer using the blue circles as waypoint constraints; the black dashed line is the actual system evolution, as driven by the feedback control layer on the continuous time dynamics.}
    \label{fig:gridworld}
\end{figure*}

\begin{my_exam}[Running example: robot navigation]
    We use a simple robot navigation problem as a running example to illustrate the concepts introduced in this section.  Suppose the system is the \emph{Dubins' car} (also called a \emph{unicycle}) with dynamics
    \begin{eqnarray}
    \label{eqn:unicycle}
    \begin{bmatrix}
        \dot{x}_1 \\
        \dot{x}_2 \\
        \dot{\theta} 
    \end{bmatrix} = 
    \begin{bmatrix}
        \cos\theta & 0 \\
        \sin\theta & 0 \\
        0 & 1  \\
    \end{bmatrix}
    \begin{bmatrix}
        u_1 \\
        u_2 
    \end{bmatrix}
    \end{eqnarray}
    and that the state and input constraints are simple box constraints
    \begin{align*}
        \mathcal{X}&=\{(x_1,x_2,\theta) \, |\, \max\{|x_1|,|x_2|\}\leq 1, \ |\theta|\leq \pi/2\},\\
        \mathcal{U}&=\{(u_1,u_2) \, | \, \max\{|u_1|,|u_2|\}\leq 1\}.
    \end{align*}
    To specify the system goals, we define the sets $\mathcal{X}_1=\{(x_1,x_2,\theta) \, | \, x_1^2 + x_2^2 \leq 0.1^2\}$ and $\mathcal{X}_2=\{(x_1,x_2,\theta) \, | \, x_1\geq 0.9, \,  x_2 \geq 0.9\}$.  The task specification is for the robot to first visit set $\mathcal X_1$, and then visit $\mathcal X_2$, i.e., to go to a small circle in the center of the space, and then go to the upper right corner (see Fig~\ref{fig:gridworld}).  This can be expressed in LTL via the specification 
    \begin{equation}
         \varphi = F(\mathcal{X}_1 \wedge F \mathcal{X}_2),
         \label{eq:gridworld-spec}
    \end{equation}
    where $\wedge$ and $F$ are the \texttt{and} and \texttt{eventually} atomic propositions, respectively.
\end{my_exam}


\subsection{Layering via Problem Decomposition and Relaxation}
Towards our goal of deriving a layered control architecture, we strategically rewrite the global problem \eqref{eq:master} by introducing redundant variables, and subsequently relaxing consistency constraints between these redundant variables to allow for modularization across layers.  \edit{We emphasize that while problem~\eqref{eq:master} is fundamental---in that it is dictated by the physics of the system as well as the specification and safety constraints of the problem---and that the proposed framework of problem decomposition and relaxation is foundational to a theory of LCAs, the particular realization of these ideas that follow are \emph{architectural design choices}.  They are by no means unique, although they are chosen to be broadly representative of approaches taken in the literature.}  


\paragraph{Decision making layer}
*TL specifications $\varphi$ are defined over a discrete state and input spaces, whereas the global problem~\eqref{eq:master} is defined over continuous state and input spaces. Towards bridging this gap, we assume that the continuous state space $\mathcal X$ admits a partitioning $\mathcal S$.  Ideally such a partition is such that for any cell $s \in \mathcal S$, we have that $L(x)=L(y)$ for all $x,y \in s$, i.e., all states in a partition satisfy the same atomic propositions and are ``semantically equivalent.''  However, if the partition of the state space is coarse, this may not hold true and partitionning may introduce conservatism. It is natural to consider this partition as defining a discrete state space for a Markov Decision Process (MDP), and with slight abuse of notation, we use $s\in\mathcal S$ to denote such a discrete state as well.  Similarly, a discrete action space $\mathcal{A}$ is induced by this partition and the system $\Sigma$, allowing us to define the MDP dynamics  $s(i+1) = f_{MDP}(s(i),a(i))$.  The MDP dynamics evolve in discrete-time, with dynamics defined to be consistent with traces of the sampled system $\Sigma_d$.  However, the MDP dynamics typically correspond to a \emph{sub-sampling} of $\Sigma_d$, i.e., $s(i)$ is determined by $x(k\delta)$, for $\delta\in\mathbb{N}_+$ a discrete sampling time.  In general, this sub-sampling may be irregular, e.g., if each time-step of the MDP is associated with a change in discrete state, but note that our model can capture this phenomena by including the null action $\emptyset$ such that $s(i+1) = f_{MDP}(s(i),\emptyset)=s(i)$ for all $s(i)\in\mathcal S$.

\begin{my_exam}[Running example: robot navigation]
A possible discrete state space $\mathcal{S}$ with $|\mathcal{S}|=14^2$ is shown in Fig.~\ref{fig:gridworld}(right), wherein each square corresponds to a discrete state $s\in\mathcal S$. This induces corresponding discrete actions $\mathcal{A}=\{\uparrow,\downarrow,\leftarrow,\rightarrow,\emptyset\}$ and MDP dynamics corresponding to those of a typical grid world problem (for clarity of exposition, we assume diagonal movement is not allowed).  In particular, letting $s=(s_1,s_2)$ denote the $(x_1,x_2)$ grid position, we then have
\begin{equation}\label{eq:gridworld-dyn}
(s_1(i+1),s_2(i+1)) = \begin{cases}
    (s_1(i),s_2(i)) & a(i) = \emptyset \\
    (s_1(i),s_2(i)+1) & a(i) = \uparrow \\
    (s_1(i),s_2(i)-1) & a(i) = \downarrow \\
    (s_1(i)+1,s_2(i)) & a(i) = \rightarrow \\
    (s_1(i)-1,s_2(i)) & a(i) = \leftarrow
\end{cases},
\end{equation}
where we only show dynamics for allowable actions $a$ that would not take the system outside of the state space $\mathcal{S}$.  Note that due to the resolution of state space discretization, the set $\mathcal{S}_1$ is an under-approximation of the corresponding continuous set $\mathcal{X}_1$.
\end{my_exam}

Given this definition and towards the goal of isolating a decision making layer we introduce redundant discrete planning variables $\vec s$ and $\vec a$, which are subject to the MDP dynamics $f_{MDP}$:

\begin{equation}\label{eq:master-redundant-0}
\begin{array}{l}
   \xi(\vec s, \vec{a}) \models \varphi, \ s(i+1) = f_{MDP}(s(i), a(i)), \ s(0)\ni x(0),\\
    x(k) \in s(\lfloor k/\delta\rfloor), \ x(k)=x_{k\tau}, \ \forall k\in\mathbb N,\\
    \dot x_t = f(x_t,u_t(x_t)), \ x_0 \text{ given},\\ 
    x_t \in \mathcal{X}, \ u_t(x_t) \in \mathcal U, \ \forall t\geq 0.\\
    
\end{array}
\end{equation}

Slightly abusing notation, here the first line replaces the trace over continuous variables $(\vec x, \vec u)$ with one defined over discrete variables $(\vec s, \vec a)$ satisfying the MDP dynamics, and we write $s(0)\ni x(0)$ to emphasize that the discrete state initial condition $s(0)$ must be consistent with the continuous state initial condition $x(0)=x_0.$   The second line enforces consistency between the remainder of the discrete plan of the first line and the underlying continuous control system $\Sigma$.  In particular, the constraint $x(k)\in s(\floor{k/\delta})$ ensures that the continuous state discrete time trace $(x(0),x(1),\dots)$ induces the correct discrete state trace $(s(0),s(1),\dots)$.  Similarly, the constraint $x(k)=x_{k\tau}$ ensures consistency between the continuous state discrete time trace $(x(0),x(1),\dots)$ and the continuous state trajectory $x_{t\geq 0}$.

Our first relaxation is to decompose problem~\eqref{eq:master-redundant-0} by isolating the discrete planning problem
\begin{equation}\label{eq:decision-layer}
    \xi(\vec s, \vec{a}) \models \varphi, \ s(i+1) = f_{MDP}(s(i), a(i)), \ s(0)\ni x(0),
\end{equation}
which we view as the decision making layer problem.  Here, the decision making layer assumes that any discrete plan $(\vec s, \vec a)$ can be realized by the underlying continuous control system $\Sigma$.  Solving the decision layer problem~\eqref{eq:decision-layer}, e.g., using *TL synthesis methods, yields a discrete state and action plan $(\vec s, \vec a)$ that satisfies the specification $\varphi$.  We show next how this high-level plan can be used to define a trajectory planning problem using a similar decomposition and relaxation technique.  We note that *TL synthesis methods are computationally expensive, and hence replanning at this layer is typically done at a slower time-scale, with the faster lower layers used to mitigate unexpected disturbances in the interim.

\begin{my_exam}[Running example: robot navigation]
    The decision making layer problem for this example becomes one of finding a state/action trace $(\vec s, \vec a)$ that satisfies the specification~\eqref{eq:gridworld-spec} subject to the discrete time discrete state dynamics~\eqref{eq:gridworld-dyn}. 
\end{my_exam}



\paragraph{Trajectory planning layer}

Removing the decision making layer problem~\eqref{eq:decision-layer} from problem~\eqref{eq:master-redundant-0} leaves us with
\begin{equation}\label{eq:master-redundant-1a}
\begin{array}{l}
    x(k) \in s(\lfloor k/\delta\rfloor), \ x(k)=x_{k\tau}, \ \forall k\in\mathbb N,\\
    \dot x_t = f(x_t,u_t(x_t)), \ x_0 \text{ given},\\ 
    x_t \in \mathcal{X}, \ u_t(x_t) \in \mathcal U, \ \forall t\geq 0.\\
\end{array}
\end{equation}
Towards the goal of isolating a trajectory planning layer, we introduce a redundant continuous state discrete time trajectory variable $\vec r=(r(0),r(1),\dots)$, constrained to be consistent with the discrete time state $\vec x = (x(0),x(1),\dots).$  The resulting strategically rewritten \emph{equivalent} problem is then given by
\begin{equation}\label{eq:master-redundant-1a}
\begin{array}{l}
    r(k) \in s(\lfloor k/\delta\rfloor), \ r(k)\in \mathcal{X},\\
    r(k) = x(k), \ x(k)=x_{k\tau},\forall k\in\mathbb N,\\
    \dot x_t = f(x_t,u_t(x_t)), \ x_0 \text{ given},\\ 
    x_t \in \mathcal{X}, \ u_t(x_t) \in \mathcal U, \ \forall t\geq 0.\\
\end{array}
\end{equation}

The first line isolates a trajectory generation problem, defined now over the reference trajectory variable $r(k)$, ensuring that (a) the trajectory is consistent with the discrete plan, as enforced by $r(k) \in s(\floor{k/\delta})$, and (b) the trajectory is safe, as enforced by $r(k)\in\mathcal{X}$.  As above, the second line enforces coupling between different layers: $r(k)=x(k)$ ensures that the trajectory and system states are consistent, and once again $x(k)=x_{k\tau}$ ensures consistency between the discrete time system $\Sigma_d$ and the continuous time system $\Sigma$.

Our next relaxation is to decompose problem~\eqref{eq:master-redundant-1a} by isolating a middle layer trajectory planning problem, which takes the form of the feasibility problem:
\begin{equation}\label{eq:mid-level}
\begin{array}{rl}
r(k) \in s(\lfloor k/\delta\rfloor), \ r(k) \in \mathcal{X}, \forall k\geq 0, \ r(0)=x(0),\\
\end{array}
\end{equation}
Here we drop all consistency constraints that $r(k) = x(k) = x_{k\tau}$ except for those enforcing the initial condition $r(0)=x(0)$, and a reference trajectory $\vec r$ that satisfies the state constraint $r(k) \in \mathcal{X}$ and specification constraints $r(k) \in s(\lfloor k/\delta\rfloor)$ is searched for. Due to this decomposition and relaxation of the constraints, the reference trajectory produced by solving feasibility problem~\eqref{eq:mid-level} is not guaranteed to be dynamically feasible, and therefore the tracking error between the true system state $x(k)$ and the planned trajectory $r(k)$ should be accounted for when enforcing the specification and safety constraints $r(k)\in \mathcal{X}\cap s(\lfloor k/\delta\rfloor)$.  Letting $\mathcal{C}(k):=\mathcal{X}\cap s(\lfloor k/\delta\rfloor)$ denote the intersection of the specification and safety constraints at discrete time step $k$, define the tightened constraint set $\overline{\mathcal{C}}(k)$, and replace the specification and safety constraints with $r(k)\in\overline{\mathcal{C}}(k)$.  How to appropriately tighten this constraint set depends on the feedback control layer, but once tracking error has been characterized, standard tools can be used.

In order to promote reference trajectories that are easy to track by the underlying continuous system $\Sigma$, the feasibility problem~\eqref{eq:mid-level} is often modified to produce approximately dynamically feasible solutions.  For example, it is common to decompose the control policy into a feedforward term $u_{ff}$, depending on the reference trajectory, and a feedback term $u_{fb}$, depending on the system state (or more specifically, on the tracking error). For example, a typical such decomposition is to simply set $u_t(x_t)= u_{ff}(r(\floor{t/\tau})) + u_{fb}(e_t)$, for tracking error $e_t:=x_t-r(\floor{t/\tau}))$.  Another standard approach is to assume that the reference trajectory obeys simplified planning dynamics $r(k+1) = f_{plan}(r(k), u_{ff}(k)).$  A common choice for these simplified planning dynamics is to use a \emph{reduced order} model defined by $y(k+1) = f_{rom}(y(k),v(k))$, where $y\in\R^p$ and $v\in\R^s$, with $p\leq n$ and $s\leq m$.  This case can be integrated into the proposed framework by replacing the consistency constraint $r(k) = x(k)$ with a reduced order consistency constraint $y(k) = \Projx(x(k))$, for $\Projx:\R^n\to\R^p$ some projection map encoding the model order reduction. To distinguish this important special case, in the sequel we reserve $r(k)$ for \emph{full order} reference trajectories, i.e., reference trajectories defined over the entire state with both $r(k),x(k)\in\R^n$, and use $y(k)$ to denote reduced order model states, as these are often used as \emph{tracked outputs} at the feedback control layer.

Finally, we make the following additional modifications.  First, towards employing a receding horizon control approach, we restrict the trajectory planning problem to be over a finite horizon $N$.  Second, we encode additional desirable properties of the trajectory, e.g., smoothness, via a running cost function $C(r,u)$ and a terminal cost $C_N(r)$. Integrating these elements with those described above yields the planning problem solved at discrete time step $k$:
\begin{equation}\label{eq:planning}
    \begin{array}{rl}
         \mathrm{minimize} & \sum_{i=k}^{k+N-1} C(r(i), u_{ff}(i)) + C_N(r(k+N)),  \\
         \text{subject to} & r(i+1) = f_{plan}(r(i), u_{ff}(i)), \ r(k)=x(k), \\
         & r(i) \in \overline{\mathcal{C}}(i),\
         i = k, k+1,\dots,k+N.
    \end{array}
\end{equation}
The planning problem~\eqref{eq:planning} is typically solved and implemented in a receding-horizon fashion, e.g., via model predictive control (MPC), and as such is limited to being resolved at an intermediate frequency (i.e., more frequently than the decision layer problem, but less frequently than the feedback control layer).

We end by noting that this is but one approach to defining a planning problem given a discrete state plan $\vec s$.  Alternative approaches consider, for example, loss functions that penalize deviations of the reference trajectory $\vec r$ from particular waypoints that are consistent with the discrete state plan $\vec s$.  It is hopefully clear that problem \eqref{eq:planning} is equivalent to such an approach, up to hard/soft constraints.

\begin{my_exam}[Running example: robot navigation]
 By converting the state trace $(s(0),s(1),\dots)$ computed by the decision making layer from discrete $(s_1,s_2)$ coordinates to continuous $(x_1,x_2)$ coordinates, e.g., by choosing the centroid of cell $(s_1,s_2)$, these can be used to define a sequence of \emph{waypoints} $((p_1(0),p_2(0)), (p_1(1),p_2(1)),\dots)$, with $(p_1(i),p_2(i))\in\R^2$, that can be used as \emph{state-constraints} within the planning layer.  These waypoints are illustrated with blue circles in Fig.~\ref{fig:gridworld}: on the left, we illustrate their use as waypoints for continuous trajectory planning, and on the right, we illustrate their use as a feasible state trace in the discretized state space $\mathcal S$.
 
We formulate the planning problem using a reduced order linear model composed of decoupled single-integrator dynamics in the $x_1$ and $x_2$ directions, i.e., we set $y(k) = (y_1(k), y_2(k))$, $v(k)=(v_1(k), v_2(k))$, and $f_{rom}(y(k),v(k))=(y_1(k) + \tau v_1(k), y_2(k) + \tau v_2(k)).$  In this case, we are using a reduced order model that projects out the angle $\theta$ and that introduces feedforward linear velocity inputs $(v_1,v_2)$.
We use the waypoints $((p_1(0),p_2(0)), (p_1(1),p_2(1)),\dots) $ to define constraints on the trajectory of the form $|y_1(k)-p_1(\floor{k/\delta})|\leq \Delta$ for $\Delta$ half the length of the square cells defining the discrete state space, and idem for the $x_2$-coordinate, i.e., we ask that the reference trajectory follow the sequence of discrete cells defined by the discrete state trace $(s_1(i),s_2(i))$, but allowing appropriate time within each cell as dictated by the sampling rate $\delta$ used at the decision making layer.  We additionally impose smoothness and control effort penalties in the objective, and constrain the reference trajectory to satisfy tightened state constraints (here we assume that the feedback control layer can guarantee a tracking error of at most .05 in either of the $(x_1,x_2)$ coordinates).  The resulting planning problem solved at time step $k$ over a horizon $N$ is then given by:
\begin{equation}
    \begin{array}{rl}
         \mathrm{minimize} & \sum_{i=k}^{k+N-1} \|y(i+1)-y(i)\|_2^2 + \|v(i)\|_2^2  \\
         \text{subject to} & y(i+1) = f_{rom}(y(i), v(i)), \\ & y(k)=(x_1(k),x_2(k)), \\
         & |y_1(i)-p_1(\floor{i/\delta})|\leq \Delta-0.05, \\
         & |y_2(i)-p_2(\floor{i/\delta})|\leq \Delta-0.05, \\ 
         & y(i) \in \mathcal{\bar X}=\{(x_1,x_2) \, | \, \max_j|x_j|\leq 0.95\}, \\
         & i=k,k+1,\dots,k+N.
    \end{array}
\end{equation}
An illustrative example of the resulting reference trajectory is shown in red in Fig.~\ref{fig:gridworld}(left).
\end{my_exam}

\paragraph{Real-time feedback control layer}
Finally, we consider the feedback control layer.  At discrete time step $k$, given a solution $(r(k:k+N),u_{ff}(k:k+N))$ to the planning layer problem~\eqref{eq:planning}, we must contend with the remaining constraints, now truncated to the planning horizon $N$:
\begin{equation}
   \begin{array}{l}
        r(i) = x(i), \ x(i)=x_{i\tau}, \ i= k, k+1, \dots, K,\\
    \dot x_t = f(x_t,u_t(x_t)), \ \text{$x_0$ given,}\\ 
    x_t \in \mathcal{X}, \ u_t(x_t) \in \mathcal U, \ \forall t\in[k\tau,(k+N)\tau].\\
   \end{array}
\end{equation}

We relax this problem by (a) removing the state constraint $x_t \in \mathcal{X}$, as this is addressed within the planning problem~\eqref{eq:planning}, and (b) allowing for the state $\vec x$ to deviate from the reference trajectory $\vec r$, as the reference trajectory $\vec r$ is not expected to be dynamically feasible:
\begin{equation}\label{eq:control0}
    \begin{array}{rl}
         \mathrm{minimize} & \int_{k\tau}^{(k+N)\tau}\left(\|e_s\|_Q^2 + \|u_{fb,s}\|_R^2\right)ds   \\
         \text{subject to:} & \dot x_t = f(x_t,u_{ff,t} + u_{fb,t}), \ \text{$x_0$ given,} \\ 
         & u_{ff,t} + u_{fb,t} \in \mathcal{U}.
    \end{array}
\end{equation}
In the above, for a positive semi-definite matrix $P$, we let $\|z\|_P^2:= z^TPz$, define the tracking error $e_t := x_t - r(\lfloor t/\tau\rfloor)$, and slightly abuse notation by letting $u_{ff,t} = u_{ff}(\lfloor t/\tau \rfloor)$ be the zero-order holds of their corresponding discrete time signal.  This problem can be viewed as a ``best-effort'' tracking controller. 

Loop rate constraints lead to either offline computed feedback policies that approximately solve the problem (e.g., LQR tracking), to simple to implement approaches such as PD control, or to myopic simplifications that can be solved in real-time via e.g., quadratic programming.  All of these approaches can be viewed as further relaxations of the above tracking problem~\eqref{eq:control0}.  While widely used, the proposed relaxation~\eqref{eq:control0} and its extensions typically lack tracking error guarantees.  To address this concern, recent work has leveraged Lyapunov-based techniques to certify tracking error bounds, which in turn allow for a principled tightening of the constraints used in the planning layer.  We highlight some of these techniques, as applied to robotic LCAs, in the next section.

\begin{my_exam}[Running example: robot navigation]
\label{ex:exampleflat}
While many feedback control approaches are possible here, we take this opportunity to briefly introduce differential flatness-based control and show how it can be used in this context.  To synthesize a feedback controller that tracks the reduced order model reference trajectory $(y_t, v_t)=(y(\floor{t/\tau}),v(\floor{t/\tau}))$, we first identify the \emph{flat outputs} \cite{van1998real} such that the state and inputs can be written as a function of these flat outputs and derivatives.  For the unicycle dynamics, a flat output is $\xi = (x_1,x_2)$, as verified by the relationship
\begin{align}
\label{eqn:flatness}
\begin{bmatrix} x_1 \\ x_2 \\ \theta \end{bmatrix} &= \begin{bmatrix}
\xi_1 \\
\xi_2 \\
\arctan(\frac{\dot\xi_2}{\dot \xi_1})
\end{bmatrix}=:x_{\flat}(\xi,\dot\xi), \\ 
\begin{bmatrix} u_1 \\ u_2 \end{bmatrix} &= 
\begin{bmatrix}
    \sqrt{\dot \xi_1^2 + \dot \xi_2^2} \\
    \frac{\dot \xi_1 \ddot \xi_2 - \dot \xi_2 \ddot \xi_1 }{\dot \xi_1^2 + \dot \xi_2^2}
\end{bmatrix}=:u_{\flat}(\xi,\dot\xi,\ddot\xi).
\end{align}

In order to synthesize a feedback controller, it is convenient to define the flat state $z := (\xi_1,\xi_2,\dot\xi_1,\dot{\xi}_2) = (x_1,x_2,\dot x_1, \dot x_2)\in \R^4$ and flat control input $a:=(\ddot\xi_1,\ddot{\xi}_2)=(\ddot x_1,\ddot x_2)\in\R^2$, 
resulting in the linear dynamics
\begin{eqnarray}
\label{eqn:unicyclelinear}
\dot{z}_t = \underbrace{
\begin{bmatrix} 
0 & 0 & 1 & 0 \\
0 & 0 & 0 & 1 \\
0 & 0 & 0 & 0 \\
0 & 0 & 0 & 0 \\
\end{bmatrix}}_{A_{rom}^c} z_t + 
\underbrace{ 
\begin{bmatrix}
0 & 0 \\
0 & 0 \\
1 & 0 \\
0 & 1 \\
\end{bmatrix}}_{B_{rom}^c} a_t
\end{eqnarray}
for $z_t \in \R^4$ and $a_t \in \R^2$.  With slight abuse of notation, we may then write $x_t = x_{\flat}(z_t)$ and $u_t=u_{\flat}(z_t,a_t)$ by making the appropriate identifications between $(\xi_t,\dot \xi_t,\ddot \xi_t)$ and $(z_t, a_t)$.  
We can then synthesize a feedforward plus feedback policy using, for example, PD control by setting
\begin{equation}\label{eqn:flat_feedback}
    a_t = \dot v_t - K^{\flat}_P((z_1,z_2)_t-(y_1,y_2)_t) - K^{\flat}_D((z_3,z_4)_t-v_t),
\end{equation}
for appropriately tuned positive definite matrices $K^{\flat}_p$ and $K^{\flat}_D$.  Note that in this case, the feedforward term $\dot v_t$ can be approximated via a finite difference, i.e., ${\dot v_t \approx (v(\floor{t/\tau})-v(\floor{t/\tau}-1))/\tau}.$

%
%
%
%

We now discuss how to translate the flat state $z_t$ and flat input $a_t$ to a control input $u_t$ composed of feedforward and feedback terms given current measurements of $z_t=(x_{1},x_{2},\dot x_{1},\dot x_{2})_t$.  We define the feedforward term to be given precisely by the mapping from flat state and input to control input:
$$u_{ff,t}(z_t,a_t)=u_{\flat}(z_t,a_t).$$  
If the mapping $u_{\flat}$ exactly captures the system dynamics, then no additional feedback term would be required.  However, in practice, the Dubins' car is often used as a reduced order model for planning trajectories for more complex systems such as quadrupeds (see Fig.~\ref{fig:multi-rate} in \nameref{sidebar:multirate} and Ex.~\ref{ex:mpcdubins}).  As such, a feedback term to ensure that flat and actual states match is additionally required.  One such option is again a PD controller:
\begin{multline*}
    u_{fb,t}(z_t)=-K_P ((x_1,x_2)_t-(z_1,z_2)_{t+\tau}) - \\ K_D((\dot{x}_1,\dot{x}_2)_t-(z_3,z_4)_{t+\tau}),
\end{multline*}
for positive definite matrices $K_P$ and $K_D$. Note here that the flat lookahead state $z_{t+\tau}$ is obtained by forward integrating the flat dynamics~\eqref{eqn:unicyclelinear} with the $a_t$ given as in~\eqref{eqn:flat_feedback}  and $z_t=(x_1,x_2,\dot{x}_1,\dot{x}_2)_t$ obtained from hardware measurements. The final control input is then the sum of both the feedforward and feedback terms, i.e.,
\begin{equation}\label{eq:dubins-feedback}
u_t(z_t,a_t)=u_{ff,t}(z_t,a_t)+u_{fb,t}(z_t),
\end{equation}
as suggested in the feedback control problem~\eqref{eq:control0}.
A conceptual illustration of the resulting evolution of the actual system state is shown with a dashed black line in Fig.~\ref{fig:gridworld}.
\end{my_exam}


\subsection{Discussion}
We obtained a LCA by introducing redundant variables and suitable relaxations to decompose the overall synthesis problem~\eqref{eq:master} into tractable decision making, trajectory generation, and feedback control subproblems, as specified in equations~\eqref{eq:decision-layer}, \eqref{eq:planning}, and~\eqref{eq:control0}, respectively. This highlights another key feature of LCAs: they allow for intractable global problems to be decomposed into tractable subproblems, often with minimal loss in performance or efficiency. We note that the proposed decomposition is only one of certainly many approaches, and is chosen to be consistent with the rest of the manuscript. Indeed, while the above framework provides a more formal perspective on LCAs, it still has an element of ``art'' to it.  In particular, \emph{how to relax the overall synthesis problem}, as well as \emph{how to bridge the simplifications between the different layers} of the resulting architecture is still up to the designer.  Nevertheless, by posing the problem in this way, there is a natural nested optimization structure that emerges that may allow for more principled methods of LCA design to be defined.  We highlight next some key open questions and concepts that we do not treat in depth, but certainly deserve more investigation.

\edit{\paragraph{What about the hardware?}
Each layer described above delineates a \emph{functional component} of an overall decision and control stack.  Equally important are the \emph{physical substrates} used to implement these functional layers. 
 For example, a motor used in a robotic system to actuate a joint is composed of different scales of components, ranging from circuit elements to microprocessors to motor components. While these physical substrates are closely related to the layers they are used to implement, they are distinct: to make this explicit we use the term \emph{levels} for physical substrates, and reserve layers for functional components.  We expand on this idea, and introduce other key concepts of LCAs not touched upon here, in \nameref{sec:LLL}.  Furthermore, in~\nameref{sec:sweet-spots}, we present a quantitative framework to inform how to choose diverse hardware to implement diverse functionality as a function of diversity in the control task at hand, and instantiate this perspective in~\nameref{sec:neuro}.}

\edit{\paragraph{How many layers should there be?}
This section presented an approach to deriving an LCA with three layers, each operating at different spatiotemporal resolutions.  While these three layers, namely decision making, trajectory generation, and feedback control, are commonly found in complex engineered systems, this pattern is by no means the only one possible.  Indeed, all of the concepts introduced above can be applied recursively, leading to layers of layers.  For example, there can be several layers of trajectory planning, operating at different loop rates, using different planning models, and planning over different horizons, see for example \nameref{sec:cts_multilayer}.  Similarly, nested control loops are a standard control design pattern that can be interpreted as different layers of real-time feedback control.  While we present a framework for deriving layers given an overall problem formulation, we still lack quantitative tools for deciding how many layers there should be, as well as what information should be exchanged between them.  This is undoubtedly a key open question.}

\paragraph{Multi-rate control}
In the above, we hint at the role of multi-rate control in LCAs that seeks to address implicit timing constraints.  Low layer feedback control~\eqref{eq:control0} operates in (near) real-time, trajectory generation~\eqref{eq:planning} at a slower rate, and decision making~\eqref{eq:decision-layer} at a slower rate still.  This suggests that ideas from singular perturbation analysis and time-scale separation (see for example~\cite{kokotovic1999singular,zhang2014singular} and references therein) may also be used to provide further rigor to the approach.  We further explore \nameref{sec:multirate} in robotic systems in the next section.

\edit{
\paragraph{Robust LCAs}\label{sec:robust-lca}
We omit uncertainty due to process noise and modeling errors in our development.  However, practical LCAs must be robust to these effects.  Promising approaches to tackling uncertainty between layers include the use of robust Lyapunov certificates for guaranteeing bounded tracking error of a reference trajectory by the feedback control layer (see for example \nameref{sec:cts_multilayer}), and more generally, the use of assume-guarantee contracts~\cite{alur1999reactive, chen2020safety}.  These approaches are intuitive and effective, but it is nevertheless of interest to investigate whether such certificates can be derived by applying similar decompositions and relaxations to a robust overall synthesis problem that explicitly acknowledges uncertainty in its initial formulation.}

\edit{
\paragraph{Layered sensing architectures}
We also emphasize that although our focus in this section has been on fully-observed state-feedback control problems, analogous layered decompositions for sensing and output feedback problems need to be developed.  A promising starting point is to recognize that different sensors, ranging from semantically rich and complex sensors (e.g., cameras and LIDAR) to simple single output sensors (e.g., IMUs and gyroscopes), are naturally assigned to each of the decision making (e.g., computer vision, semantic segmentation), trajectory generation (e.g., VIO + SLAM), and feedback (e.g., IMUs) layers. }

\paragraph{Learning in LCAs}
The use of rich perceptual sensors such as cameras invariably introduces learning into the resulting LCAs, which is a topic we cannot hope to do justice to within the scope of this paper.  This is however an exciting and important direction to be explored, with learning and data-driven techniques poised to make significant impact in designing effective LCAs.

\subsection{Related work} 
The framework proposed above is inspired by a rich literature seeking to establish principles of LCA design, although it is not necessarily explicitly identified as such.

\paragraph{Layering as optimization decomposition}
The layering as optimization decomposition (cf. \cite{chiang2007layering,palomar2006tutorial} and references therein) and the reverse/forward engineering (cf. \cite{zhao2014design,cai2017distributed} and references therein) paradigms have been particularly fruitful in tackling internet and power-grid control problems, respectively.  Both of these frameworks can be loosely viewed as using the dynamics of the system to implement a distributed optimization algorithm through vertical (layering) and horizontal (distributed) decomposition.  These methods ensure that the state of the system converges to a set-point that optimizes a utility function.  These approaches can scale to large systems by taking advantage of the structure underlying the utility optimization problem, and can simultaneously identify and guarantee stability around an optimal equilibrium point.  Nevertheless, they do not explicitly consider optimal control, and in particular transients, in their analysis, making them an important but incomplete first step towards a theory of LCAs.

\paragraph{Decision making and continuous control}
This line of work seeks to make explicit that although formal specifications are inherently discrete, in order to ensure that a system satisfies them, designers must account for continuous dynamics and control.  One line of work seeks to reformulate *TL specifications into continuous control tasks through the use of control Lyapunov functions (CLFs) and control barrier functions (CBFs)~\cite{ames2014control,ames2016control,ames19cbf}, see for example \cite{lindemann2018control,dimitrova2014deductive}.  An alternative approach is to encode *TL constraints via mixed integer linear constraints in robust/optimal control problems \cite{raman2014model,sadraddini2015robust,wolff2014optimization}, or to abstract the continuous control problem into an uncertain finite-state MDP and use robust dynamic programming \cite{wolff2012robust}.  Closely related is the work of \citet{fan2020fast}, wherein decision making is done via SAT-based trajectory planning methods which solve a satisfiability (SAT)
problem over quantifier free linear real arithmetic.  Other representative works that explicitly acknowledge the inherently hybrid (discrete/continuous) nature of the decision making and control problem, and that seeks to bridge them in a principled way, include reactive planning approaches~\cite{kress2009temporal} and the use of CBFs for determining the magnitude of disturbance that a system can be subject to while still ensure satisfaction of STL specifications~\cite{akella2021disturbance}.
Barrier functions have also been applied to POMDPs in the context of distribution temporal logic (DTL) \cite{ahmadi2020barrier} and coherent risk measures (e.g., CVaR) \cite{ahmadi2021risk,singletary2022safe} to enforce safety constraints at a planning level.  More broadly, risk-aware planning and control is considered in~\cite{majumdar2020should,singh2018framework,hakobyan2019risk}.  Implicit in all of the above is a layered architecture wherein the high-layer decision making component operates on a discrete abstraction of the underlying continuous time system, and similarly, the underlying continuous time system implements planning/control layers in order to meet the plan specified by the top decision making layer.

\paragraph{Trajectory generation and continuous control}
Approaches to dynamics aware trajectory generation typically follow a classic two phase approach, wherein first a graph is constructed whose nodes are collision free configurations and whose edges correspond to feasible paths between these configurations, see for example \cite{kavraki1996probabilistic} and references therein.  More recent approaches based on graphs of convex sets~\cite{marcucci2022motion}, motion primitives~\cite{ubellacker2022robust,majumdar2017funnel}, optimization-based methods~\cite{yin2020optimization}, control lyapunov, barrier, and contraction metrics~\cite{singletary2020safety, csomay2022multi, singh2017robust, singh2021safe, sun2021learning, cohen2020approximate}, and reachability techniques~\cite{herbert2017fastrack} have also been proposed to bridge the gap between low-layer fast-time scale control and middle-layer intermediate-time scale trajectory generation.  The common theme in all of these approaches is the use of simplified dynamics in the trajectory generation layer, allowing for fast replanning, and a low-layer feedback controller that provides certifiable guarantees on tracking error.  Finally, most closely related to the framework presented in the previous discussion are the results found in~\cite{matni2016theory,srikanthan2023augmented,srikanthan2023data,zhang2024change}, wherein it is shown that a two-layer trajectory generation/feedback control LCA can be obtained by suitably relaxing consistency constraints between the state and reference trajectory.  A key feature of this approach is that the trajectory planning problem is augmented with a tracking penalty regularizer that promotes dynamic feasibility of the synthesized reference signal.

\section[Layered Control Architectures for Robotic Systems]{Part 1.2: Layered Control Architectures for Robotic Systems}
\label{sec:robotics}

\begin{figure*}[t!]
\centering
\includegraphics[width=0.88\textwidth]{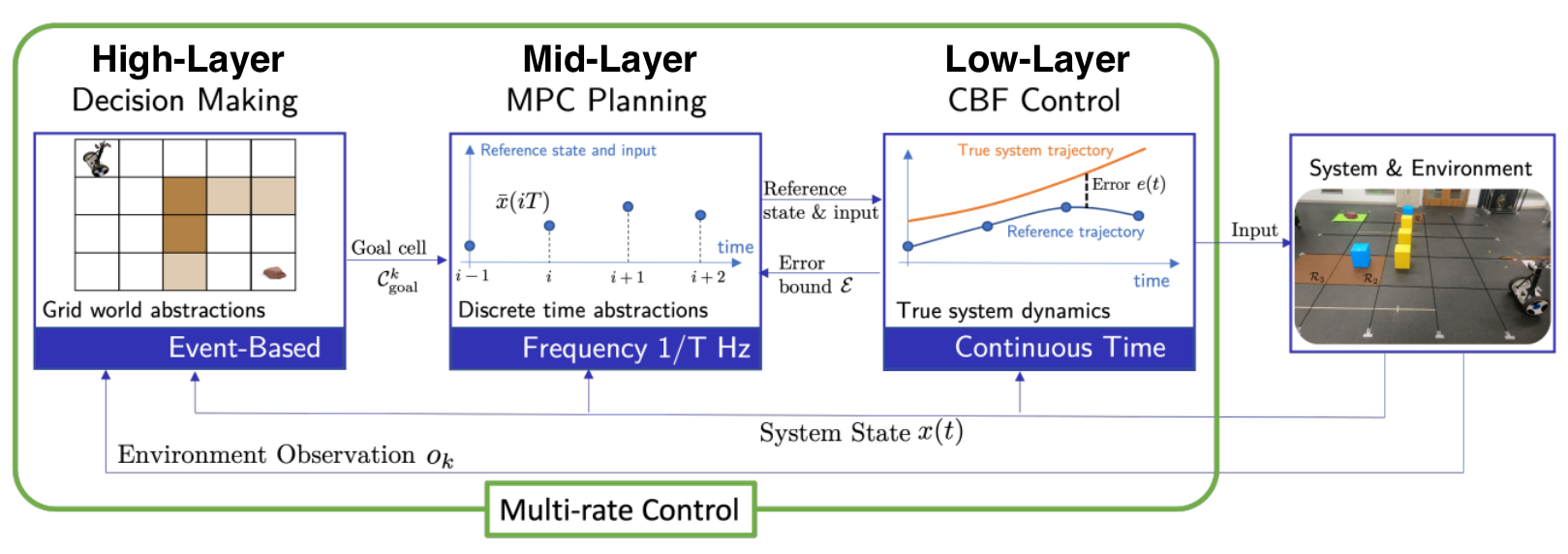}
\caption{Multi-rate robotic control can be viewed through the lens of layered architectures (figure from \cite{rosolia2022unified}).}
\label{fig:multi-layered}
\end{figure*}

LCAs have long found use on robotic systems---empirically, it is well-known that this is the best (and arguably only) way to implement controllers in practice.  Yet, despite this empirical evidence there is very little analysis of LCAs.  Conversely, while the control community applies rigorous approaches to controller synthesis, it is often only applied to a single layer.  This points to a unique opportunity for the controls community: reverse-engineering and analyzing the LCAs deployed on robotic systems that have proven useful in practice.  

To put the central role of LCAs on robotic systems in context, one should first consider the hardware itself.  A robotic system, broadly defined, typically consists of three main components: sensors used for perception, a central processor, and motor controllers used for driving actuators.  Concrete examples of this include: cameras mounted on a legged robot for localization and mapping, or proximity sensors on a vehicle for advanced driver assistance.   In this context, LCAs (as shown in Fig. \ref{fig:multi-layered}) are often deployed relative to these physical levels on hardware (see \nameref{sec:LLL} for more examples of how levels and layers interact in LCAs).  Perception leads to a decision making layer operating at a discrete/semantic level of abstraction, the central processor leads to reference signal generation using reduced order models, and finally at the actuator level real-time algorithms instantiate feedback control.\footnote{The robotics literature refers to these different layers as high-level, mid-level, and low-level control.  We however argue that these are better viewed through the lens of LCAs, and hence use the layered terminology defined in the previous sections.}  


This section gives concrete instantiations of LCAs for robotic systems.  We start by defining robotic system dynamics, and subsequently work our way up the layers of a typical instantiantion of LCAs for robotic control.  A goal of this section is to highlight the importance of multi-rate control in the context of LCAs. We start with the feedback control layer---termed the \emph{real-time feedback control layer} herein to highlight the fast loop rate at which it is implemented.  We then discuss trajectory planning paradigms, and end by highlighting approaches to their integration with the real-time feedback control layer.  We forego a discussion on the decision making layer for the sake of brevity.  See \cite{rosolia2022unified} for a formal inclusion of decision making with the methods presented in this section. 

\subsection{Robot Dynamics}

Robotic systems are inherently governed by nonlinear equations of motion.  These represent the physical evolution of the system and are typically obtained from the Euler-Lagrange equations: 
\begin{equation}
\label{eqn:eom}
D(q_t) \ddot{q_t} + C(q_t,\dot{q}_t) \dot{q}_t + G(q_t) = B u_t, 
\end{equation}
where here $q_t \in Q$ are the configuration variables of the system, $\dot{q}_t \in T_q Q$ is a vector of velocities (which take values in the tangent space to the configuration space), $D(q_T)$ is the inertia matrix, $C(q_t,\dot{q}_t)$ the Coriolis matrix, $G(q_t)$ contains the gravity related terms, and $B$ is the actuation matrix.  Defining the state vector: $x_t = (q_t, \dot{q}_t) \in TQ$, where for simplicity we can work with a local coordinate chart of $Q$ wherein $TQ \cong \R^n$ for $n$ even, allows for the formulation of a control system affine in the control input: 
\begin{eqnarray}
\label{eqn:affinecontrolsys}
\dot{x}_t = f(x_t) + g(x_t) u_t . 
\end{eqnarray}
where $f : \R^n \to \R^n$ and $g : \R^n \to \R^n \times \R^m$ can be directly obtained from \eqref{eqn:eom}. 
The end result is a control system of the form \eqref{eqn:continuoustime} with the additional structural property that the control input appears in an affine fashion, an observation which has important ramifications for controller synthesis. 

\subsection{Real-Time Feedback Control Layer}
We begin by tackling the feedback control problem~\eqref{eq:control0} for robotic systems.  We assume that a trajectory, containing both a (reduced order) reference trajectory and a feedforward control term, is available.  We discuss approaches to solving this trajectory problem after addressing the feedback control problem.  An emphasis is placed throughout on the need for \emph{real-time} feedback control.

\paragraph{Linear Control}  
For robotic systems, the most common form of real-time controller\footnote{Historically called ``low-level control.''} is a simple linear feedback controller acting on the error, e.g., a PD controller implemented at the motor control layer.  These controllers are \emph{highly} effective in practice when implemented properly.  It is important to stress that this is not due to the dynamics being linear (they are not), nor does it imply that the dynamics are even locally linear (again, they are not).  Rather, these controllers work well exactly because of their use within a LCA, running at a fast loop rate and actuating as a function of tracking error.  This second observation, that linear control actuates only on the error (and thus is independent of model information) is critical.  
To provide a precise instantiation of 
real-time linear controllers, consider a continuous time reference signal $\vec r = (\vec r_{q},\dot{\vec r}_{q})$ that we wish to track (decomposed into a reference position and velocity) and an associated feedforward control input ${\vec u_{ff}}$.  Then the simplest form of feedforward and feedback control becomes: 
\begin{eqnarray}
\label{eqn:PDcontrol}
u_t &= u_{ff,t} + 
\ddot{r}_{q,t} - K_P (q_t - r_{q,t}) - K_D (\dot{q}_t - \dot{{r}}_{q,t}) \\
 &=: u_{ff,t} + \ddot{r}_{q,t} + u_{fb,t}(q_t - r_{q,t}, \dot{q}_t - \dot{{r}}_{q,t}),
\end{eqnarray}
for $K_P$, $K_D$ positive definite matrices (see for example the Dubins' car control policy \eqref{eq:dubins-feedback}).  Given the decoupled nature of this controller, it can be deployed in a decentralized fashion, i.e., actuator by actuator, on the motor controllers at a very fast loop rate (faster than $1 kHz$).  It is important to note that many variations of this controller are possible.  For example, an integral term can be added, or a reference velocity signal from a higher layer can be tracked in which case the $K_P$ term might be removed. 


That linear controllers can stabilize the nonlinear dynamics associated with a robotic system can be made rigorous in certain cases.  To this end, assume that the robotic system is fully actuated, i.e., $B$ is invertible, and for simplicity take $B = I$.  Then Picking $u _{ff,t} = G(q_t)$ in equation \eqref{eqn:PDcontrol} results in asymptotic stability of the tracking of the reference signal $\vec r$.  To see this, let $\vec e_{q} = \vec q - \vec{r}_{q}$ and $\dot{\vec e}_{q} = \dot{\vec q} - \dot{\vec r}_{q}$ be the position and velocity error in tracking the reference signal $\vec r$.  Define the error signal $\vec{e} := (\vec e_{q},\dot{\vec e}_{q})$ and consider the Lyapunov function candidate: 
$$
V({e}) =  \frac{1}{2} \dot{{e}}^T_q D({e}_q + {r}_q) \dot{{e}}_q + \frac{1}{2} {e}^T_q K_P {e}_q . 
$$
which is positive definite since $D(q)$ is symmetric positive definite.
Then differentiating $V$ along solutions of \eqref{eqn:eom} yields: 
$$
\dot{V}({e}) \leq - \dot{{e}}_q^T K_D \dot{{e}}_q \leq 0
$$
as $\dot{D}(q,\dot{q}) - 2 C(q,\dot{q})$ is skew symmetric. 
Invoking LaSalle's Principle then shows asymptotic stability of $({e}_q,\dot{{e}}_q) = (0,0)$, i.e., shows that the reference signal is asymptotically tracked. 


\paragraph{Nonlinear Control}
The use of Lyapunov functions in certifying linear controllers' ability to asymptotically track reference signals points to nonlinear controllers, based on Lyapunov functions, that can achieve improved performance.  Indeed, to maximize performance on robotic systems it is necessary to exploit the full nonlinear dynamics of the system, which can only be done with nonlinear controllers.  These controllers must however be synthesized in a way that yields both theoretical guarantees while also being deployable in practice, i.e., they must be implementable at fast loop rates ($>1 KHz$).   With this in mind, a key attribute of the nonlinear controllers we define next is that they can be expressed as convex optimization problems that can be solved quickly, e.g., linear and quadratic programs. 

With the goal of driving the error signal $\vec{e} = (\vec{e}_q,\dot{\vec{e}}_q)$ to zero exponentially, consider a CLF $V$ satisfying: 
\begin{align}
k_1 \| {e} \|^c \leq V({e}) \leq k_2 \| {e} \|^c & \nonumber\\
\inf_{u \in \mathcal{U}} 
\dot{V}({e}, {u})  \leq - \lambda V({e})  & 
\label{eqn:Vdotcond}
\end{align}
for $c, k_1, k_2, \lambda  > 0$. Importantly, due to the affine nature of the dynamics \eqref{eqn:affinecontrolsys}, $\dot{V}$ is affine in the input $u$:
$$
\dot{V}({e}, {u})   = \frac{\partial V}{\partial {e}} 
\left(  
f(x) + g(x) u - \dot{{r}}
\right)
$$
and can therefore be expressed as a quadratic program (QP) when $\mathcal{U} = \R^p$: 
\begin{align}
\label{eqn:QPV}
 { u_{fb}}({e}) = \argmin_{{u} \in \R^p} & ~ \| {u} - { u_{ff}} \|^2  \\
 \text{s.t.} \quad  & ~ 
\dot{V}({e}, {u})  \leq - \lambda V({e})  \nonumber
\end{align}
which computes a minimal deviation from the desired feedforward control input $u_{ff,t}$ while tracking the reference signal exponentially: ${e}_t = x_t - {r}_t \to 0$. Importantly, there are many variations of QP-based controllers that are used for real-time control in robotic systems, i.e., those utilizing the dynamics as a constraint and can, therefore, account for constraints on forces and moments in real-time.  In all cases, the fact that these are QPs means that they can be implemented in real-time at loop rates of $1 kHz$ or greater, even on complex robotic systems like walking robots.



\subsection{Multi-Rate Layered Control Architectures}
\label{sec:multirate}

In LCAs, planning layers typically operate using \emph{reduced order models}.  These are simpler, usually lower dimensional, representations of the components of the full order dynamics~\eqref{eqn:eom} of interest, designed to capture essential behavior needed for the control task. As such, reduced order models are often application dependent, and their generation is often heuristic in nature: some of the most common reduced order models used for robotic control are the single integrator, double integrator, and unicycle.  These are often leveraged for control synthesis in the context of kinematic models, e.g., for mobile robots.  The overarching goal in the design of these reduced order models is the ability to generate reference signals that are (approximately) dynamically feasible and can be tracked well by a real-time feedback controller. 

This trajectory generation is typically performed over a longer horizon, requiring more computation time: it is therefore natural to view the interplay of trajectory generation and feedback control through the lens of \emph{multi-rate layered control architectures} i.e., through the lens of LCAs for which the controllers at different layers operate at different frequencies or rates.  This can be captured using continuous models (e.g., singular perturbation theory \cite{kokotovic1999singular,zhang2014singular}).  Yet in the case of robotic systems it is advantageous to be more concrete about the time scale separations present between layers.  One way to explicitly capture this is through the use of discrete time reduced order models at the planning layer, and continuous time full order models at the real-time feedback control layer.
An additional advantage of discrete time models at the planning layer is that they can provide an effective means of generating reference trajectories---this discrete instantiation better allows for planning forward in time, e.g., through MPC.   

To that end, consider a linear discrete time reduced order model
\begin{equation}
\label{eqn:romdiscrete}
y(k + 1) = A_{rom} y(k) + B_{rom} v(k)
\end{equation}
that will be used by the trajectory planning layer, where the reduced order model state $y \in \R^{p}$ and inputs $v \in \R^s$ are typically (but not necessarily) of a lower dimension than full order state $x \in \R^n$ and inputs $u \in \R^m$, i.e., $p \leq n$ and $s \leq m$.   The reduced order state is often related to the full system state via a projection map: $\Projx(x) = y$.  For robotic systems a commonly used projection is $\Projx(x) = \Projx(q,\dot{q}) = q$, i.e., one considers reduced order models on the configuration variables only.   Additionally,  $v \in \R^{s}$ is an auxiliary input to the reduced order model that is used to generate a reference signal sent to the real-time controller---we use $v$ to denote this auxiliary input, as they can often be interpreted as velocity commands. Analogously, we typically require an embedding of the auxiliary input $v$ into the full order dynamics~\eqref{eqn:affinecontrolsys} via $\Projv(v) = u$. 
While this subsection focuses on discrete time linear reduced order models, neither feature (discrete time, linear) is essential.  In the next subsection~\nameref{sec:cts_multilayer}, we explore the use of nonlinear continuous time reduced order models within multi-rate LCAs.

Proceeding with the discrete time linear reduced order model~\eqref{eqn:romdiscrete}, we follow the approach proposed in intermediate subproblem~\eqref{eq:master-redundant-1a} to couple the discrete time reduced order model state $y(k)$ and the underlying continuous time dynamics~\eqref{eqn:affinecontrolsys}.  
As in~\cite{rosolia2020multi}, we model the dynamics of the combined multi-rate LCA composed of a planning layer operating in discrete time with sampling period $\tau$ together with the full order control system~\eqref{eqn:affinecontrolsys} defined on $\mathcal{T} := \bigcup_{k\in \mathbb{N}_{\geq 0}} \mathcal{T}_k$, with $\mathcal{T}_k := (k\tau, (k+1) \tau)$,  as follows: 
\begin{align}
\label{eqn:multiratedynamics}
\textbf{Slow:} ~ & ~ y(k+1)  =   A_{rom} y(k) + B_{rom} v(k) ,   \qquad k \in \mathbb{N}_{\geq 0} \nonumber\\
\textbf{Fast:} ~ & ~  \dot{x}_t  =   f(x_t) + g(x_t) (u_t + u_{ff}(k)) ,   \qquad \quad t \in \mathcal{T}_k \nonumber\\
\textbf{Coupling:} ~ & ~ y(k) = \Projx(x_{k \tau }), \quad u_{ff}(k) = \Projv(v(k))
\end{align}
Here the planning layer operates at a slow loop rate, defined by the sampling rate $\tau$, in discrete time on a reduced order model, while the real-time feedback controller operates at a fast time scale (represented by a continuous time evolution).  Analogous to the coupling constraints in subproblem~\eqref{eq:master-redundant-1a},  the reduced and full order models are coupled via $\Projx$ and $\Projv$ where $\Projv(v(k))$ is held constant over the interval $(k\tau, (k+1)\tau)$ on which the ``fast'' low layer feedback controller operates, and $\Projx(x_{k\tau})$ is used to update the current state of the reduced order model every discrete step.  

The goal is to synthesize controllers $u$ and $v$ for the fast and slow dynamics in a synergistic fashion to achieve an overall control objective.  In particular, suppose we synthesize a controller, $v(k) = v_{fb}(y(k))$, that achieves a control objective for the reduced order linear model~\eqref{eqn:romdiscrete}, i.e., the closed loop system $y(k+1) = A_{rom} y(k) + B_{rom} v_{fb}(y(k))$ drives $y(k) \to y_{g}$ for a desired goal state $y_{g}$.  The evolution of the discrete time system can also be used to give a set of tracking goals for the real-time feedback controller expressed by the error%
\footnote{Note that this causes a discrete jump in the error every discrete step---to avoid this, a smooth function of time ${r}_t(k)$ on $\mathcal{T}_k$ can be defined such that $\vec{r}_{k\tau}(k) = y(k)$ and $\vec{r}_{(k+1)\tau}(k)(k) = y(k+1)$ wherein the error term becomes $\vec{e}_k = (\Projx(x) - \vec{r}_k)$. This can be achieved by converting the discrete time system $y(k+1) = A_{rom} y(k) + B_{rom} v(k)$ into a continuous time system $\dot{y} = A^c_{rom} y + B^c_{rom} v$ with $v$ implemented in a sample and hold fashion, i.e.,  $A^c_{rom}$ and $B^c_{rom}$ are defined by: 
$$
A^c_{rom} = \frac{1}{T}\log(A_{rom}), \qquad
B^c_{rom} = A^c_{rom}  (A_{rom} - I)^{-1} B_{rom}
$$
when the $\log$ and inverse are well defined.  In the case when they are not, one can assume the discrete time system came from Euler integration, $A_{rom} = (I + A^c_{rom} T)$ and $B_{rom} = B^c_{rom} T$, to obtain $A^c_{rom} = (A_{rom} - I)/T$ and $B^c_{rom} = B_{rom}/T$.  In either case, the result is therefore a continuous reference signal: 
$$
{r}_t(k) = 
e^{A^c_{rom}(t-k\tau)} y(k)
+ 
\left( 
\int_{k\tau}^{t} e^{A^c (t - s)} d s
\right) B^c v_{fb}(y(k)), \qquad t \in \mathcal{T}_k
$$
Alternatively, the discrete time system \eqref{eqn:romdiscrete} can be replaced by a continuous time system at the slow control layer, as was done in \cite{rosolia2020multi}, to generate a smooth reference signal \emph{a priori}.} %
terms: ${e}_{k,t} = (\Projx(x_t) - y(k+1))$. 
A controller can be synthesized that ideally drives this error to zero, e.g., a linear controller as in \eqref{eqn:PDcontrol} which here takes the form: 
\begin{eqnarray}
\label{eqn:linearmulti}
u({e}_{k,t},x(k\tau)) & = & \overbrace{\Projv \circ v_{fb} \circ \Projx(x_{k\tau})}^{u_{ff}(x_{k\tau})}  \\
&& \quad + \underbrace{\left(\frac{\partial \Projx}{\partial x}\right)^T K_P(\Projx(x_t) - y(k+1))}_{u_{fb}(e_{k,t})},  \nonumber
\end{eqnarray}
or using the Lyapunov controller in \eqref{eqn:QPV} with ${e}$ replaced by ${e}_{k,t}$ and ${u}_{ff}$ replaced by $u_{ff}(x(k\tau))$ for $t \in \mathcal{T}_k$.  The end result is the closed loop multi-rate system:
\begin{align}
\label{eqn:slowclosedloop}
&\textbf{Slow:} ~ \textrm{For }  k \in \mathbb{N}_{\geq 0}: \\
& \qquad  y(k+1)  =   A_{rom} \Projx(x_{k\tau}) + B_{rom} v_{fb} \circ \Projx(x_{k\tau}),   \quad  \nonumber\\
\label{eqn:fastclosedloop}
& \textbf{Fast:} ~ \textrm{For } \vec{e}_{k,t} = (\Projx(x_t) - y(k+1)) \textrm{ and }  t \in \mathcal{T}_k  \\
&  \qquad  \dot{x}_t  =   f(x_t) + g(x_t)(u_{ff}(x_{k\tau})+u_{fb}({e}_{k,t})).  ~  \nonumber
\end{align}
Here the state $x_{k\tau}$ at the beginning of the sampling period informs the next iteration of the slow dynamics, while the slow dynamics informs the fast dynamics through the feedforward input (which depends on $v_{fb}$) and ${e}_{k,t}$ which drives the system to the next desired setpoint $y(k+1)$ over the interval $\mathcal{T}_k$. 
To provide a specific example of the generation of closed-loop policies, we begin with slow controller synthesis viewed as a planning problem.

\paragraph{Slow trajectory generation}
We view trajectory generation as a planning problem~\eqref{eq:planning} which can be solved using MPC.  In particular, consider a goal state $y_g$ for the reduced order model, obtained for example from a higher decision layer, with the objective of synthesizing a controller that achieves this objective subject to a safety constraint expressed as state constraints $\mathcal{S} = \{ y \in \R^p : h(y) \geq 0\}$ and input constraints $\mathcal{V}$, i.e., the system must evolve such that $y(k) \in \mathcal{S}$ and $v(k) \in \mathcal{V}$ for all $k\geq 0$.  To this end, we can formulate a MPC problem resembling that proposed in equation~\eqref{eq:planning} with a $N\geq 1$ planning horizon and positive (semi)definite cost matrices $Q,Q_K\succeq 0$ and $R\succ 0$ as a QP: 
\begin{align}\label{eq:MPC}
    \mathrm{minimize} ~ &  ~ \sum_{i=k}^{k+N-1} \left(\| y(i) - y_g \|^2_Q  + \| v(i) \|_R^2 \right) \\
    & \qquad \quad +  \|y(k+N) - y_g\|^2_{Q_N}  \nonumber \\
         \text{s.t.} \quad & ~  y(i+1) = A_{rom}y(i) + B_{rom} v(i),   \nonumber\\
         & ~  y(k) = \pi_x(x_{\tau k})  \nonumber\\
         & ~  y(i),y(k+N) \in \mathcal{S}, \quad v(i),v(k+N) \in \mathcal{V}, \nonumber\\
         &  ~  i = k, \ldots, k+N - 1,  \nonumber
\end{align}
with $\| y\|^2_Q = y^T Q y$.  At each discrete time step $k\geq 0$, problem~\eqref{eq:MPC} is solved to produce a sequence of nominal reduced order model states $\vec y^\star=y^\star(k:k+N)$ and inputs $\vec v^\star = v^\star(k:k+N-1)$.  The planning layer controller is then chosen as $v(k) = v_{MPC}(y(k)) := v^\star(k)$, as is the standard approach in MPC.  

\paragraph{Fast and safe feedback control}
The solution $(\vec y^\star, \vec v^\star)$ to each MPC subproblem solved at discrete time step $k$ can further be leveraged to synthesize a lower layer real-time feedback controller.  In particular, we can transmit these solutions to the fast lower layer to define the error signal ${e}_{k,t} = (\Projx(x_t) - y^\star(k+1))$, which can in turn be driven to zero using the Lyapunov controller~\eqref{eqn:QPV} (replacing ${e}$ with ${e}_{k,t}$) and using $u_{ff}(k) = \Projv(v_{MPC}(y(k)))$.  While this will drive $\Projx(x) \to y^\star(k+1))$ (the next step produced by the MPC problem) with the input from the MPC problem $v_{MPC}(y(k))$ as a reference, there is no guarantee that the safety constraints will be satisfied over the time interval $\mathcal{T}_k$.   To address this shortcoming, we can combine the Lyapunov controller~\eqref{eqn:QPV} with the CBF controller~\eqref{eqn:QPh} into an optimization problem that resembles optimization problem~\eqref{eq:control0}.  However, in this case we can exploit the control-affine structure of the dynamics~\eqref{eqn:affinecontrolsys} to obtain a QP through the use of Lyapunov and barrier functions: 
\begin{align}
\label{eqn:QPVh}
 {u_{fb}}({e}_{k,t}) = \argmin_{{u} \in \R^p, ~\delta \in \R} & ~\| {u} - \Projv(v_{MPC}(y(k))) \|^2  + p \delta^2 \\
 \text{s.t.} \quad  & ~
 \dot{V}({e}_{k,t}, x_t,  {u}) \leq - \lambda V({e}_{k,t}) + \delta  \nonumber\\
 & ~ \dot{h}_{\Projx}(x_t, {u}) \geq - \alpha h_{\Projx}(x_t)  \nonumber
\end{align}
with $h_{\Projx} := h \circ \Projx$ assumed to be a valid CBF: 
$$
\sup_{u \in \R^{m}}  \underbrace{\left[
\frac{\partial h}{\partial y}\frac{\partial \Projx}{\partial x} \left( f(x) + g(x) u \right)
\right]}_{\dot{h}_{\Projx}(x_t,u)} \geq - \alpha h_{\Projx}(x)
$$
Here $p > 0$ is a penalty associated with the relaxation term $\delta>0$, which can be interpreted as an instantaneous analog to the tracking cost~\eqref{eq:control0} that penalizes $\| x - r\|_Q^2$.   Note that input constraints can also be added to the QP \eqref{eqn:QPVh}, but this would require a relaxation of the CBF condition or the input constraints to ensure feasibility (see \cite{garg2021multi} which considers the interplay between continuous and discrete dynamics in the context of input constraints). 


The solution ${u_{fb}}({e}_{k,t})$ to the QP~\eqref{eqn:QPVh} can be used in the multi-rate dynamics~\eqref{eqn:fastclosedloop}, which when combined with the MPC problem~\eqref{eqn:slowclosedloop}, yields a closed loop multi-rate controller instantiated via a LCA. The power of the closed loop multi-rate LCA, as opposed to a single layer feedback controller using Lyapunov and barrier functions, is that information from the MPC problem encodes knowledge of future system behavior via the reduced order model. Conversely, the fast layer accounts for the full nonlinear dynamics of the system that are not present in the slow layer. Thus, we are able to synthesize a LCA that leverages the nonlinear dynamics of the system at the real-time feedback control layer, while looking ahead to future behaviors defined in the trajectory planning layer, via a synergistic coupling of the two.   Formal guarantees for the multi-rate LCA presented here can be found in~\cite{rosolia2022unified,rosolia2020multi}.

\begin{sidebar}{Multi-Rate LCAs in Practice}
\section[Multi-Rate LCAs in Practice]{}\phantomsection
\label{sidebar:multirate}
\setcounter{sequation}{0}
\renewcommand{\thesequation}{S\arabic{sequation}}
\setcounter{stable}{0}
\renewcommand{\thestable}{S\arabic{stable}}
\setcounter{sfigure}{0}
\renewcommand{\thesfigure}{S\arabic{sfigure}}

To highlight the practical impact of the multi-rate LCAs we presented, we provide an overview of successful experimental implementations in safe navigation, safe locomotion, and data-driven locomotion.  In all cases, a multi-rate LCA facilitates the ability to realize controllers in practice.  The commonality of approaches, and architectures more specifically, in these disparate applications on different hardware platforms shows the broad applicability of LCAs for robot control.

\subsection{Safe Navigation}  

Consider the problem of safe navigation with a ground robot \cite{agrawal2021constructive}---in this case, both a wheeled vehicle and a quadruped.  Following Example~\ref{ex:mpcdubins}, differential flatness of the Dubins' car is utilized to create a linear system~\eqref{eqn:unicyclelinear} that is discretized as in equation~\eqref{eqn:discretedubins}.  Using this discrete time linear reduced order model, an MPC problem is formulated as in equation~\eqref{eq:MPC}, wherein safety is enforced (avoiding obstacles) while planning a path towards a goal.  This generates a discrete time reference trajectory that is sent to the Dubins' car and tracked with a real-time controller as in equation~\eqref{eqn:feedbackdubinsflat}.  This paradigm is illustrated in Fig.~\ref{fig:multi-rate}.  In particular, the discretely updated reference signals are shown (in red) along with the tracking of these reference signals by the real-time controller (in green).  Importantly, the input to the Dubins' car model can be viewed as a reference velocity that can be tracked on hardware with onboard controllers.  This further layering allows for the experimental deployment on both a wheeled vehicle and a quadruped, as again shown in Fig.~\ref{fig:multi-rate}. 

\sdbarfig{\label{fig:multi-rate}
    \includegraphics[width=1\columnwidth]{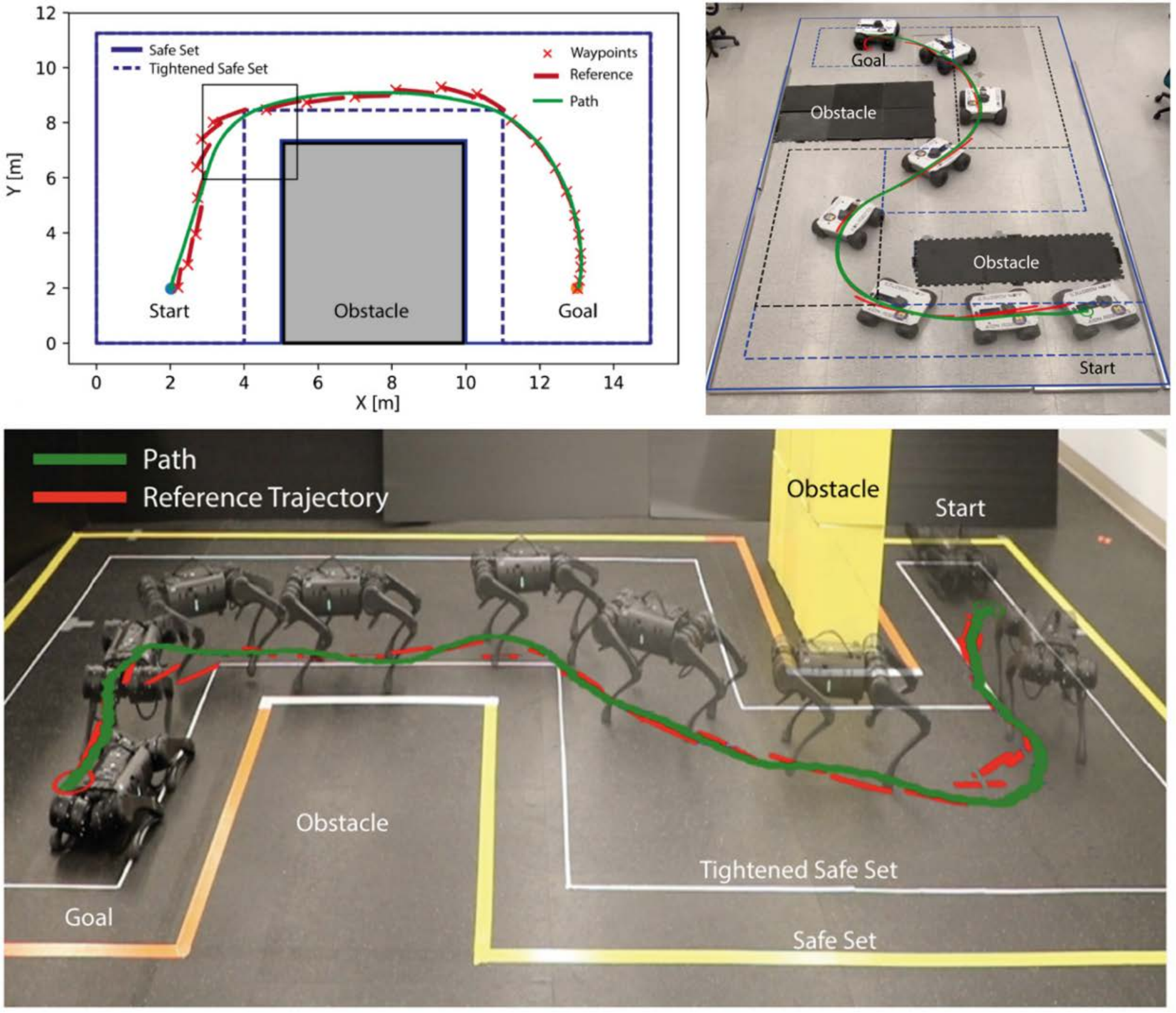}
    }{\footnotesize Experimental demonstration of a multi-rate LCA using the Dubins' car and differential flatness (from~\cite{agrawal2021constructive}).  Reference trajectories are generated with a discrete time linear reduced order model and tracked by the Dubins' car model (top left).  Additionally, these signals can be passed to hardware and tracked onboard with real-time controllers, both on a wheeled vehicle (top right) and a quadruped (bottom). }

\subsection{Safe Locomotion}

As noted throughout this section, LCAs provide an effective paradigm for enforcing safety constraints on complex robotic systems by enforcing safety (framed as set invariance) at both the planning layer (e.g, in the MPC problem \eqref{eq:MPC} as a state constraint), and at the real-time control layer via a CBF (e.g., as in the QP \eqref{eqn:QPVh}).  To demonstrate this, consider the \emph{stepping stone problem} where the goal is for a legged robot to precisely place its feet on a series of stepping stones.  This is safety-critical in that if this foot placement target is missed by the feet, the robot will fall.  Additionally, it requires a layered approach, in that the system must maintain safety while also remaining dynamically balanced.  

In \cite{grandia2021multi}, a LCA formulation was implemented on a quadrupedal robot (ANYmal) to realize stepping stone behavior experimentally.  This is illustrated in Fig.~\ref{fig:steppingstones}.  In particular, a safety constraint is implemented at the planning layer via a CBF (via MPC with a kinematic reduced order model), and at the real-time control level (as a CBF constraint in a whole body controller).  Implementing CBFs at both layers resulted in no failures (missing the stepping stone) over 140 steps.  Without the LCA framework more failures were observed, i.e., just enforcing a CBF constraint at the real-time layer leads to 5 failures, while implementing CBFs at only the planning layer leads to 6 failures.  This demonstrates the practical utility of LCAs for safety-critical systems.

\sdbarfig{\label{fig:steppingstones}
    \includegraphics[width=1\columnwidth]{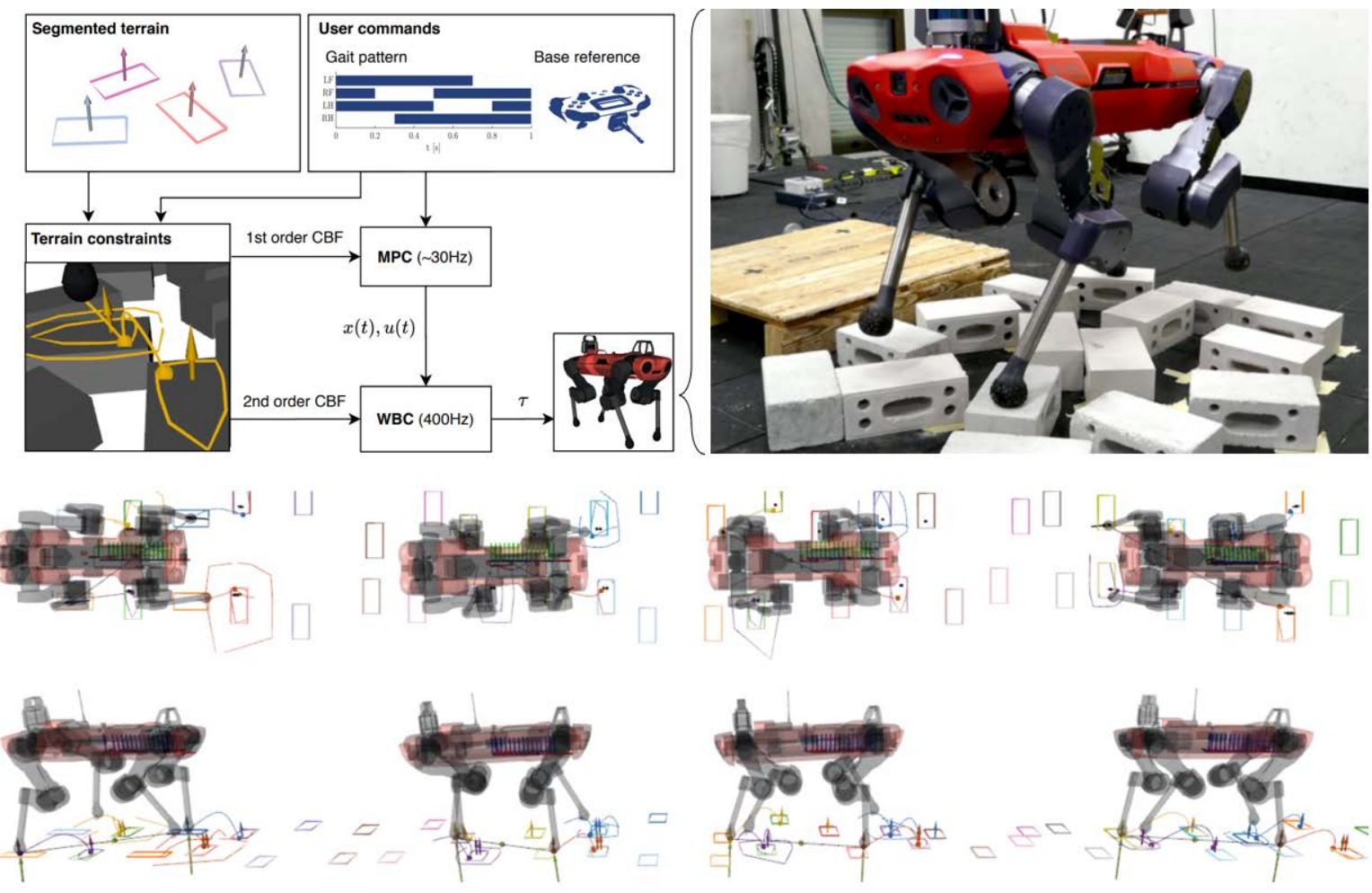}
    }{\footnotesize Experimental demonstration of an LCA (top left) used to realize dynamic walking on stepping stones (from \cite{grandia2021multi}).  This framework was implemented on ANYmal (top right) with the results being walking of the form illustrated (bottom).}

\subsection{Data-Driven Locomotion} 

Finding reduced order models on which to instantiate multi-rate control can be challenging, often requiring domain specific knowledge.  This can be addressed by learning reduced order models via data-driven methods.  To provide an example of this paradigm, consider again the problem of legged locomotion where reduced order models are used at the planning layer---the goal is to learn this model and deploy the learned model experimentally.  
\hfill \emph{(continued on next page)}
   \end{sidebar}

   \begin{sidebar}{\continuesidebar}
   \renewcommand{\thesequation}{S\arabic{equation}}
   \renewcommand{\thestable}{S\arabic{stable}}
   \renewcommand{\thesfigure}{S\arabic{sfigure}}

\subsection{Multi-rate LCAs in Practice}
Following~\cite{fawcett2022toward,fawcett2022distributed}, we learn a linear reduced order model at the planning layer and leverage this model to pose an MPC problem.  In particular, given sufficiently rich (persistently excited) data collected from the robot, Hankel matrices can be used to exactly determine the forward evolution of linear time invariant systems~\cite{willems2005note}.  Following~\cite{coulson2019data}, given a user defined reduced order model state $y \in \R^p$ and input $v \in \R^s$, one can define the \emph{data-driven state transition matrix} $\mathcal{G}(\texttt{data})$ over N-steps: 
\begin{sequation}
\label{eqn:datadrivendyn}
y_{k:k+N} = \mathcal{G}(\texttt{data})
\begin{bmatrix}
    v_{\rm ini} \\
    y_{\rm ini} \\
    v_{k:k+N} 
\end{bmatrix}
\end{sequation}
Here, $v_{\rm ini}$ and $y_{\rm ini}$ are the reduced order inputs and state observed in the past over an estimation horizon $T_{\rm ini}$.   Equation~\eqref{eqn:datadrivendyn}, which defines linear relationships between the control input and state over the next $N$ steps ,can be used as a constraint in the MPC problem \eqref{eq:MPC} instead of the explicit dynamic constraint $y(k+1) = A_{rom} y(k) + B_{rom} v(k)$.

This approach was experimentally deployed on a quadruped robot, as shown in Fig.~\ref{fig:datadrivenquad}.  The robot considered has 18 degrees of freedom, and thus the state is 36 dimensional ($x \in \R^{36}$).  A reduced order model is considered with a 10 dimensional state $y \in \R^{10}$ consisting of select body positions, velocities, and orientations.  The reduced order input $v \in \R^{12}$ consists of ground reaction forces (GRFs).  With this reduced order model, data is collected, $\mathcal{G}(\texttt{data})$ is computed and equation~\eqref{eqn:datadrivendyn} is leveraged in a LCA via an MPC problem at the planning layer, which is subsequently implemented on the robot via a nonlinear real-time controller similar to \eqref{eqn:QPVh}.  The end result is robust data-driven dynamic walking. 

\sdbarfig{\label{fig:datadrivenquad}
    \includegraphics[width=1\columnwidth]{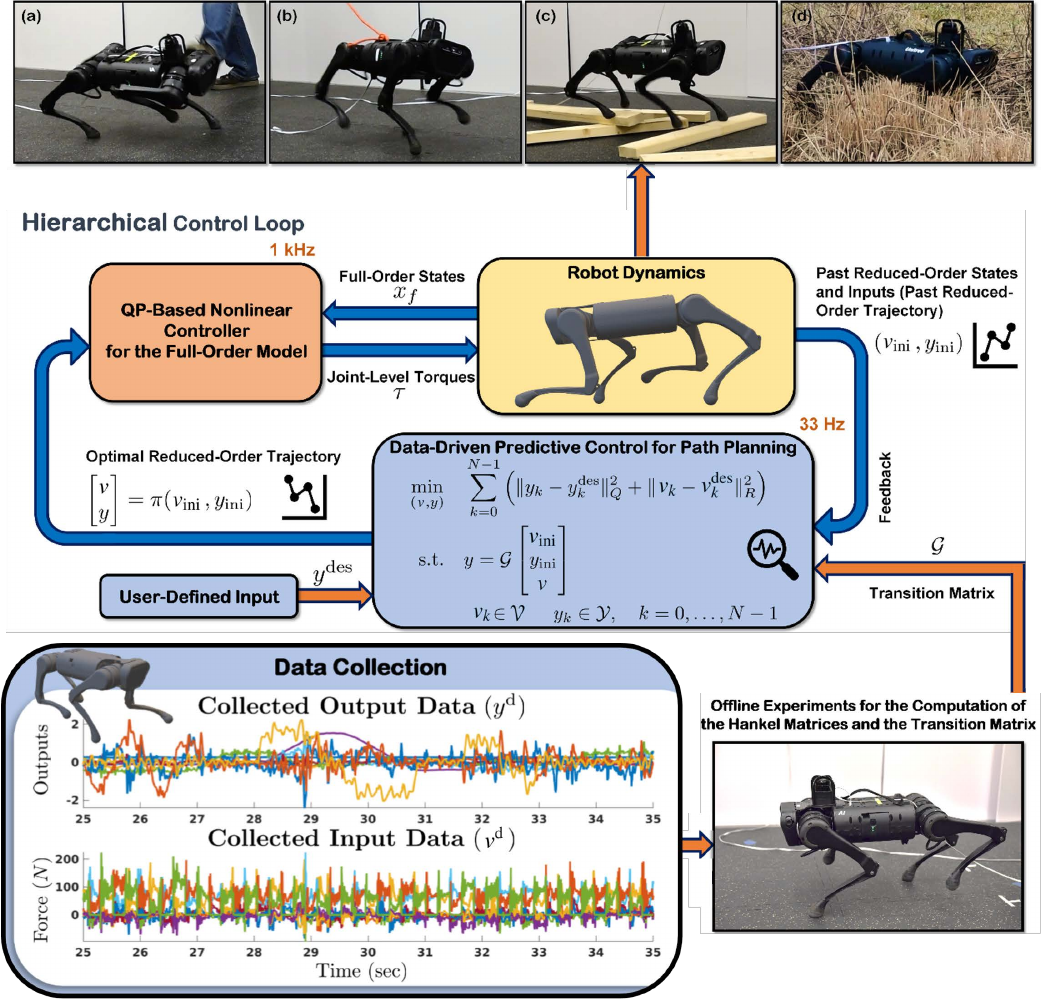}
    }{\footnotesize Realization of a data-driven LCA on a quadrupedal robot (from \cite{fawcett2022toward}).  Data is collected from the robot (bottom) to determine the data-driven state transition matrix $\mathcal{G}(\texttt{data})$, which is then used in an MPC problem to plan trajectories (middle).  The reduced order model state and input is passed to a nonlinear optimization based controller at the real-time layer to control the robot. The result is dynamic walking (top) that is robust to pushes (a), external pulls (b), rough terrain (c), and natural terrain (d). }

\end{sidebar}

\begin{my_exam}[Running example: robot navigation]
\label{ex:mpcdubins}
Consider again the Dubins' car which, for the moment, we view as the full order dynamics. We aim to instantiate a discrete time planning layer via MPC.  To determine the corresponding reduced order model, we leverage that the dynamics of the Dubins' car $\dot{q}_t = f_{rom}(q_t) u_t$ with $q = (x_1,x_2,\theta)$ and $u=(u_1,u_2)$, as given in \eqref{eqn:unicycle}, are differentially flat per Example \ref{ex:exampleflat}.  For the flat output $\xi = (x_1,x_2)$, denote the relationships between the states, inputs and flat outputs by: $q = q_{\flat}(\xi,\dot{\xi})$ and $u = u_{\flat}(\xi,\dot{\xi},\ddot{\xi})$. Note we make a slight deviation from the notation used in Example~\ref{ex:exampleflat} to be consistent with the configuration space notation defined in this section, and use $q$ to denote the state rather than $x$.

To apply the approach outlined in this section, we can forward integrate the flat continuous time linear dynamics~\eqref{eqn:unicyclelinear} over the time interval $\mathcal{T}_k=[k\tau,(k+1)\tau]$ to obtain the discrete time reduced order model
\begin{eqnarray}
\label{eqn:discretedubins}
y(k+1) =  \underbrace{e^{A_{rom}^c \tau}}_{A_{rom}} y(k) + 
\underbrace{\int_{0}^\tau e^{A_{rom}^c s} B_{rom}^c ds}_{B_{rom}} v(k).
\end{eqnarray}
Here $y(k) = (\xi(k),\dot{\xi}(k)) \in \R^4$ and $v(k)  = \ddot{\xi}(k) \in \R^2$.  We note that here the reduced order model is actually higher dimensional but is ``reduced'' in complexity by being linear. Utilizing this system, a feedback controller $v_{fb}(y(k))$ can be synthesized.  For example, this can be chosen to be the result of the MPC problem~\eqref{eq:MPC}, i.e., $v_{fb}(y(k)) = v_{MPC}(y(k))$.  The result is the error ${e}_{k,t} = (q_t - q_{\flat}(y(k+1)))$ and a feedforward input $u_{ff}(y(k)) = u_{\flat}(y(k),v_{fb}(y(k)))$.  This can be used to synthesize a linear feedback controller of the form~\eqref{eqn:linearmulti}, modified slightly to exploit differential flatness, as in~\eqref{eq:dubins-feedback}: 
\begin{equation}
\label{eqn:feedbackdubinsflat}
\begin{array}{rl}
u({e}_{k,t},y(k)) &= u_{ff}(y(k)) + 
K_P {e}_{k,t} \\
&=: u_{ff}(y(k)) + 
u_{fb}({e}_{k,t}),
\end{array}
\end{equation}
for $K_P$ a positive definte matrix.

The final feedback controller is given by setting $y(k) = (x_{1,k\tau},x_{2,k\tau},\dot{x}_{1,k\tau},\dot{x}_{2,k\tau})$.  This paradigm is deployed experimentally in~\nameref{sidebar:multirate}.  Alternatively, in the expression above, we could consider a continuous reference signal $\vec{e}_t = (q_t - q_{\flat}(y_t))$ where $y_t$ is the solution to \eqref{eqn:unicyclelinear} given a feedback controller $v = K_{fb} y $, i.e., by solving $\dot{y} = (A + B K_{fb}) y $ with initial condition $y(k\tau) = (x_{1,k\tau},x_{2,k\tau},\dot{x}_{1,k\tau},\dot{x}_{2,k\tau})$.   This paradigm in deployed experimentally as described in~\nameref{sidebar:contmidlevel}, and is discussed in more detail in the next subsection.
 
\end{my_exam}

\subsection{Continuous Time LCAs}  
\label{sec:cts_multilayer}
We observe that we can further layer the control architecture defined in Example~\ref{ex:mpcdubins} by viewing the unicycle as a reduced order model wherein continuous multi-rate control can be applied.  That is, we can view $\mu_{fb}({e}_{k,t},y(k))$ as a reference velocity (see Example \ref{ex:cbfdubins} below) which we want a more complex robot to track, i.e., a quadruped or drone as described in~\nameref{sidebar:contmidlevel}. The end result is a three layer architecture with two planning layers and one feedback control layer: (1) a slower discrete time planning layer using a linear model, (2) an intermediate reference signal generation layer using a continuous time unicycle model, and (3) a fast feedback control layer for tracking of the reference signal by the underlying complex robotic system.  This observation highlights that layers can often be added in a fairly modular fashion, allowing for the benefits of each layer to be enjoyed.  This subsection explores the bottom two layers of the LCA described above, namely the interplay between a continuous time reference signal generation layer and a real-time feedback control layer, in more detail.  Implicit throughout is the assumption that the loop rates at each layer, and the communication between layers, happens sufficiently fast.  We focus on safety-critical systems, wherein safe reference signals are generated by the trajectory generation layer to be tracked by the real-time control layer.  We show that formal guarantees of safety can be obtained for these LCAs and, importantly, this architecture enables theory to be widely deployed in practice.


\paragraph{Safe reference signal generation}
To that end, consider a continuous time reduced order model
\begin{eqnarray}
\label{eqn:from}
\dot{y}_t = f_{rom}(y_t,v_t),
\end{eqnarray}
where as above, the reduced order dynamics $f_{rom}$, state $y_t\in\R^p$, and auxiliary input $v_t\in\R^s$, are chosen to be simpler than the full equations of motion~\eqref{eqn:eom}, while nevertheless capturing the essential features of the control problem at hand.  We recall that we assume that the reduced order state $y$ is related to the full system state $x$ via the projection $\Projx(x)=y$, and that the auxiliary input $v$ can be embedded into the full dimensional input space via the embedding $\Projv(v)$.

\begin{figure}[t]
    \centering
    \includegraphics[width=1\columnwidth]{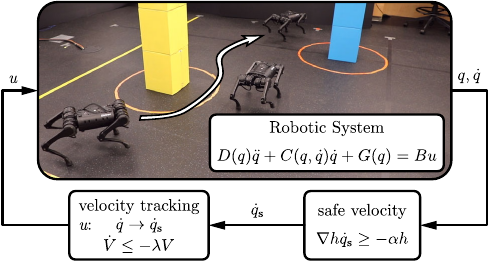}
    \caption{An LCA (from \cite{molnar2021model}) for controlling robotic systems via continuous reference signal generation.  In this case, a safe reference velocity is generated that is tracked by a real-time controller realized on the robot. }
    \label{fig:modelfree}
\end{figure}

Assume now that the reduced order model~\eqref{eqn:from} is used to generate a desired reduced order state trajectory $\vec y^d$ and a corresponding feedback law $v^d(y_t)$ for the auxiliary input, e.g., via the techniques described in the previous subsection.  Now consider the objective of ensuring the reference signal satisfies a safety constraint, encoded by making the set $\mathcal{S} = \{ y \in \R^{p} : h(y) \geq 0\}$ forward invariant, for some differentiable function $h:\R^p\to\R.$  We can leverage CBFs if $h$ satisfies the CBF condition with respect to the reduced-order dynamics: 
\begin{eqnarray}
\label{eqn:CBFcond}
\sup_{v \in \R^{s}} \dot{h}(y,v)
= \sup_{v \in \R^{s}} \left[
\frac{\partial h}{\partial y} f_{rom}(y,v)
\right] \geq - \alpha h(x)
\end{eqnarray}
where $\alpha > 0$ is a positive constant.\footnote{Note that, more generally, $\alpha$ can be chosen to be an extended class $\mathcal{K}$ function---we opt for a positive constant for simplicity of exposition}  If $f_{rom}$ is affine in the auxiliary control input $v$ this inequality can be expressed as a QP of the form \eqref{eqn:QPV}, with the result being a safety filter operating on the reduced order model within the reference signal generating layer:
\begin{align}
\label{eqn:QPh}
 v_{safe}(y) = \argmin_{v \in \R^{s}} & \| v - v^d(y) \|^2  \\
 \text{s.t.} \quad  & 
\dot{h}(y,v) \geq - \alpha h(y)  \nonumber
\end{align}
This safety filter can then be integrated into reference signal generation in a variety of ways. 
%
%
For example, forward integrating the closed loop dynamics $\dot{y}_t = f_{rom}(y_t,v_{safe}(y_t))$ to generate a reference signal $\vec r_{y}$ with corresponding error: $e_{y,t} = y_t - {r}_{y,t}$, which can then be tracked with a linear or nonlinear controller, i.e., replacing ${e}$ in equation~\eqref{eqn:QPV} with ${e}_y$.  
Alternatively, the safe input $v_{safe}$ can be tracked by the real-time controller as described below.  

\paragraph{Real-time feedback control}
To provide a concrete example of the use of continuous time reduced order models at the planning layer being coupled with real-time feedback controllers, consider the case when we have a kinematic reduced order model with $y = q$, i.e., our reduced order model operates on the configuration variables, $\dot{q}_t = f_{rom}(q_t,v_t)$, where now the auxiliary input $v_t$ is naturally associated with the generalized velocities $\dot q_t$.  Consider a safe velocity, $v_{safe}(q)$, generated from the QP \eqref{eqn:QPh}.  
Following~\cite{molnar2021model}, assume that this velocity is passed to a real-time controller via the error signal $\dot{e}_{safe,t} := \dot{q}_t - v_{safe}(q_t)$, i.e., the real-time controller takes the safe velocity from the reduced order model as a reference with the goal of tracking this reference signal.  Assume a real-time feedback controller $u_t = {u}_{fb}(x_t,v_{safe}(q_t))$ that can exponentially track this reference velocity, e.g., via the controller \eqref{eqn:QPV} with ${e}_t$ replaced by $\dot{{e}}_{safe,t}$, resulting in exponentially fast tracking:
\begin{eqnarray}
    \label{eqn:expconvergence}
\| \dot{{e}}_{safe,t} \| \leq M e^{-\lambda t} \| \dot{\vec{e}}_{safe,0} \| 
\end{eqnarray}
for $M,\lambda > 0$, with the error calculated along solutions of the closed loop system: $\dot{x}_t = f(x_t) + g(x_t){u}_{fb}(x_t,v_{safe}(q_t))$.  The following theorem adapted from~\cite{singletary2022safety} provides formal guarantees for the reduced order model-based LCA applied to the full order dynamics~\eqref{eqn:eom}.  For experimental implementations related to this formal result, see \nameref{sidebar:contmidlevel}.


\begin{my_thm}
\label{thm:kintofull}
Consider a control system \eqref{eqn:affinecontrolsys}, where $x = (q,\dot{q})$, and
a safe set $\mathcal{S} = \{ q \in Q : h(q) \geq 0\}$.  Assume that $h$ has bounded gradient,
i.e., there exists $K_h>0$ s.t.
${\left\| \frac{\partial h}{\partial q}\right\|}_2 \leq K_h$ for all
$q \in \mathcal{S}$.  Let $v_{safe}(q)$ be the safe velocity given by the
QP~\eqref{eqn:QPh}, with corresponding error $\dot{e}_{safe} = \dot{q} - v_{safe}(q)$ satisfying~\eqref{eqn:expconvergence}.  If $\lambda > \alpha$, safety is achieved for the full-order
dynamics \eqref{eqn:affinecontrolsys}:
\begin{eqnarray}
  \label{eqn:SMsafety}
  (q(0),\dot{e}_{safe}(0)) \in \mathcal{S}_M ~ \quad \implies ~ \quad q(t)
  \in \mathcal{S}, \quad \forall t \geq t_0, 
\end{eqnarray}
where
\begin{equation}
\label{eqn:SM}
  \mathcal{S}_M  = \left\{ (q,\dot{e}_{safe}) \in \R^{2n} ~ : ~ h(q) -
    \frac{K_h M}{\lambda - \alpha} \| \dot{e}_{safe}\|_2 \geq 0 \right\}. 
\end{equation}
\end{my_thm} 

Note that to certify ``fast-enough'' tracking by the real-time controller, a Lyapunov certificate can be used \cite{molnar2021model} (see Fig. \ref{fig:modelfree}):  assume the real-time feedback controller tracks the error $\dot{e}_{safe}$ per a Lyapunov function as in equation~\eqref{eqn:Vdotcond}: $\dot{V}(\dot{e}_{safe}) \leq - \lambda V(\dot{e}_{safe})$.  Then, for any differentiable $v_{safe}(q)$ satisfying the CBF condition \eqref{eqn:CBFcond}, safety for the full-order dynamics is achieved if $\lambda > \alpha$: 
$$
(q(0),\dot{e}_{safe}(0)) \in \mathcal{S}_M ~ \quad \implies ~ \quad q(t)
  \in \mathcal{S}_V, \quad \forall t \geq t_0
$$
where
\begin{align}
\mathcal{S}_V & = \left\{ (q,\dot{e}_{safe}) \in \R^{2n} ~ : ~ h_V(q,\dot{e}_{safe}) \geq 0 \right\} \nonumber\\
h_v(q,\dot{e}_{safe}) &  := -V(q,\dot{e}_{safe}) + \alpha_e h(q), 
\quad 
\alpha_e = \frac{(\lambda - \alpha)k_1}{K_h}
\nonumber
\end{align}
Interestingly, this result is established by synthesizing a CBF for the full system dynamics, $h_V$, using the CBF for the reduced order model, $h$, together with the Lyapunov function for the tracking controller, $V$. 

\begin{sidebar}{Continuous Time LCAs in Practice}
\section[Continuous Time LCAs in Practice]{}\phantomsection
\label{sidebar:contmidlevel}
\setcounter{sequation}{0}
\renewcommand{\thesequation}{S\arabic{sequation}}
\setcounter{stable}{0}
\renewcommand{\thestable}{S\arabic{stable}}
\setcounter{sfigure}{0}
\renewcommand{\thesfigure}{S\arabic{sfigure}}

To illustrate the practical consequences of Theorem  \ref{thm:kintofull}, we highlight how a common reduced order model can be used to achieve safety across a variety of robotic systems, including a drone, quadrupedal robot, manipulator, and full-scale automotive system.  This diverse set of robotic systems have very different underlying dynamics, and yet deploying a well designed LCA does not require direct knowledge of these dynamics, rather only ``good'' onboard tracking controllers that allow planning layers to operate on reduced order models rather than on the underlying complex system dynamics (as illustrated in Fig. \ref{fig:modelfree}).  

\subsection{Drones and Quadrupeds}
We wish to enforce collision avoidance to obstacles in the environment. To begin, consider the ``simplest'' kinematic model of a robot, a single integrator: 
\begin{sequation}
\label{eqn:singleint}
\dot{q}_t = f_{rom}(q_t,v_t) = v_t, 
\end{sequation}
obtained by setting $y_t = q_t$.  Collision avoidance  is encoded by the safety constraint $\mathcal{S} = \{ q \in \R^n : h(q) \geq 0\}$ for 
$$
h(q) = \| q - q_0 \| - r. 
$$
Here $q_0 \in \R^n$ is the centroid of the obstacle and $r$ its radius.  The safety filter~\eqref{eqn:QPh} can be expressed as the QP: 
\begin{sequation}
\label{eqn:QPsingleint}
\begin{aligned}
     v_{safe}(q) = \argmin_{v \in \R^{n}} & \| v - v^d(q) \|^2  \\
 \text{s.t.} \quad  & 
\frac{(q - q_0)^T}{\|q - q_0 \|} v \geq - \alpha (\| q - q_0 \| - r)  
\end{aligned}
\end{sequation}
where $v^d(q) = - K_P(q - q_g)$ is a desired velocity which drives the system to a goal position $q_g \in \R^n$.  In the case of planar collision avoidance $q,q_0,q_g \in \R^2$ (representing the spacial position in the plane)---this is the case that will be considered in the context of the experimental implementation on a drone and quadruped.

\sdbarfig{\label{fig:modelfreehardware}
    \includegraphics[width=1\columnwidth]{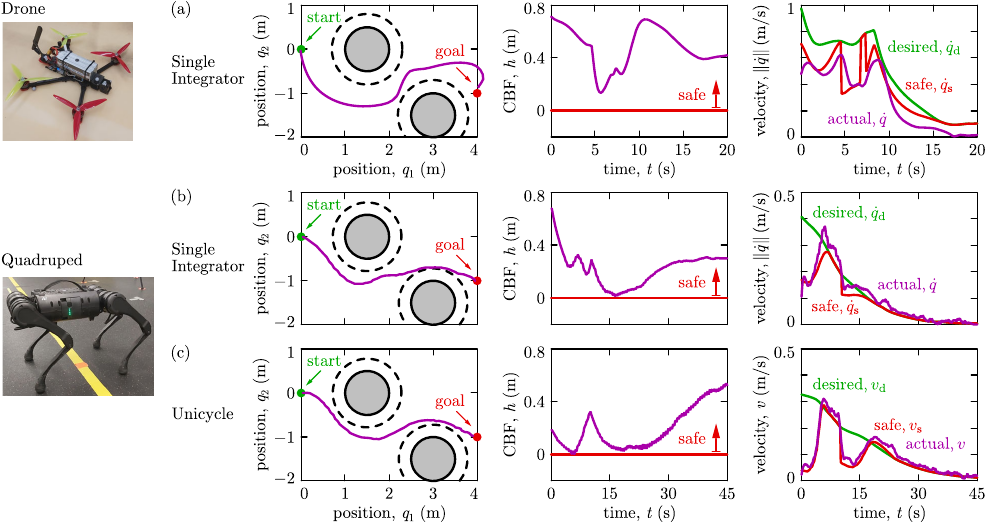}
    }{\footnotesize Experimental results on a drone and quadruped (from \cite{molnar2021model}).  Both the drone and quadruped use a single integrator \eqref{eqn:singleint} reduced order model and the corresponding QP \eqref{eqn:QPsingleint}.  Suitable integration into an LCA yields safe behavior.  Additionally, a unicycle (Dubins' car) model is used on the quadruped to also achieve safe behavior that is less conservative. }

The result of \eqref{eqn:QPsingleint} is a safe velocity $v_{safe}(q)$ that can be tracking with existing onboard tracking controllers.  This was implemented on both a drone and a quadruped hardware platform.  We highlight that while these platforms have dramatically different underlying dynamics, by leveraging a well designed LCA, the exact same safe reference velocity, $v_{safe}(q)$ can be used and tracked on both platforms.  The results can be seen in Fig. \ref{fig:modelfreehardware}: safety is achieved (as certified by $h(q)\geq 0$) for both the drone and quadruped tracking reference signals $v_{safe}(q)$ produced by \eqref{eqn:QPsingleint}.  For the quadruped, we can also use the Dubins' car as a reduced order model instead of the single integrator, as described in Example \ref{ex:cbfdubins}.  In this case, the resulting safe velocity $u_{safe}(q)$ produced by the QP \eqref{eqn:QPh} is tracked as a reference signal.  The resulting behavior is again safe, but is less conservative due to the Dubins' car being a better representation of the movement of the quadruped in the plane, i.e., a better reduced order model produces less conservative behavior while still maintaining safety. 

\subsection{Manipulators}
Consider a robot manipulator, as illustrated in Fig. \ref{fig:manipulator}.  The control task is to achieve collision free behavior between the robot and environment while accomplishing a task (in this case, flipping a burger).  
Importantly, there is no access to the proprietary onboard real-time controllers of the commercial robot arm, and therefore safety \emph{must be achieved through a LCA}.  

Let $A(q) \subset \R^3$ be the set of all points on the robot (which depends on the configuration of the robot $q \in \R^n$) and $B \subset \R^3$ be the set of all points in the environment.  Collision free behavior between the robot and environment, captured by $A(q) \cap B = \emptyset$ or $A(q) \subset \overline{B}$ with $\overline{B}$ the complement of $B$, is encoded by a barrier function $\mathcal{S} = \{ q \in \R^n : \mathrm{sd}_{AB}(q)\geq 0\}$ defined in terms of the \emph{signed distance} (\cite{schulman2014motion}): 
\begin{sequation}
\begin{aligned}
h(q)  & =  \mathrm{sd}_{AB}(q)   := 
\underbrace{\inf_{\substack{p_A \in A(q) \\p_B \in B}} \norm{p_A-p_B}_2}_{\mathrm{distance}(A(q),B)} - 
\underbrace{\inf_{\substack{p_A \in A(q) \\p_B \in \overline{B}}} \norm{p_A-p_B}_2}_{\mathrm{penetration}(A(q),B)} 
\end{aligned}  
\end{sequation}
The advantage of using the signed distance, as opposed to the distance, is that the addition of the ``penetration'' term which gives a negative value when this occurs---as opposed to the distance which is strictly non-negative.  This negative value allows for convergence back to the safe set $\mathcal{S}$ per the fact that CBFs render $\mathcal{S}$ attractive.  

The challenge with using the signed distance as a barrier function is that it is discontinuous on a set of measure zero \cite{sakai1996riemannian}.  
To accommodate for the discontinuities, consider: 
\begin{sequation}
\label{eqn:sddecomp}
\frac{\partial h}{\partial q} = \frac{\partial  \mathrm{sd}_{AB}}{\partial q} =  \frac{\partial \mathrm{sd}_{AB}^{C^1}}{\partial q} + \delta(q)  
\end{sequation}
which decomposes $\mathrm{sd}_{AB}(q)$ into its differentiable and non-differentiable component,
where the gradient of the non-differentiable component, $\delta$, is viewed as a disturbance that is non-smooth on a set of measure zero; as a result, we can design a controller that is robust to adversarial disturbances of magnitude matching the essential supremum $\| \delta \|_{\infty} = \esssup_{t \geq 0} \| \delta(q_t) \|$.

\hfill \emph{(continued on next page)}
   \end{sidebar}

   \begin{sidebar}{\continuesidebar}
   \renewcommand{\thesequation}{S\arabic{equation}}
   \renewcommand{\thestable}{S\arabic{stable}}
   \renewcommand{\thesfigure}{S\arabic{sfigure}}

\subsection{Continuous Time LCAs in Practice}
For the reduced order model, we consider a kinematic model of the robot arm \eqref{eqn:singleint}, i.e., $\dot{q}_t = v_t$ with $q_t \in \R^n$ for $n$ the number of degrees of freedom (in this case, $n = 6$).  To enforce a safety filter on the reduced order model, the goal is to leverage a QP of the form \eqref{eqn:QPh}.  Yet in this case, due to the fact that the signed distance is not continuously differentiable, we leverage the decomposition in \eqref{eqn:sddecomp} to obtain the QP: 
\begin{sequation}\label{eqn:QPmanipulator}
\begin{aligned}
 v_{safe}(q,t) = \argmin_{v \in \R^{n}} & \| v - v^d(q) \|^2  \\
 \text{s.t.} \quad  & 
\frac{\partial \mathrm{sd}_{AB}^{C^1}}{\partial q} v \geq - \alpha (\mathrm{sd}_{AB}(q)) + \| \delta \|_{\infty} \dot{q}_{\max}  
\end{aligned}
\end{sequation}
with $\dot{q}_{\max} = \| \dot{q} \|_{\infty}$, and $\| \delta \|_{\infty}$ defined as above.   
Here, $v_d(q)$ is obtained from a series of preplanned trajectories that must be executed while avoiding collisions, i.e., $v_{d}(q) = K_P(q^i_{d}-q)$ with $q^i_{d}$ the next waypoint (in time) of the preplanned trajectory.  

The QP in \eqref{eqn:QPmanipulator} was implemented experimentally on FANUC robotic manipulator in a kitchen scenario \cite{singletary2022safety}, i.e., the robot was required to do a variety of cooking related tasks while avoiding collisions with the environment.  As illustrated in Fig. \ref{fig:manipulator}, the robot was able to perform a variety of complex tasks while maintaining safety $h(q) = \mathrm{sd}_{AB}(q) \geq 0$.

\sdbarfig{\label{fig:manipulator}
    \includegraphics[width=1\columnwidth]{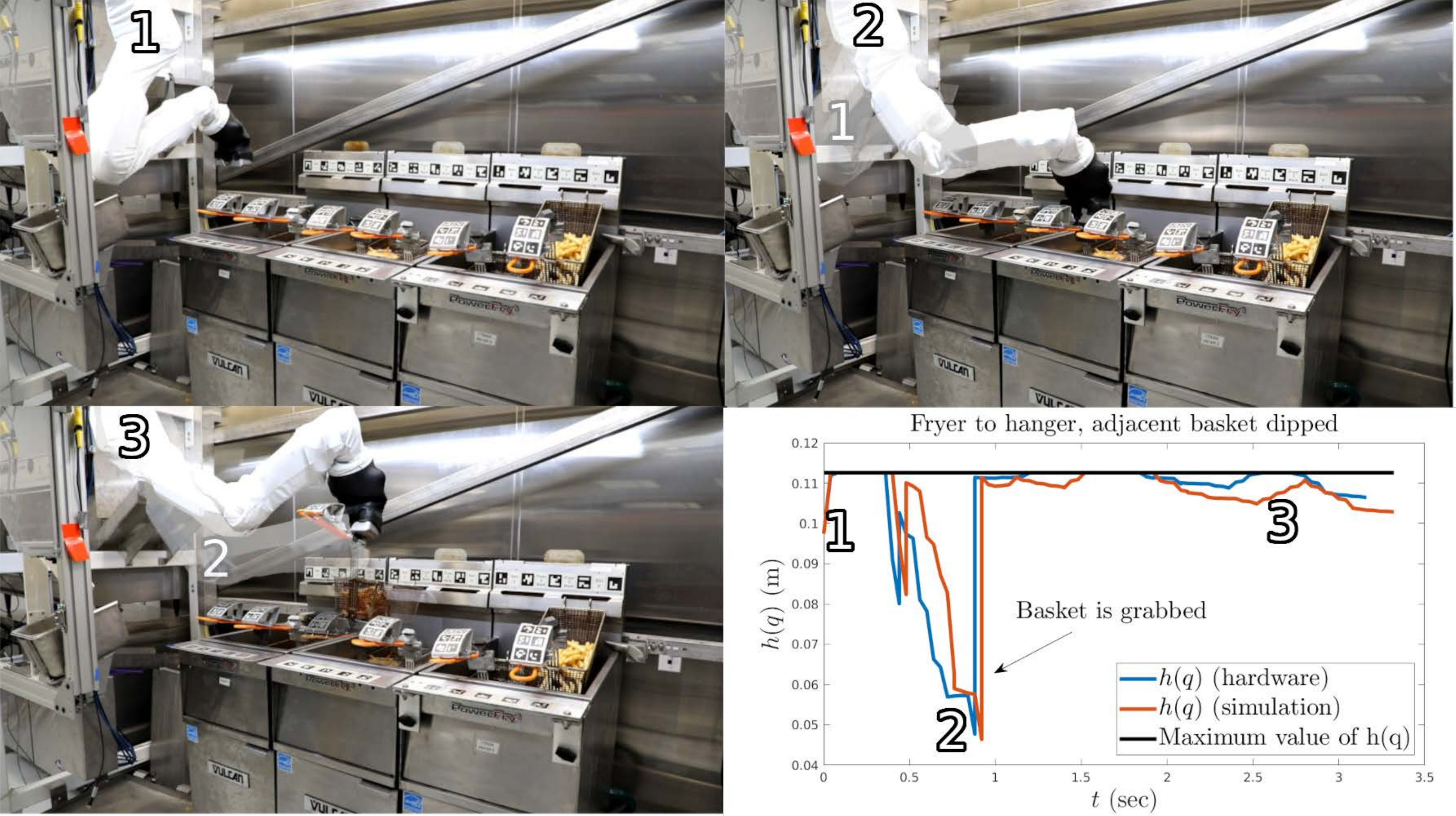}
    }{\footnotesize Achieving safety on a robot manipulator (from \cite{singletary2022safety}).  The manipulator executes a series of preplanned trajectories, and a safety filter is intantiated via a reduced order model to prevent collisions with the environment.  The value of the barrier function is shown, wherein non-negative values imply collision free behavior. }

\vspace{-0.6cm} 
\subsection{Automotive Systems}  


For complex real-world applications, domain specific reduced models are needed.  Additionally, as in the application to manipulators, real-world settings also require extended notions of safety to account for differences between the reduced and full order dynamics. 

To provide an example of this, consider adaptive cruise control (ACC) where the control objective is to achieve a desired speed subject to maintaining a safe distance from a lead car.   In this setting, consider a reduced order model~\eqref{eqn:from} defined by a point-mass model of a vehicle moving in a straight line:  
\begin{no-sequation}
\label{eqn:fgdynamics}
\begin{aligned}
\dot{y} & =  \underbrace{\left[ \begin{array}{c} y_2 \\ -\frac{1}{m}F_r(y) \end{array}\right]}_{f_{rom}(y)} + \underbrace{\left[ \begin{array}{c} 0 \\ \frac{1}{m}\end{array}\right]}_{g_{rom}(y)} u , \qquad  F_r(y) = c_0 + c_1 y_2 + c_2 y_2^2 \\
\dot{z} & = v_0 - y_2, 
\end{aligned} \nonumber
\end{no-sequation}
where  $y_1$ (in $m$) is the position, $y_2 = \dot{y}_1$ (in $m/s$) the velocity, $m$ is the mass of the car (in $kg$), the input $u$ (in Newtons) represents the wheel force, $F_w$, and $F_r$ is the rolling resistance; typically, $c_0$, $c_1$ and $c_2$ are determined empirically.  Finally, $z$ is the distance between the vehicles, wherein it is assumed that the lead vehicle is traveling at a constant speed $v_0$.

The key safety constraint is {\it ``keep a safe distance from the car in front of you.''}  This is generally encoded by the ``half the speedometer'' rule which states that $D \geq \frac{v}{2}$ (with $D$ in $m$ and $v$ in $km/hr$), i.e., the distance between the two vehicles should be at least half the current speed.  Converting this to $m$ and $s$ results in the safety constraint, $z \geq 1.8 y_2$, which can be translated to a barrier function $h(y,z) = z - 1.8 y_2 \geq 0$.  It is easy to verify that this is a valid CBF and can be implemented in practice \cite{mehra2015adaptive}, but we will consider the generalization: 
$$
h(y,z) = z - (a_0 + a_1 y_2 + a_2 v_0 + a_3 y_2^2 + a_4 y_2 v_0 + a_5 v_0^2).
$$
for which the parameters, $a_i$, can be determined such that $z \geq 1.8 y_2$ is satisfied while allowing for actuation limits and other practical considerations to be enforced. 

Let $v^d(y)$ be the ``nominal'' ACC system, i.e., the current algorithm on the vehicle, that drives the velocity $y_2  \to v_g$.  We can then instantiate a safety filter in the form of a QP: 
\begin{sequation}
\label{eqn:QPhISSf}
\begin{aligned}
 v_{safe}(y) = &  \argmin_{v \in \R^{s}} \| v - v^d(y) \|^2  \\
 & \text{s.t.} \quad   
\dot{h}(y,v) \geq - \alpha h(y)  + \frac{\| \frac{\partial h}{\partial y} g_{rom}(y) \|^2}{\epsilon(h(y))},  \quad
  \frac{\partial \epsilon}{\partial h } \geq 0.
\end{aligned}
\end{sequation}
The added term $\epsilon(h(y))$ is a ``tunable'' term that enforces a generalization of input-to-state safety \cite{kolathaya2018input} termed tunable input-to-state safety \cite{alan2021safe}.  Here $\epsilon$ is a function that can be tuned and must have a positive derivative; we pick $\epsilon(h(y)) = \epsilon_0 e^{\beta h(y)}$.  
The safety filter \eqref{eqn:QPhISSf} was implemented on a class-8 truck without a trailer \cite{alan2022control}.  As shown in Fig. \ref{fig:truck}, the nominal ACC controller $v^d$ results in a safety violation and, in fact, a collision.  Using the safety filter on this nominal controller results in safe system behavior.

\vspace{-0.4cm}
\sdbarfig{\label{fig:truck}
    \includegraphics[width=1\columnwidth]{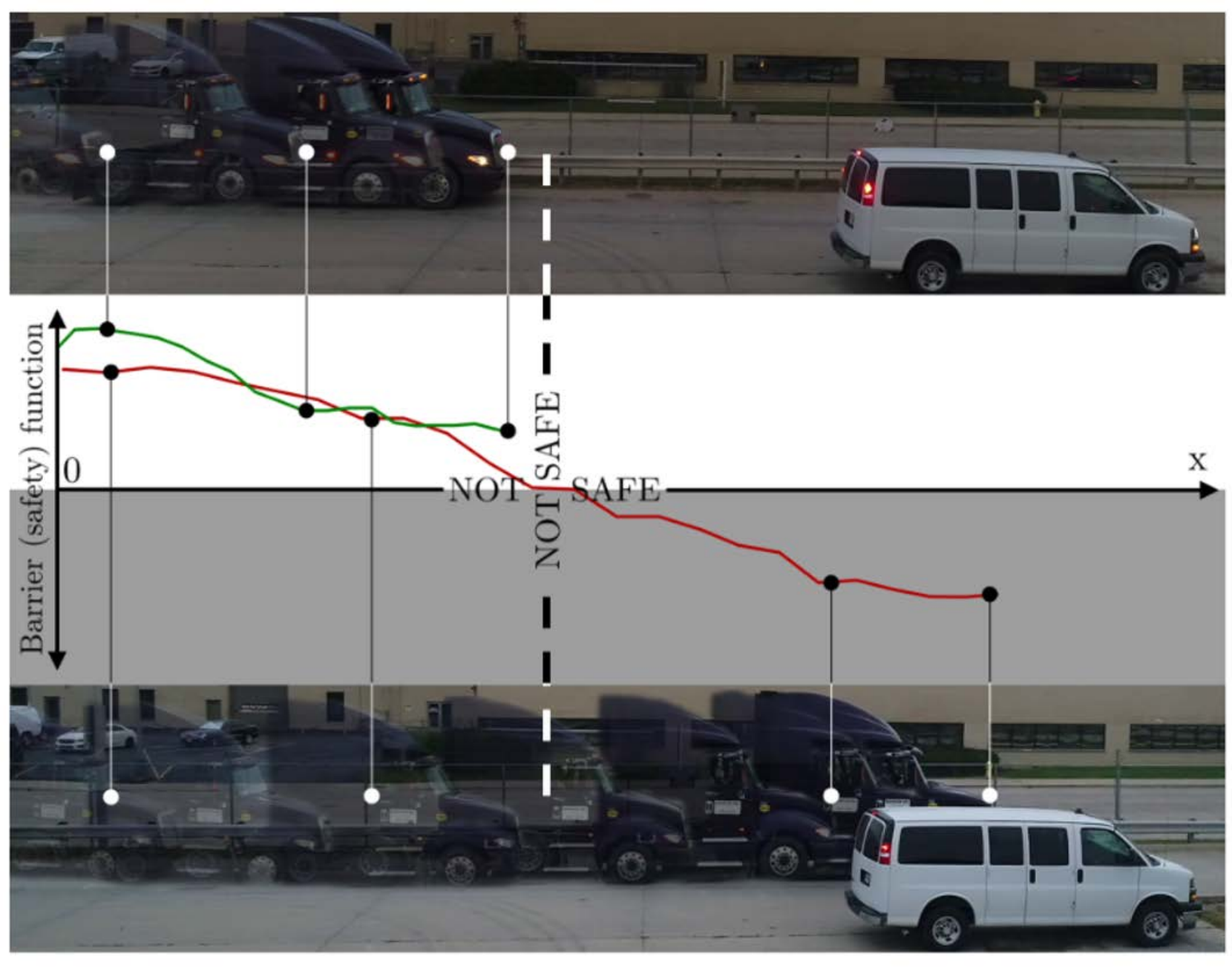}
    }{\footnotesize Safety filter implemented on a full-scale truck \cite{alan2022control}.  Shown is the nominal controller $v^d$ which violates the safety condition (bottom).  When the safety filter \eqref{eqn:QPhISSf} is implemented safety is achieved (top)}
    \vspace{-0.4cm} 
\end{sidebar}

\begin{my_exam}[Running example: robot navigation]
\label{ex:cbfdubins}
Returning to the running example, we can view the Dubins' car as a reduced order model used to enforce safety constraints on a complex mobile robot, e.g., a quadruped.  Recall that the Dubins' car dynamics~\eqref{eqn:unicycle} take the form: $\dot{q}_t = f_{rom}(q_t) u_t$ where $q = (x_1,x_2,\theta)$ and $u = (u_1,u_2)$.  Consider a barrier function defined on the Dubins' car dynamics aimed at avoiding collisions with obstacles:
$$
h(q) = d_0 - r - \kappa \cos(\theta - \theta_0)
$$
where $d_0 = \| (x_1-x^0_1,x_2-x_2^0) \|$ with $(x_1^0,x_2^0)$ the position of the obstacle, $\theta_0 = \arctan((x_2^0-x_2)/(x_1^0 - x_1))$, and $\kappa > 0$ a tunable parameter. 

Let $u^d(q)$ be the feedback controller~\eqref{eq:dubins-feedback} synthesized in the previous section for tracking a nominal trajectory. 
Then the safety filter~\eqref{eqn:QPh} yields a QP on the Dubins' car: 
\begin{align}
\label{eqn:QPhmu}
 u_{safe}(q) = \argmin_{u \in \R^{2}} & \| u - u^d(q) \|^2_{\Gamma}  \\
 \text{s.t.} \quad  & 
\frac{\partial h}{\partial q} f_{rom}(q) u \geq - \alpha h(q)  \nonumber
\end{align}
where $\|u\|^2 = u^T \Gamma u$ with $\Gamma = \mathrm{diag}(1,R)$ where $R > 0$ is a control cost parameter. 
The result is a reference velocity $u_{safe}(q) = (u_{1,safe}(q), u_{2,safe}(q))$ on the forward velocity and change in heading.  These reference signals can be sent to a robot with more complex dynamics as if they were joystick commands.  Theorem \ref{thm:kintofull} guarantees safety for the more complex system under the assumption of good tracking.
\end{my_exam}

\subsection{Discussion}

An underlying principle in designing LCAs for robotic systems is to synergistically leverage the strengths at each layer.  For example, the lowest layer can handle high dimensional nonlinear systems (e.g., via Lyapunov and barrier functions), yet model uncertainty and ''looking ahead'' in nonlinear systems is challenging.   Adding a reference signal generation layer that uses continuous time reduced order models mitigates model uncertainty while still yielding formal guarantees, e.g., on safety.  Adding a discrete planning layer above the reference signal generating layer allows for longer horizon planning, e.g., via MPC with a discrete time linear reduced order model.  Combining these together mitigates the weaknesses at each layer while enjoying their strengths.  This use of diverse models, time-scales, and control approaches was highlighted through experimental demonstration on a wide variety of robotic systems.  We return to the idea of diversity across layers enabling behavior that cannot be achieved by any single layer in \nameref{sec:sweet-spots}, where we introduce a quantitative notion of a \emph{diversity enabled sweet spot} in LCAs.


\paragraph{Final remarks}
The success of LCAs in robotic systems, and the ability to add and remove layers as needed, points to the power of these methods.  It also conveys their complexity---different models at different layers, and the interfacing between these models, results in complex and notationally intensive mathematical models, and makes establishing formal guarantees becomes daunting.  Yet the fact that these approaches work in practice, and are widely understood as the ``way to control robots,'' points to the value in formalizing and analyzing LCAs.   It can be argued that this is a central challenge for the control community moving forward: going beyond homogeneous system models, and analyzing heterogeneous models interacting within an LCA.

\section[Architecture Design as Multi-Criterion Optimization]{Part 2.1: Architecture Design as Multi-Criterion Optimization}
\label{sec:sweet-spots}
The previous sections illustrate how an LCA can be naturally derived from a global decision and control problem, and provide a concrete instantiation of these ideas in the context of robotic systems.  These results highlight both the power of LCAs, as well as the art and complexity involved in designing them.  \edit{We highlight that many idealized assumptions were made in Part 1: we assumed that the control system hardware was already fixed,  that we knew how many layers were needed, what each layer should do, and how layers should interact within the LCA.  In this section, which marks the start of Part 2 of the paper, we try to address some of these idealized assumptions, and propose a  framework rooted in multi-criterion optimization for quantitative reasoning about architecture design choices such those described in the previous two sections.}  A key theme that we explore in this section is that while each layer may be subject to specific constraints and tradeoffs, by leveraging \emph{diversity across layers}, these tradeoffs can be mitigated to yield high-performing LCAs such as those highlighted in the previous sections.

\begin{figure}
    \centering
    \includegraphics[width=0.9\columnwidth]{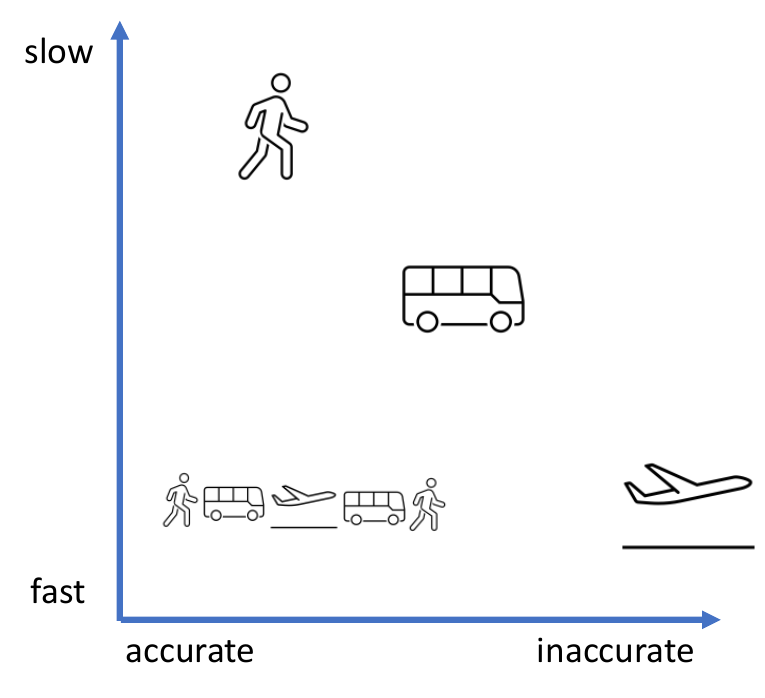}
    \caption{Each individual ``travel layer'' is subject to speed/accuracy tradeoffs, but combining them appropriately in an LCA enables an overall transportation system with minimal tradeoffs in either speed or accuracy.}
    \label{fig:travel-sat}
\end{figure}

We begin with a familiar illustrative example: long-distance travel.  We consider three possible ``travel layers,'' namely air travel (implemented via aircraft and airports), public transit (implemented via busses and bus stops), and walking (implemented via human sensorimotor control).  Each of these travel layers are subject to speed-accuracy tradeoffs, which are themselves a function of architectural design choices (but we will not focus on these here): air travel is fast but inaccurate since we can only fly between airports, public transit is moderately fast and moderately accurate as we are limited to bus stops, and walking is slow but extremely accurate.  These travel layers can be placed in a speed/accuracy plot as shown in Fig.~\ref{fig:travel-sat}.

However, as we all know, when traveling long distances, it is most efficient to appropriately combine these travel layers: we walk to the bus stop, take the bus to the airport, fly to the nearest airport to our destination, take the bus to the stop nearest our destination, and then walk to our destination.  Although not usually thought of in this way, this is an LCA for travel, with air travel serving as a fast but inaccurate layer, public transit serving as an intermediate layer, and walking as a slow but accurate layer.  The resulting LCA, which implements diverse layers using diverse components, is nearly as fast as flying, and just as accurate as walking.  We call such an LCA that leads to minimal tradeoffs between speed and accuracy an architectural \emph{sweet spot}.

It is our claim that such sweet spots are ubiquitously enabled through \emph{diverse layers}\footnote{Diverse layers typically require diverse hardware, or levels.  We discuss levels in more detail in \nameref{sec:LLL}.} being appropriately combined in LCAs.  Indeed, we see comparable diversity in sensorimotor control, robotics, computer networks, and biology, in order to mitigate what appears to be a universal constraint on individual layers, namely that the lower the layer in the ``stack,'' the faster it must operate, but the more limited it is in its capabilities.  Nevertheless, by appropriately combining slow decision making with moderate speed trajectory generation and fast feedback control, we are able to design autonomous systems that are as flexible as the decision making layer, and as accurate and fast as the feedback layer.
In the remainder of this section, we propose a quantitative framework for reasoning about such \emph{Diversity enabled Sweet Spots (DeSS)}.

\subsection{Pareto Surfaces and Pareto Minimax Points}

Our goal is to both characterize the fundamental tradeoffs that different control architectures induce, and to determine whether a control architecture enjoys a (diversity enabled) sweet spot.  To formalize these concepts, we turn to multi-criterion optimization.

\begin{figure}
    \centering
    \includegraphics[width=\columnwidth]{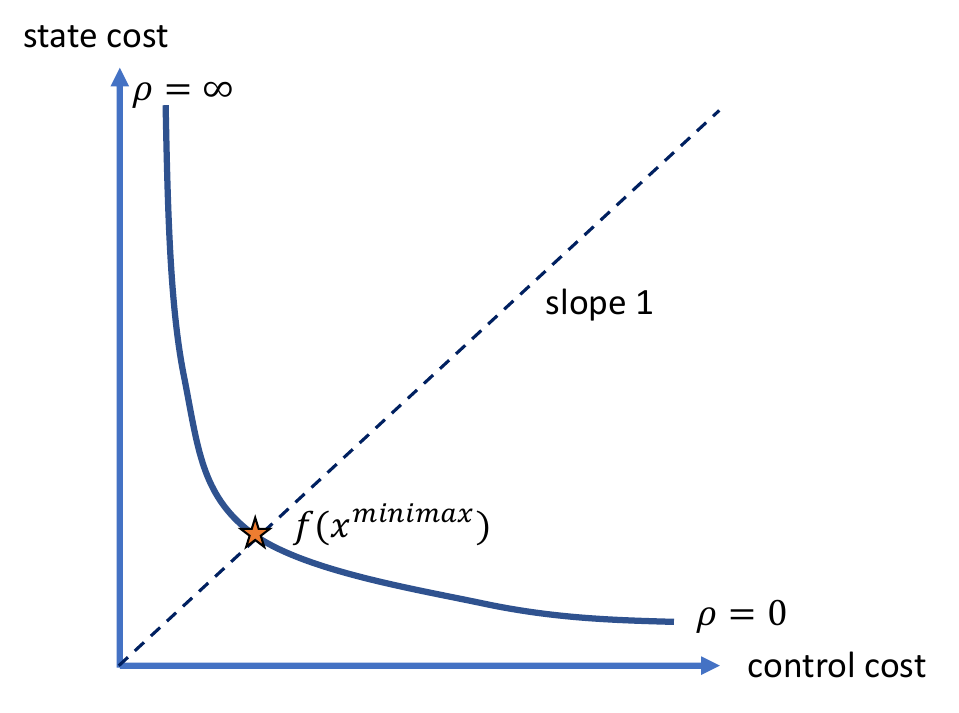}
    \caption{Pareto surface and minimax Pareto optimal point for LQR.  By varying the weight $\rho$ on the control cost, the Pareto surface is traced out.}
    \label{fig:lqr-sat}
\end{figure}

\paragraph{Multi-criterion optimization}  Multi-criterion optimization problems seek to minimize a \emph{vector-valued objective function}. Following \cite[Ch. 4]{boyd2004convex}, we consider a vector-optimization problem which seeks to minimize the vector-valued objective
$$
C(x) = (C_1(x), \dots, C_d(x))
$$
with respect to the positive orthant $\R^d_+$. Such an optimization problem should be interpreted as having $d$ different objectives $C_i$, each of which we would like to make small. 

In contrast to scalar-valued objectives, we must take care in defining appropriate notions of optimality. In particular, we may define both optimal and Pareto optimal points.  A feasible point $x^\star$ is \emph{optimal} if it is \emph{unambiguously better} than any other feasible point, where better is defined in terms of the partial order induced by the positive orthant, i.e., a feasible $x^\star$ is optimal
 if for any other feasible $y$, $C(x^\star)\preceq C(y)$, i.e., if $C_i(x^\star)\leq C_i(y)$ for all $i=1,\dots,d$.  Most engineering design problems are subject to fundamental tradeoffs between optimization criteria $C_i$, and such an optimal point typically does not exist.  Instead, a family of \emph{Pareto optimal} points can be defined, wherein a feasible point $x^{po}$ is Pareto optimal if for any feasible $y$, if $C(y) \preceq C(x^{po})$, then $C(y) = C(x^{po})$, i.e., a feasible point $x^{po}$ is Pareto optimal if no other point exists that is unambiguously better.  Indeed, the existence of multiple Pareto optimal points imply that there is a fundamental tradeoff between the different objectives.

 The standard approach to solving such a multi-criterion optimization problem is via \emph{scalarization}. A common approach to scalarization is to take a weighted sum of the objectives, i.e., for $\lambda \in \R^d_{++}$, define the scalarized objective $C_\lambda(x)=\lambda^TC(x) = \sum_{i=1}^d \lambda_i C_i(x)$.  By sweeping over weighting parameters $\lambda \succ 0$, we obtain a family of Pareto optimal points $x^{po}(\lambda)$, which in turn defines a Pareto surface $(C_1(x^{po}(\lambda)), \dots, C_d(x^{po}(\lambda))) \subset \R^d$.\footnote{Up to boundary points, such an approach is guaranteed to recover all Pareto optimal points if the objective functions $C_i$ are convex in $x$, see~\cite[Ch. 4]{boyd2004convex}.}  An alternative, but also important, scalarization approach is to consider minimizing the the maximum of the objectives, i.e., $C_{\max}(x) = \max\{C_1(x),\dots, C_q(x)\}.$  The resulting solution $x^{mm}$ is called the minimax Pareto optimal point.

A familiar example of bi-criterion optimization in control is LQR optimal control.  Indeed, defining the vector-valued objective $(\sum_{k=0}^N \|x(k)\|_2^2, \sum_{k=0}^{N-1}\|u(k)\|_2^2)$, we recognize the LQR objective $\sum_{k=0}^{N-1} \|x(k)\|_2^2 + \rho \|u(k)\|_2^2 + \|x(N)\|_2^2$ as a scalarization of the competing objectives of small state and control cost.  An alternative, albeit less common, scalarization would be to consider the maximum objective $\max\{\sum_{k=0}^N \|x(k)\|_2^2, \sum_{k=0}^{N-1}\|u(k)\|_2^2\}$.  See Fig.~\ref{fig:lqr-sat} for an example of a typical Pareto curve for an LQR problem.

\subsection{Sweet Spots are Nearly Optimal Points}

\begin{figure*}[t]
    \centering
    \includegraphics[width=0.8\textwidth]{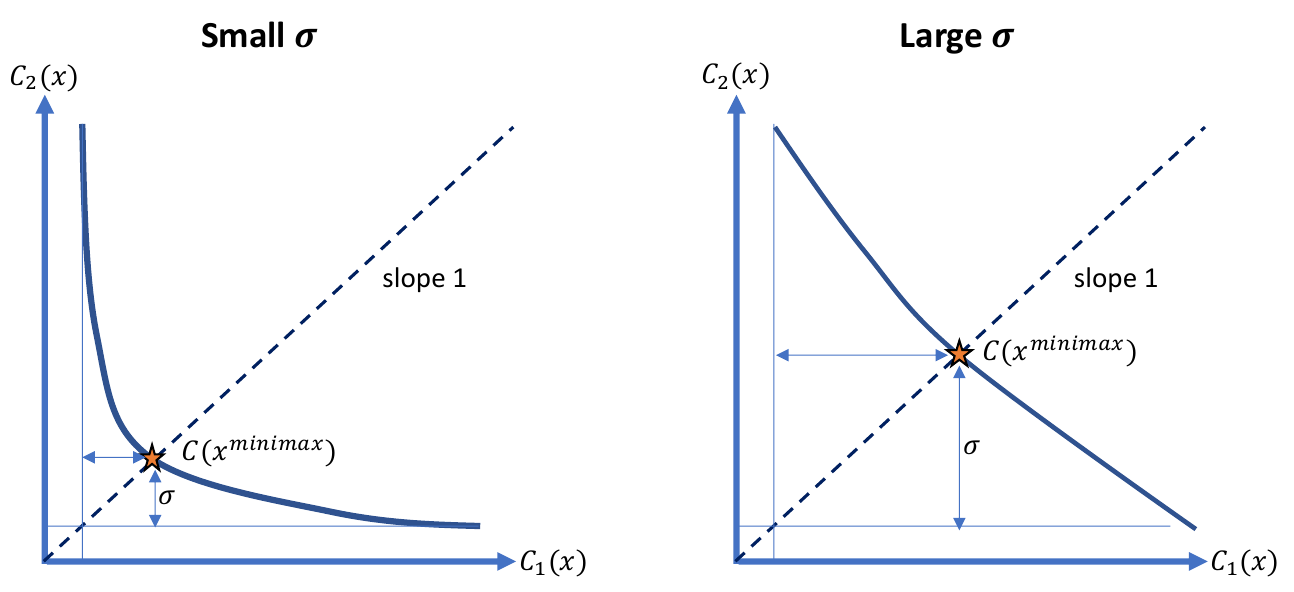}
    \caption{The measure $\sigma$ of a minimax Pareto point quantifies how much of a tradeoff there is between competing objectives by comparing the minimax Pareto optimal point to all other Pareto optimal points.  A small $\sigma$ indicates that there exists a point that is \emph{nearly optimal}.}
    \label{fig:sig-sss}
\end{figure*}

We now have the required concepts to formally define a \emph{sweet spot}.  Intuitively, a sweet spot is a point on the Pareto surface that is \emph{nearly optimal}.  We quantify this notion of near optimality by defining a \emph{$\sigma$-sweet-spot} to be a minimax Pareto optimal point that is $\sigma$ away from being an optimal point in the following sense:
\begin{multline}
    \sigma := \sup_{\lambda \succ 0} \max\{C_1(x^{mm})-C_1(x^{po}(\lambda)), \dots, \\  C_d(x^{mm})-C_d(x^{po}(\lambda))\}.
    \label{eq:sigma}
\end{multline}
In words, the measure $
\sigma$ characterizes the biggest loss in optimality \emph{in any of the criterion $C_i$} of a minimax Pareto optimal point relative to \emph{any other Pareto optimal point}.  Note that if there exists an optimal point $x^\star$ then $\sigma = 0$, and that $\sigma$ increases as the tradeoff between objectives becomes more severe.  See Fig.~\ref{fig:sig-sss} for a qualitative illustration of when $\sigma$ is small or large as a function of the geometry of the Pareto surface.

\subsection{Diversity Enables $\sigma$-Sweet-Spots}
One of our key claims, which is broadly supported by examples in engineering, science, and biology, is that diversity enables nearly optimal sweet-spots, despite individual layers being subject to strict and at time severe tradeoffs.  We begin with a simple stylized example for which the sub-optimality measure $\sigma$ can be computed exactly.  We then explore \nameref{sec:neuro} 
in the next section.

\paragraph{Illustrative Example: Bi-Criterion Least-Squares} We study the bi-criterion least-squares problem
\begin{equation}\label{eq:bicriterion-ls}
    \mathrm{minimize}_x \text{ (w.r.t. $\R^2_+)$} \ (\|A_1 x - b_1\|_2^2, \|A_2 x - b_2\|_2^2)
\end{equation}
through the lens of DeSS.  We assume that $b_1, b_2 \in \R^m$ $A_1, A_2 \in \R^{m \times 2m}$, and $x\in \R^{2m}$.  Our stylized architecture design problem is to design the matrices $A_1$ and $A_2$ by selecting their rows, possibly with replacement, from a palette of $2m$ linearly independent rows $\mathcal{V} = \{v_1,\dots, v_{2m}\} \subset \R^{2m}.$  Our goal is to quantify how \emph{diversity in the row-spaces of $A_1$ and $A_2$ affects the resulting $\sigma$-sweet-spot} of the bi-criterion problem.  We begin with some simple observations:
\begin{itemize}
    \item If we assume that the we design $A_1$ and $A_2$ to each respectively have full row rank, then it is clear that each individual objective can be made $0$. 
 We make this assumption going forward, and hence we have $\sigma = \max\{\|A_1 x^{mm} - b_1\|_2^2, \|A_2 x^{mm} - b_2\|_2^2\}$.
   \item If we further assume that the stacked matrix $\bar A = [A_1^T, A_2^T]^T$ has full row rank, i.e., that $A_1$ and $A_2$ do not share any rows selected from $\mathcal V$, then $\sigma = 0$.  This is easily verified by setting $x^{mm} = \bar A^{-1}\bar{b}$, with $\bar{b}=(b_1, b_2).$
\end{itemize}

Thus, our remaining task is to characterize the $\sigma$-sweet-spot for optimization problem~\eqref{eq:bicriterion-ls} when $A_1$ and $A_2$ share a common row space.  Towards that end, we consider the minimax scalarization~\eqref{eq:bicriterion-ls}:
\begin{equation}\label{eq:minimax-ls}
    \begin{array}{rl}
    \mathrm{minimize}_{x} & \max\{\|A_1 x - b_1\|_2^2, \|A_2 x - b_2\|_2^2\},
    \end{array}
\end{equation}
and its dual (see Appendix~\ref{sec-app:dual} for details):
\begin{equation}
\label{eq:minimax-ls-dual}
    \begin{array}{rl}
    \mathrm{maximize}_{\lambda_1, \lambda_2, \mu_1, \mu2} & 2\mu_1^Tb_1 - \frac{\|\mu_1\|_2^2}{\lambda_1} - 2\mu_2^Tb_2 - \frac{\|\mu_2\|_2^2}{\lambda_2}\\
    \text{subject to} & A_1^T \mu_1 = A_2^T \mu_2, \\ & \lambda_1 + \lambda_2 = 1, \, \lambda_1, \, \lambda_2 >0.
    \end{array}
\end{equation}

This allows us to immediately reconfirm that $\sigma=0$ if $A_1$ and $A_2$ do not share any rows, as in this case any dual feasible solution has $\mu_1=\mu_2=0$.  Similarly, when $A_1 = A_2$, a simple argument shows that $\sigma = 1/4\|b_1-b_2\|_2^2$.  A generalization of this argument is presented in the next theorem, proved in the Appendix, which allows us to characterize the solution when $A_1$ and $A_2$ share $k$ rows.
\begin{my_thm}\label{thm:sigma}
    Consider the bi-criterion least-squares problem~\eqref{eq:bicriterion-ls}.  Suppose that $A_1$ and $A_2$ are both full row-rank, and assume without loss of generality, reordering rows in $A_i$ and elements in $b_i$ if necessary, that $A_1$ and $A_2$ share their first $k$ rows.  Then the minimax solution $x^{mm}$ to the scalarized problem~\eqref{eq:minimax-ls} defines a $\left(\frac{1}{4}\|E_k^T(b_1-b_2)\|_2^2\right)$-sweet-spot, as defined in equation~\eqref{eq:sigma}.  Here, $E_k = [e_1, \dots, e_k]$ with $e_i\in\R^m$ the standard basis elements.
\end{my_thm}

\begin{figure*}
    \centering
    ~\hfill~\begin{subfigure}[t]{.4\textwidth}
         \includegraphics[width=.9\textwidth]{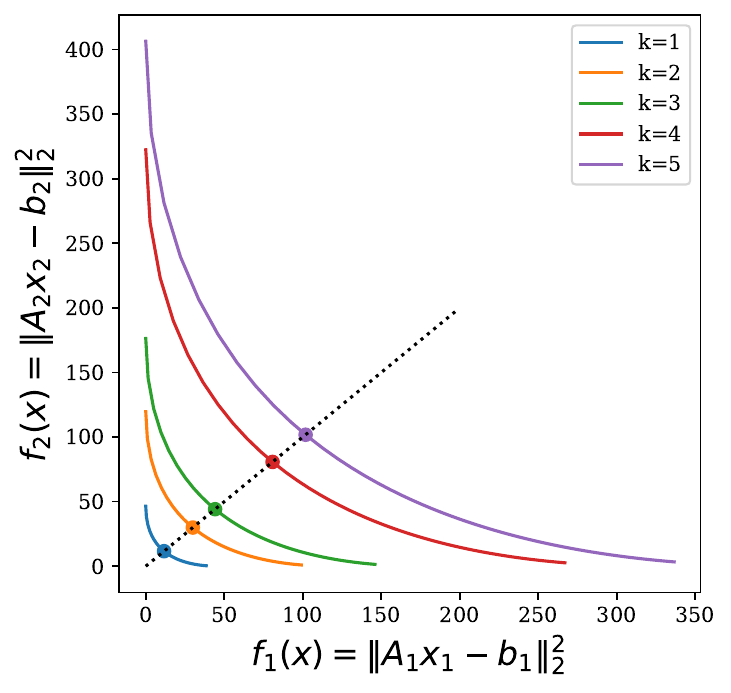}
    \caption{We observe how the Pareto surface for the bi-criterion least-squares problem~\eqref{eq:bicriterion-ls} becomes increasingly unfavorable as we decrease the diversity across $A_1$ and $A_2$.  This is true both in terms of the overall Pareto surface, as well as the sub-optimality measure $\sigma$.}
    \label{fig:ls-dess}
    \end{subfigure}~\hfill~\begin{subfigure}[t]{.4\textwidth}
       \includegraphics[width=.9\textwidth]{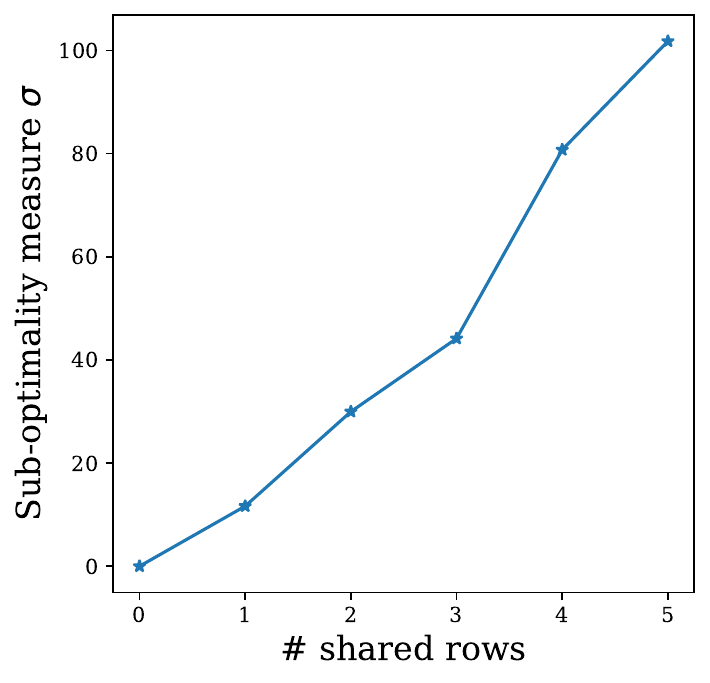}
    \caption{We plot the sub-optimality measure $\sigma$ as a function of shared rows across $A_1$ and $A_2$ for the bi-criterion least-squares problem~\eqref{eq:bicriterion-ls}.  We observe that $\sigma$ deteriorates as we decrease diversity across $A_1$ and $A_2$.}
    \label{fig:sweet-spot-ls}
   \end{subfigure}~\hfill~
   \caption{Increased diversity in row spaces provably leads to less severe tradeoffs, as quantified by a smaller suboptimality measure $\sigma$, in the bi-criterion least-squares problem~\eqref{eq:bicriterion-ls}.}
   \label{fig:bicriterion}
\end{figure*}

Theorem~\ref{thm:sigma} makes clear that the more diverse the matrices $A_1$ and $A_2$, i.e., the smaller the number of shared rows $k$, the less severe the tradeoff; similarly, the less diverse the matrices $A_1$ and $A_2$, i.e., the larger the number of shared rows $k$, the more severe the tradeoff.  We compute a family of the resulting Pareto curves and minimax optimal points for $m=5$ in Fig.~\ref{fig:ls-dess}, and plot the evolution of the sub-optimality measure $\sigma$ as a function of the number of shared rows in Fig.~\ref{fig:sweet-spot-ls}. Parameters are randomly generated so as to ensure the requisite linear independence conditions, and such that $|(b_1-b_2)_i|$ is approximately even for all $i$: details of how the parameters are generated can be found in the Appendix, and the code used to create these plots can be found \href{https://colab.research.google.com/drive/1jK0fJbSzxBb78yHcru8EsRle9V-aZYYx?usp=sharing}{here}.

To further gain insight into the LCA design problem, let us view $A_1$ and $A_2$ as defining two layers, with layer $i$ aimed at addressing control subtask $b_i$.  This analogy reinforces that diversity is not enough to ensure a small sub-optimality measure $\sigma$: the control subtasks, here characterized by $b_1$ and $b_2$, must themselves also be compatible with system diversity (or lack thereof).  For example, even if $k=1$, a very large $(e_1^T(b_{1}-b_{2}))^2$ will nevertheless lead to a severe tradeoff between optimizing the two objectives, resulting in a large $\sigma$.  Conversely, diversity is only needed in $A_1$ and $A_2$ if the control subtasks $b_1$ and $b_2$ are also diverse: if $b_1=b_2$ then $A_1=A_2$ will still yield $\sigma=0$.  Connecting this back to the travel example, if a destination is just a block away, then diversity is not required, and just walking is an optimal travel LCA. 
 Conversely, if the destination is extremely remote, then the three layers of commercial air travel, public transit, and walking will still be very slow. Thus this simple example hints at an explanation as to why diverse layers are needed by systems that must accomplish diverse tasks across diverse environments at diverse spatiotemporal resolutions.  We explore a (still stylized) control problem in the next subsection that further reinforces this concept.



\section[A Case Study in Sensorimotor Control]{Part 2.2: A Case Study in Sensorimotor Control}
\label{sec:neuro}

\newcommand{\eg}{{\it{e.g.}}}
\newcommand{\ie}{{\it{i.e.}}}
 \newcommand{\tm}{\text{max}}
\newcommand{\inr}{\mathbb R}
\newcommand{\inn}{\mathbb N}
\newcommand{\mi}{\mathcal I}
\newcommand{\mj}{\mathcal I}
\newcommand{\eq}{&=&}
\newcommand{\lleq}{&\leq&}
\newcommand{\BEQ}{\begin{eqnarray*}}
\newcommand{\EEQ}{\end{eqnarray*}}
\newcommand{\BEQL}{\begin{eqnarray}}
\newcommand{\EEQL}{\end{eqnarray}}
\newcommand{\BDEF}{\begin{definition}}
\newcommand{\EDEF}{\end{definition}}
\newcommand{\BTHM}{\begin{theorem}}
\newcommand{\ETHM}{\end{theorem}}
\newcommand{\BLEM}{\begin{lemma}}
\newcommand{\ELEM}{\end{lemma}}
\newcommand{\BPF}{\begin{proof}}
\newcommand{\EPF}{\end{proof}}
\newcommand{\BCO}{\begin{corollary}}
\newcommand{\ECO}{\end{corollary}}
\newcommand{\BRM}{\begin{remark}}
\newcommand{\ERM}{\end{remark}}
\newcommand{\lb}{\llbracket}
\newcommand{\rb}{\rrbracket}
\newcommand{\bbm}{\begin{bmatrix}}
\newcommand{\ebm}{\end{bmatrix}}

We adapt the following from~\citet{nakahira2015hard} and ~\citet{Nakahira2021}.  Our goal in this section is to highlight how diverse layers, and the diverse hardware used to implement them, in the human sensorimotor LCA (see Fig.~\ref{fig:axonsize}) enable astonishingly efficient DeSS in spite of severe speed/accuracy tradeoffs.  To that end, we first derive robust performance limits for a simplified model of sensorimotor control subject to communication that is delayed and quantized due to its implementation using physiological hardware composed of axons.  We then identify a simple layered architecture composed of delayed but accurate vision (planning) and fast but inaccurate reflex control (feedback) layers, and show that this architecture is optimal for the aforementioned sensorimotor control model, and leads to a DeSS.  Finally, we show that despite the simplicity of the model and analysis, it is shockingly predictive of real-world behavior as confirmed in~\nameref{sidebar:biking}.

\subsection{A Simplified Model} Consider an initial minimal model with discrete time dynamics

\BEQL
\label{eq:plant}
\begin{array}{rcl}
x(k+1) &=& a x(k) + w(k-T_w) + Q(u(k-T_u))\\
u(k) &=&  K(x(0:k),w(0:k),u(0:k-1)),
\end{array}
\EEQL
where $x(k) \in \inr$ is the state, $w(k) \in \inr$ is the disturbance, $u(k) \in \inr$ is the control action generated by the controller $ K$, and $ Q:\mathbb R\to \mathscr{S}_R$, for $\mathscr S_R \subset \inr$ a finite set of cardinality $2^R$, is a quantizer that limits communication between the controller and the actuator to $R$ bits/sampling interval.  The form of the control law in system \eqref{eq:plant} implies that the controller is \textit{Full Information} (FI), as the control signal $u(k)$ is allowed to depend on all current and past states $x(0:k)$, current and past disturbances $w(0:k)$ and past control actions $u(0:k-1)$.

\begin{figure}[h!]
\centering
\includegraphics[width=\columnwidth]{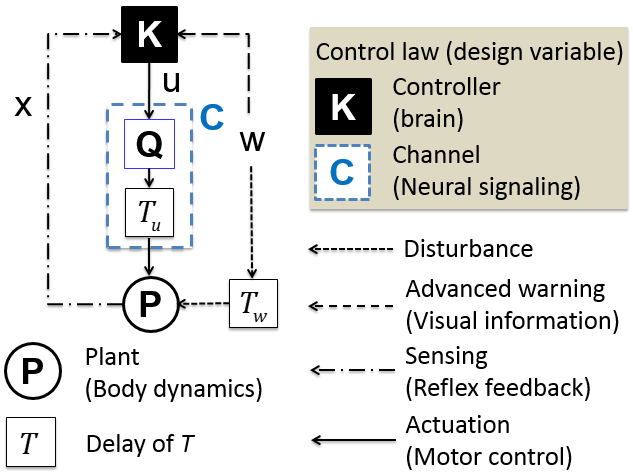}
\caption{Feedback system model for sensorimotor control.}
\label{fig:systemfigure}
\end{figure}
A schematic for this model is shown in Fig. \ref{fig:systemfigure}, where we use $P$ to denote the plant defined by equation \eqref{eq:plant}.  The control signal $\vec u$ is transmitted to the actuator (colocated with the physical plant $P$) via the communication channel $C$, which is defined by the composition of the quantizer $Q$ with the delay block $T_u$. This delay block implies that the controller command $u(k)$ takes $T_u (\geq 0)$ sampling intervals to reach and be executed by the actuators, i.e.,  $u(k)$ only affects the plant $T_u+1$ sampling intervals later.\footnote{We assume that the channel $C$ is memoryless and stationary with rate $R$, allowing us to restrict the quantizer $Q$ to be memoryless and static as well.  Generalizations that lift this assumption can be found in~\cite{nakahira2015hard}.} Note that because $Q$ and the delay block commute the dynamics \eqref{eq:plant} and Fig. \ref{fig:systemfigure} are indeed consistent. We assume $\|w(k)\|_\infty \leq 1$ and $x(0)=0$. The disturbance is known to the controller with an advanced warning of $T_w (\geq 0)$ sampling intervals, i.e., the controller has access to $w(0:k)$ even though the disturbance only affects the plant $T_w +1$ sampling intervals later.

The robust control problem can then be posed as
\BEQL
\label{eq:optimal}
\begin{array}{rl}
\underset{( K, Q) \in \mathcal Q_R}{\text{minimize }} & \sup_{ k\geq 0,\, \| w(k) \|_\infty \leq 1} \| x(k)\|_\infty \\
\text{s.t.} & \text{dynamics \eqref{eq:plant}}
\end{array}
\EEQL
where $\mathcal Q_R$ is the space of control laws defined by the pair of mappings $( K,  Q)$, with $ Q$ constrained to be a static memoryless quantizer of rate $R$, i.e., $ Q:\mathbb R \to \mathscr{S}_R$. This cost function is standard in $L_1$ robust control \cite{dahleh1994control}, except that a communication channel $ C$, composed of a quantizer $ Q$ and a delay $T_u$, is inserted into the feedback loop.  Perhaps surprisingly, this problem formulation still allows for a simple and intuitive analytic solution. 
Indeed, without quantization or delay, the control law $$u(k)=-ax(k)-w(k)$$ ensures that $x(k+1)=0$. Thus any errors in the state is a direct consequences of quantization and/or delay, or to saturation of the control signal $\vec u$.


\subsection{Fundamental Limits due to Delay and Quantization}  In this subsection, we provide an exact solution to the robust control problem \eqref{eq:optimal} for fixed advanced warning $T_w$ and actuation delay $T_u$. In particular, we show that the worst-case state deviation can be expressed as a function of the plant pole $a$, the channel rate $R$, and the \emph{net delay} of the system $T := T_u - T_w$.  The achievable performance takes a different form depending on the net delay regime that the system is operating under. When the net delay $T$ is positive ($T>0$), this corresponds to a system in which the control action $u(k)$ can only affect the plant $T$ sampling intervals \emph{after} the disturbance $w(k)$ affects the state. Conversely, when the net delay  $T$ is non-positive ($T \leq 0$), this corresponds to a system in which there is \emph{advanced warning} of the disturbance, allowing the controller to act in advance. These two qualitatively different cases are treated separately.  We then use these insights in the next section to pose a LCA design problem that seeks to identify an appropriate combination of fast but inaccurate and slow but accurate neural signaling to enable a DeSS.

\begin{my_thm}
    \label{theorem:single_pole}
Suppose that $|a|<2^R$. Then the minimal state deviation achievable in robust control problem \eqref{eq:optimal} is 
\vspace{-0.3mm}
\begin{equation}
\label{eq:single_pole_lemma}
\begin{array}{rl}
  \sum_{i = 1}^{T} |a^{i-1}|+ |a^{T}| \Big(2^R - |a| \Big)^{-1}  & \text{if } T > 0\\ 
     \Big(2^R  - |a| \Big)^{-1} & \text{if }  T \leq 0.
\end{array}
\end{equation}
Conversely, if $|a| \geq 2^R $, then the system cannot be stabilized, and the optimal value to optimization problem \eqref{eq:optimal} is infinite.
\end{my_thm}

The performance limits \eqref{eq:single_pole_lemma} are remarkably simple and intuitive. The net warning case ($T\leq 0$) has only one term due to quantization, with the stabilizability condition $|a|<2^R$ well-known from the networked control system literature \cite{nair2007feedback}. With no dynamics ($a=0$) this reduces to a trivial rate distortion theorem with error $2^{-R}$.
The net delayed case ($T<0$) is more interesting, with the first term due to the delay alone, and the second term an additional contribution due to quantization.  As expected, both grow rapidly with increased net delay $T$ and unstable $a>1$, for reasons familiar and intuitive. 

\subsection{Speed Accuracy Tradeoffs in Neural Signaling} 
\begin{figure}
    \centering
        \includegraphics[width=0.45\textwidth]{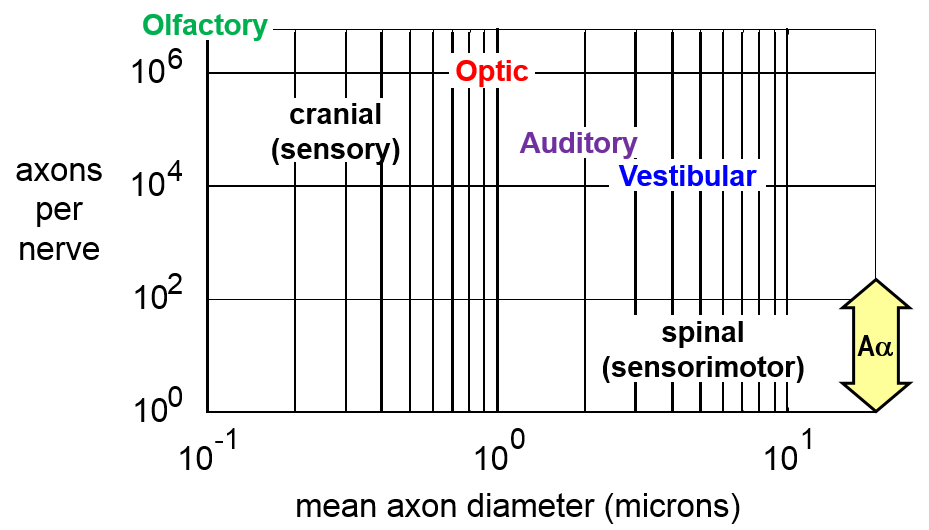}
        \caption{Axons per nerve ($\propto$ resolution) versus mean axon diameter ($\propto$ speed) for four key cranial nerves, and the largest $(A\alpha)$ sensorimotor axon which occurs in spinal and peripheral nerves (in copies from 1 to hundreds). }
        \label{fig:axonsize}
    \end{figure}

We now add a tradeoff between temporal and spatial resolution in neural signaling to our model via the net delay $T$ and data rate $R$. We believe this is the first important constraint in explaining the extreme heterogeneity found in the nervous system, and is analogous to the speed/accuracy tradeoff highlighted in the travel example above. The nervous system communicates between components and the body with a variety of nerves, which are bundles of axons. Axons are the wiring by which spiking neurons communicate long range using action potentials, and it is possible to derive some rough tradeoffs from well-known physiology.  Fig.~\ref{fig:axonsize} shows some of the tremendous diversity of axon numbers and sizes among the cranial and peripheral nerves. We argue that much of this arises due to hard constraints on speed versus accuracy. 

We suppose that our channel $ C$ (see Fig.~\ref{fig:systemfigure}) is a single nerve with uniform signaling delay $T_s$, and assume that the total delay $T_u$ is the sum $T_u = T_s + T_c$ with an additional fixed delay $T_c$ due to grey matter computation and other communications. Initially we assume that $T_c$ is fixed and given, and that $T_s$ is variable and depends on the nerve composition, as in Fig.~\ref{fig:axonsize}. Following the arguments provided in~\cite{nakahira2015hard,Nakahira2021}, we use the physiologically plausible yet remarkably simple relationship between data rate $R$ and signalling delay $T_s$:
\begin{equation}
\label{eq:delay_rate}
R = \lambda_\alpha T_s
\end{equation}
where $\lambda_\alpha $ is a resource measure that scales with the axon area $\alpha$. 

Next we explore the surprisingly rich consequences of the constraint $R= \lambda_\alpha T_s$ on our minimal model of sensorimotor control using Theorem~\ref{theorem:single_pole}. For simplicity, we write $\lambda$ from now on as the resource dependence is understood.  One can verify that if $R=\lambda T_s$ and $T_u := T_s + T_c$, then the optimal optimal performance specified in Theorem~\ref{theorem:single_pole} becomes 
\begin{equation*}
\begin{array}{rl}
  \sum_{i = 1}^{T} |a^{i-1}|+ |a^{T}| \Big(2^{\lambda T_s} - |a| \Big)^{-1}  & \text{if $T>0$}\\ 
     \Big(2^{\lambda T_s}  - |a| \Big)^{-1} &  \text{if $T\leq 0$}.
   \end{array}
\end{equation*}

\begin{figure}
 \centering \includegraphics[width=0.6\columnwidth]%
    {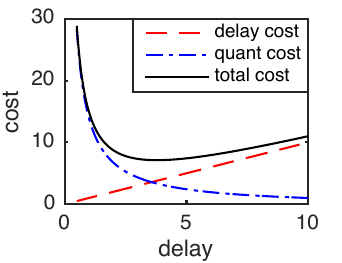}
    \caption{Impact of speed versus accuracy. The cost of delay $\sup_{\|w\|_\infty \leq 1} \|x_d\|_\infty$, the cost of quantization $\sup_{\|w\|_\infty \leq 1} \|x_q \|_\infty$, and the total cost $\sup_{\|w\|_\infty \leq 1} \| x_d\|_\infty$ is shown with varying delay $T_s$ when $\lambda = .1$, $a = 1$. }
    \label{fig:kradeoff}
\end{figure}

Fig.~\ref{fig:kradeoff} shows the system performance when varying delay $T_s$ (and thus channel rate $R$) for $T_c=T_w=0$ and a fixed resource level $\alpha$. Increased delay increases the delay error term $\sup_{\|w\|_\infty \leq 1} \|x_d\|_\infty:=\sum_{i = 1}^{T} |a^{i-1}|$ but reduces the quantization error term $\sup_{\| w\|_\infty \leq 1} \| x_q\|_\infty:=\Big(2^{\lambda T_s} - |a| \Big)^{-1}$. Consequently, the optimal system level performance is achieved at intermediate levels of delay and channel rate. Because of the exponential dependence there is no analytic formula for the optimum, but the error is convex and the minimum easily found numerically.  Next we consider in more detail the consequences of these formulas by varying the additional delays and plotting the resulting optimal errors, bits, and delay.

\begin{figure}[h!]
\centering
\includegraphics[width=\columnwidth]{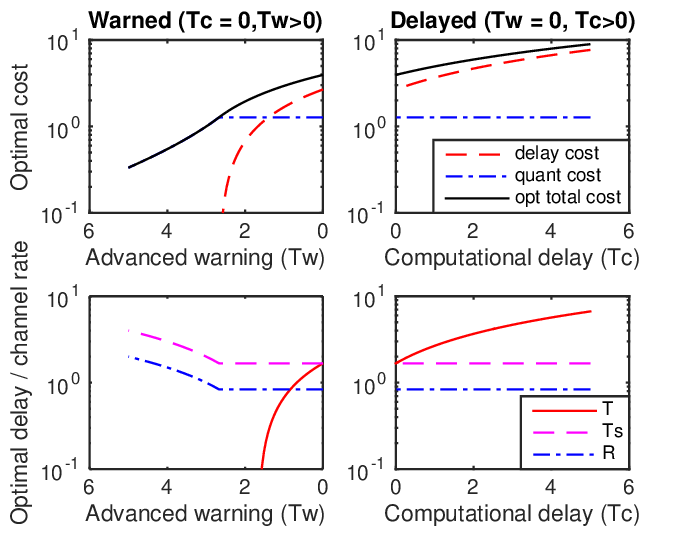}
\caption{Delayed versus warned system.}
\label{fig:phasetransition}
\end{figure}

Fig.~\ref{fig:phasetransition} shows the optimal delays $T_s$ (and resulting net delay $T$) and channel rate $R=\lambda T_s$ that achieves the minimum total error when varying $T_w \ge 0$ and  $T_c \ge 0$ separately in the two special cases (i) $T = T_u- T_w \leq 0$ (warned)  and (ii) $T = T_s +T_c > 0$ (delayed).  What results are clearly two distinct regimes with distinct physiology.  When the computation delay $T_c$ is greater than  $0$, the system has a net delay $T$ and the delay cost increasingly dominates the total cost, leading to both the data rate $R$ and signaling delay $T_s$ becoming constant (i.e., suggesting axons of a large and constant radius $\rho$), independent of $T_c$.  This corresponds to the reflexes on the right half of Fig.~\ref{fig:axonsize} with nerves having relatively few large axons---these are the physiological analogs to aircraft and airports from our travel example. The total error, due mostly to delay, can be much larger than the disturbance.  Concretely, in running or cycling on rough terrain or through heavy traffic, a relatively small but well placed perturbation to the foot or wheel can be amplified into a crash, even a fatal one---this effect gets worse at high speeds when the delay is relatively larger. Our nervous system invests in large nerves, axons, and muscles to avoid such crashes, consistent with the theory. 

With increasing advanced warning $T_w > 0$ the net delay $T$ becomes non-positive, and in this case the errors due to quantization increasingly dominate the total cost.  Further, this total cost goes to zero as $T_w$ increases, exactly the opposite of the delayed case.  Further, as the advanced warning $T_w$ increases, so does the data rate $R$, and consequently the axon radius $\rho$ decreases (as $\alpha\approx \pi R \rho^2$ is fixed).  This corresponds to the left side of Fig.~\ref{fig:axonsize} with many relatively small axons---these are the physiological analogs to walking in our travel example. In running or cycling we can start with huge errors to remotely located objects, and given enough time drive them to zero.  Here we are limited largely by the resolution of our vision in accurately locating the object, again consistent with the theory.

Thus we have an extremely simple model that connects the high layer requirements of advanced warning and planning (e.g. as enabled by vision) to the low layer control implemented by fast reflexes. 
In the sequel we explore further aspects of this model, and introduce additional constraints and generalizations.

\subsection{A Minimal LCA}

\begin{figure}
\includegraphics[width=.9\columnwidth]{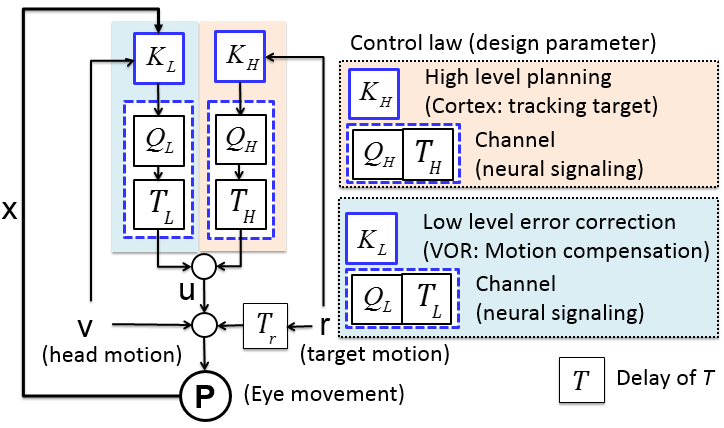}

\caption{Modeling visual processes as a LCA.
}
\label{fig:visualmodel}
\end{figure}

One of the most important features of a visual system is its distributed nature, in which sensors, actuators, and computational components are interconnected via sparse communication.  Fig.~\ref{fig:visualmodel} sketches a minimal model of this kind that is composed of two copies of each component in Fig.~\ref{fig:systemfigure}. The plant dynamics are given by $x(k+1) = a x(k) + u(k) + w(k)$ except the disturbance is now composed of two terms $w(k) = v(k) + r(k-T_r)$, as is the control action $u(k) =u_L(k-T_L)+u_H(k-T_H)$, each generated by their own sensors, computing, and communication components. Visual trajectory planning is done through the control loop involving $ Q_H$ which is responsible for tracking, via the control signal $u_H(k)$, a visual target whose change in position is captured by $\vec r$. We assume a very simplified view of vision whereby remote (in space) sensing means that $r(k)$ is seen but it takes $T_r$ for the disturbance to arrive, effectively creating an advanced warning of $T_r$, though the physical details are all causal.

On the other hand, local (reflex) compensation is done through the control loop involving $ Q_L$. Disturbances such as those caused by body and head motion are captured by $\vec v$, and are sensed directly by the Vestibular Occular Reflex (VOR), which computes a control action $u_L(k)$ to compensate. The control commands $(u_H(k), u_L(k))$ from both loops are sent to the plant through different signaling pathways, modeled by channels with rates $R_H$ and $R_L$ and delays $T_H$ and $T_L$, respectively, after which their gains are summed to produced the final previously described control action $u(k) = u_L(k-T_L)+u_H(k-T_H)$. Connecting this LCA back to the formalism introduced in \nameref{sec:layering} and \nameref{sec:robotics}, we immediately recognize $u_H$ as a \emph{feedforward} control term computed at the \emph{planning layer} which provides advanced warning of the coming reference position $\vec r$, and $u_L$ as a \emph{feedback} control term, computed at the \emph{feedback control layer} and executing in near real-time to compensate for unforeseen disturbances $\vec v$.

Using the tradeoff \eqref{eq:delay_rate} in both signaling pathways, and bounding $\|\vec v\|_\infty$ and $\|\vec r\|_\infty$ from above by $1$ and $\delta$, respectively, the optimal performance is then given by
\BEQ
\left\{ \sum_{i = 1}^{T_L} |a^{i-1}|+ |a^{T_L}| \Big(2^{R_L} - |a| \Big)^{-1}  \right\}
+   \delta \Big(2^{R_H}  - |a| \Big)^{-1}. 
\EEQ 
This result follows by noting that the total system can be decomposed into two independent subsystems, corresponding to the $ Q_H$ and $ Q_L$ loops, and thus so can its performance. The first subsystem is a delayed system driven by $\vec v$ and controlled by $\vec u_L$, while the second subsystem is a warned system driven by $\vec r$ and controlled by $\vec u_H$. From our previous analysis, it is expected that the first system achieves better performance when its nerves are composed of a few large and fast axons, whereas the second system achieves better performance when its nerves are composed of many small and slow axons. This phenomena can be indeed observed in the real visual systems \cite{sterling2015principles}. Specifically, the optic nerve has approximately 1M axons of mean diameter $0.64 \mu m$ with CV $0.46 \mu m$, while the 20K vestibular axons have mean diameter $2.88 \mu m$ with CV $0.41$, significantly larger and less numerous and slightly less variable.

We conclude by emphasizing that a key enabler for DeSS is diversity in hardware used to implement diverse layers to address diverse system tasks.  For example, in the biking example discussed in \nameref{sidebar:biking}, if the trail planning layer had to update the nominal trajectory faster than vision could handle, the LCA would fail to enable a DeSS.  It is this \emph{multi-rate} nature of control tasks, characterized by local fast corrections and global slow updates, and which seems to be ubiquitous across engineered and natural complex systems, that allows for corresponding multi-rate LCAs to be designed that enable DeSS.  Developing a general quantitative design framework for multi-rate LCAs that enable DeSS is arguably the most important open problem in engineering today, and one that control theorists are particularly well-suited to tackle.

\begin{sidebar}{Experimental Validation in a Biking Simulator}
\section[Experimental Validation in a Biking Simulator]{}\phantomsection
\label{sidebar:biking}
\setcounter{sequation}{0}
\renewcommand{\thesequation}{S\arabic{sequation}}
\setcounter{stable}{0}
\renewcommand{\thestable}{S\arabic{stable}}
\setcounter{sfigure}{0}
\renewcommand{\thesfigure}{S\arabic{sfigure}}

Sensorimotor control was studied in the context of the multisensory task of mountain bike riding using a video game as the experimental platform \cite{Nakahira2021}. The game captures tunable requirements on player performance which require layered architectures in the nervous system to create DeSS due to the constraints imposed by physiology, see Fig.~\ref{fig:visualmodel}.  Naively, success in the biking task seems to require speed and accuracy that the raw hardware lacks, making non-layered solutions infeasible. The layered nervous system breaks the overall biking problem into a high \textit{trails} (trajectory planning) layer of slow but accurate vision with trail look-ahead for advanced warning, and a low \textit{bumps} (feedback control) layer that uses fast but inaccurate muscle spindles and proprioception to sense and reject bump disturbances. The motor commands from these two control loops to the muscles simply add in the optimal case, as well as in experiments \cite{Nakahira2021,Liu2019}, though muscles have their own constraints, as demonstrated by Fitts Law~\cite{Nakahira2019b}. 


\citet{Nakahira2021} developed experimental tasks and corresponding sensorimotor control models that mimicked three aspects of mountain biking: compensation by the spinal cord for the random shaking coming down the trail, the anticipation of turns in the trail by the visual system, and the stabilization of images on the retina by the oculomotor system to compensate bouncing. 
Two driving experiments were performed: the first is to test the interactions between layers, and the second is to test the errors caused by delays and rate limits in control within a layer. In the two experiments, subjects follow the trail on a computer screen and control a cursor with a wheel to stay on the trail. The goal of the subjects is to minimize the errors between the desired and actual trajectories shown in a computer monitor by moving the steering wheel (see Figs.~\ref{fig:interface} and~\ref{fig:wheel}). 

In the first experiment, the higher-layer and the lower-layer are coordinated, and the authors compared how subjects' control behaviors and the resulting errors differ in three settings: 1) when there are random force disturbances to the steering wheel due to bumps on the ground (denote as 'Bump only'), 2) when the trail trajectory is curved and changes direction (denote as `Trail only'), and 3) when both exist (denote as `Both'). 
Rejection of bump disturbance in the first and last settings is likely to be performed at the lower layer reflex, while trajectory following in the second and last settings is likely to be performed at the higher layer planning. 

\sdbarfig{\centering \includegraphics[width=\linewidth]{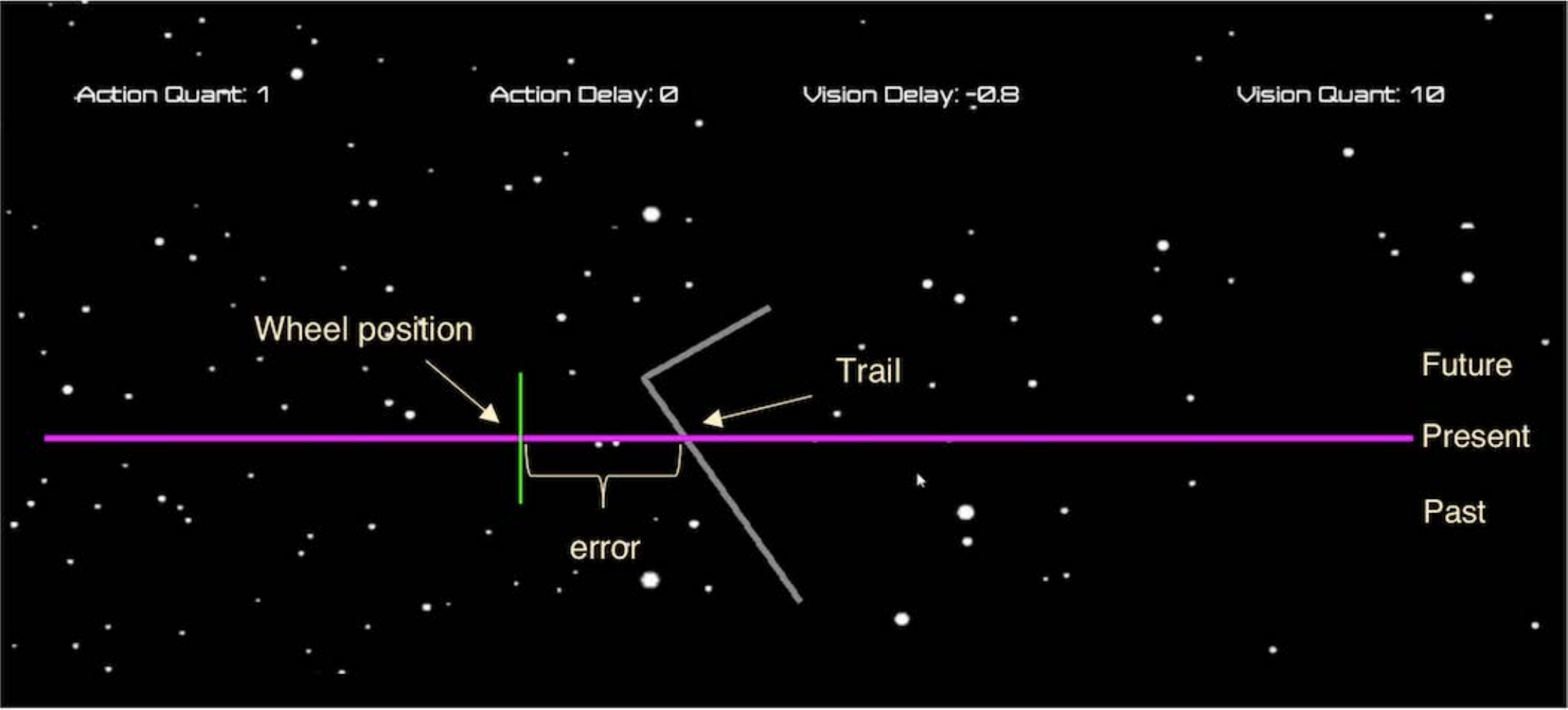}}{\footnotesize Players see a winding trail scrolling down the screen at a fixed speed, and with a fixed advanced-warning (the visible trial ahead), both of which can be varied widely. The player aims to minimize the error between the desired trajectory and their actual position using a gaming steering wheel.\label{fig:interface}}

\sdbarfig{\centering \includegraphics[width=\linewidth]{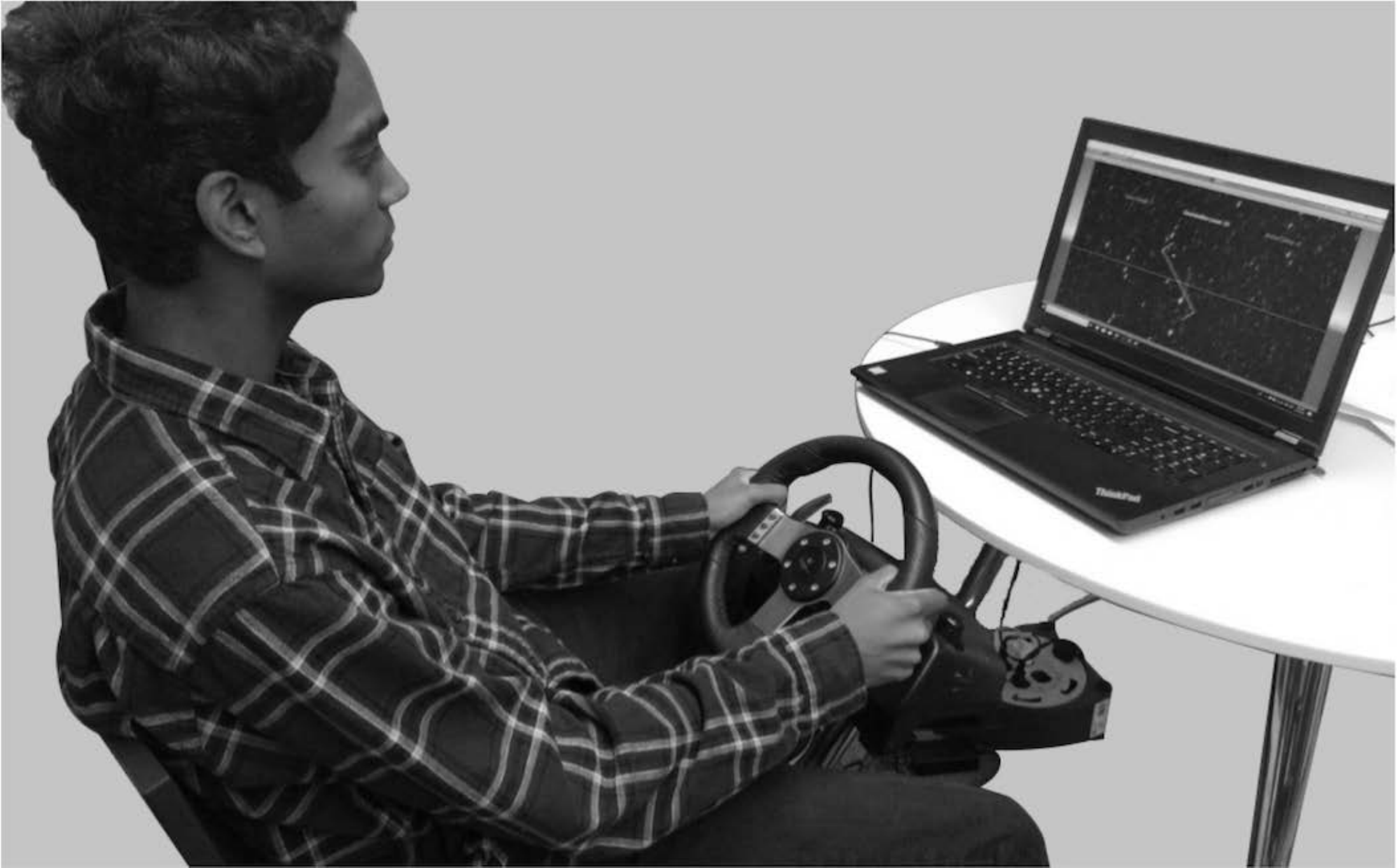}}{\footnotesize Bumps are added using a motor torque in the wheel. Experiments can be done with bumps only or trails only, or both together, and with varying trail speed and/or advanced-warning, and with additional quantization and/or time delay in the map from wheel position to players' actual position.\label{fig:wheel}}



\hfill \emph{(continued on next page)}
\end{sidebar}
\begin{sidebar}{\continuesidebar}
\renewcommand{\thesequation}{S\arabic{sequation}}
\renewcommand{\thestable}{S\arabic{stable}}
\renewcommand{\thesfigure}{S\arabic{sfigure}}

\subsection{Experimental Validation in a Biking Simulator}
The experimental results are shown in
Fig.~\ref{fig:experiment-bumps-trail}. The observed error in setting 3 (with both bumps and trail curvature) positively correlated with the sum of the errors from the first two settings with either bumps or trail curvature (Pearson correlation coefficient $=0.57$), suggesting the two signals tended to have consistent sign and amplitude. Moreover, the two signals showed no significant difference in the two-side t-test analysis. 
The results suggest that the two layers could be analyzed separately. This separability motivates the modeling of each layer separately and to further decompose the errors into those caused by neural signaling delays or rate limits in the control loop. 

The impact of neurophysiological limits was studied in the second experiment. We observed changes in lateral control error in three settings: when external delays are added in the display, when external quantizers are added in the actuation effect of the steering wheel, and when both are added. These manipulations served as noninvasive probes for how component constraints affect system behavior. The lateral errors in the three settings are shown in Fig.~\ref{fig:experiment-SAT}, and their corresponding theoretical prediction is shown in Fig.~\ref{fig:SAT-system-theory} (see the modeling details in the next section). The bridge between the constraints at the two levels highlights the benefits of the heterogeneity observed in nerves  (Fig.~\ref{fig:axonsize}) and the advantages of layering in sensorimotor control (as in Fig.~\ref{fig:visualmodel}).

\sdbarfig{\includegraphics[width=0.8\linewidth]{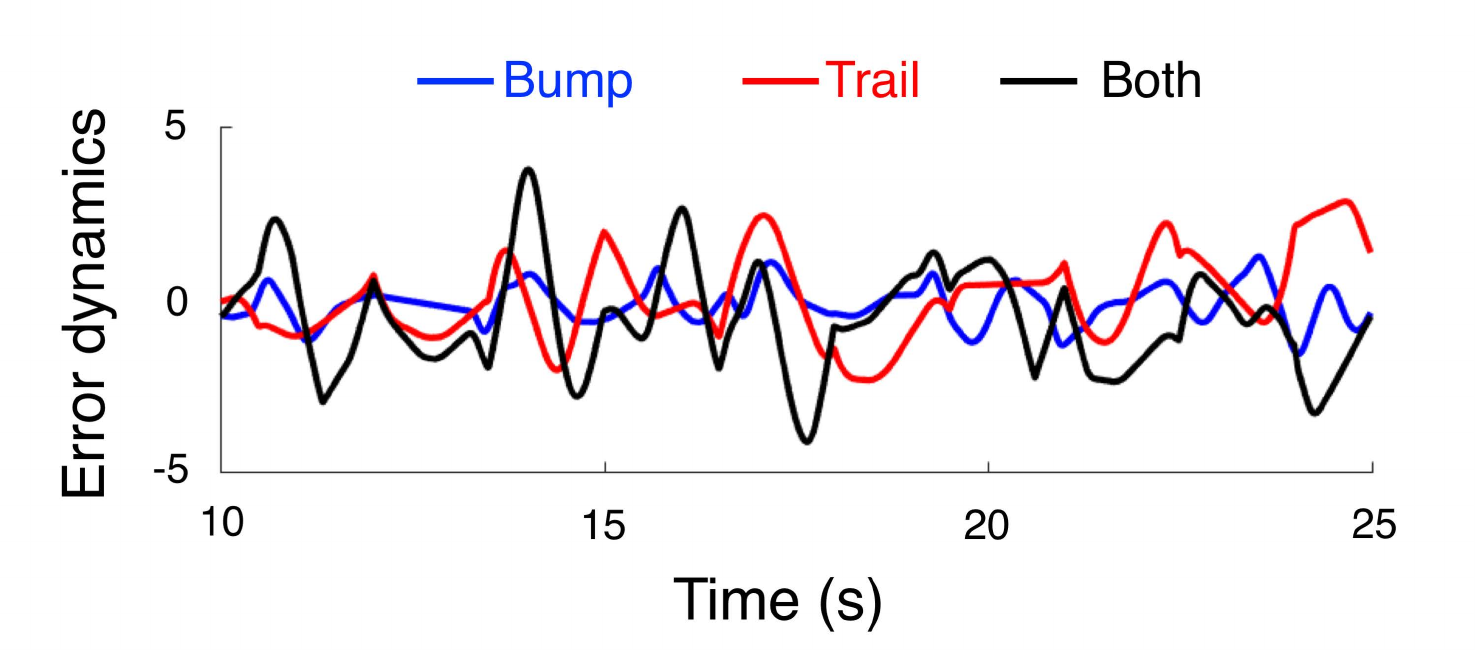}}{\footnotesize Errors in the case of bump only, trail only, and both.\label{fig:experiment-bumps-trail}}

\sdbarfig{\includegraphics[width=0.9\linewidth]{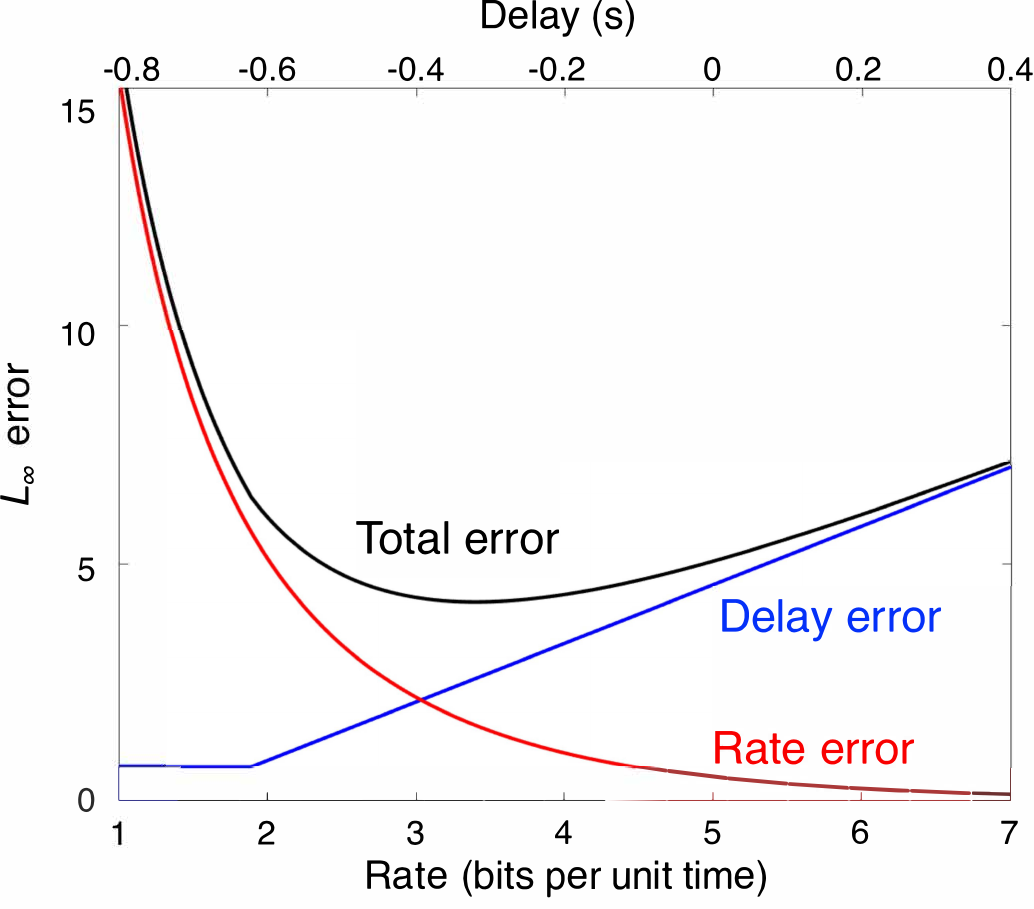}}{\footnotesize The delay error $\max(0,T)$ (blue), rate error $(2^{R}-1)^{-1}$ (red), and the total error $\max(0,T) + (2^{R}-1)^{-1}$ (black) are shown with varying component signaling delay $T_s$ and rate $R$ subject to the component constraint $T = (R-5)/20$.\label{fig:SAT-system-theory}}

\sdbarfig{\includegraphics[width=0.9\linewidth]{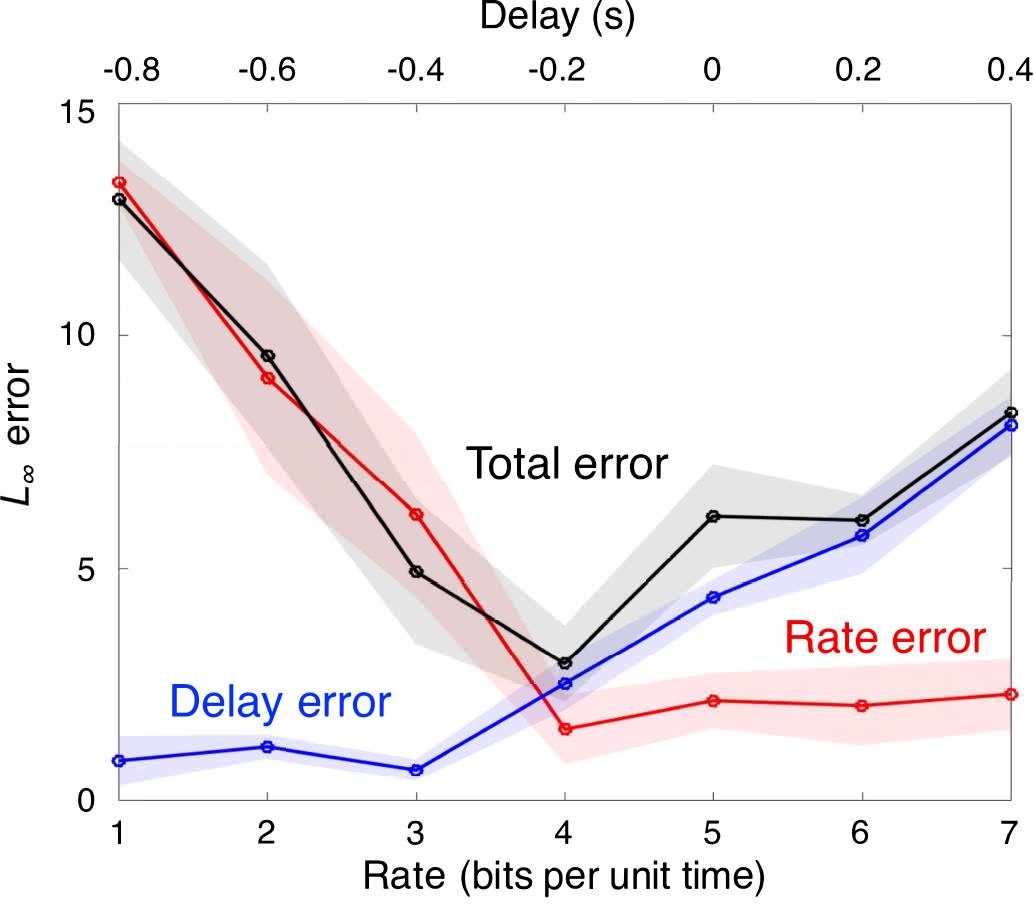}}{\footnotesize The error under an added delay (blue), the error under added quantization (red), and the error under added delayed plus quantization (black) are shown. In the last case, the added delay $T$ and quantization rate $R$ subject to the component constraint $T = (R-5)/20$. The dot shows the averaged error of 4 subjects, and the shadowed area indicates the standard error of the mean for these subjects.\label{fig:experiment-SAT} }

\end{sidebar}

\section[Key Concepts in Control Architecture]{Part 3: Key Concepts in Control Architecture}
\label{sec:LLL}

\begin{table*}
\centering
\begin{tabular}{C{1.5cm}C{2.5cm}C{2.75cm}C{1.5cm}C{3cm}}\toprule
& \textbf{Levels} & \textbf{Layers} & \textbf{Laws} & \textbf{Diversity Enabled Sweet Spots}\\ 
          & Physical  & Functional & Pareto surface & Near optimal Pareto point 
          \\\midrule
Clothing \includegraphics[width=1.5cm]{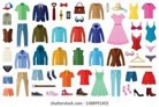}& garment\par fabric\par fibers\par thread& inner (soft, comfort)\par middle (insulation) \par outer (windproof) & warm vs. waterproof vs. soft & inner + middle + outer =  waterproof, warm, soft \\\midrule
Sensorimotor \includegraphics[width=1.5cm]{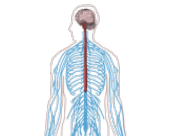}  & nerve, muscle\par axon, muscle fiber & goals\par planning\par reflex & speed vs. accuracy & vison + reflex =\par fast accurate motion\\\midrule
Power grid \includegraphics[width=1.5cm]{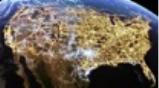} & grid\par local distribution\par tx lines, substations & economic dispatch\par $2^{\text{nd}}$-ary freq. control\par primary freq control & sustainable vs. resilient vs. efficient & Traditional + renewable + active control =\par sustainable, resilient, efficient power \\\bottomrule
\end{tabular}
\caption{Key concepts in control architectures are present across all engineered systems.  Here we illustrate these concepts using clothing, the human sensorimotor control system, and the power grid.}
\label{tab:key_concepts}
\end{table*}
 

In Parts 1 and 2, we introduced two concepts core to LCAs, namely layers and DeSS, and proposed quantitative frameworks for their analysis and design.  In this final part and section, which should be viewed as a glossary of LCA terminology, we highlight that these are but a subset of the components that can be universally found in LCAs across domains.  Although we do not have quantitative techniques for reasoning about them, we present qualitative descriptions, and illustrate their importance using various case-studies.

Table~\ref{tab:key_concepts} illustrates concepts we believe essential to the study of universal control architectures in the context of three familiar examples: \nameref{sidebar:clothing}, Sensorimotor Control, and the Power-Grid.  We also indulge in a more fanciful digression that frames \nameref{sidebar:language}.  These were introduced and developed in \cite{doyle2011architecture, Nakahira2021}, and conceptually underpin much of the previous discussion.

 
 \noindent\textbf{Levels.} Conceptually, levels can be thought of as the (usually physical) substrates or components used to implement a system. All complex systems have many levels or scales:  for example, in biology, levels range from molecules to synapses, cells, circuits, systems, and organisms. Analogous levels can be identified in familiar engineered systems.  For example, in circuits, levels range from atoms to wires, resistors, capacitors, transistors, to integrated circuits, to PCB boards. Deducing the levels experimentally is often necessary for understanding (reverse engineering) the design of existing control architectures found in nature and legacy engineered systems.

\begin{sidebar}{Clothing as a Layered Control Architecture}
\section[Clothing]{}\phantomsection
\label{sidebar:clothing}
Clothing is a familiar example that surprisingly highlights many universal concepts of control architecture \cite{doyle2011architecture}. The levels are familiar--from thread to fibers to fabric to garment to outfit and we'll focus on the latter with notation garment$\backslash$outfit to denote levels. The layers for making clothing for harsh conditions are the outer/middle/inner garments: Outer layers provide water- and wind-proofing, middle layers are insulating, and inner layers are compatible (soft) for interfacing with skin. So, layers and levels are orthogonal decompositions of outfits, and both can have further decompositions within. This architecture of clothing creates a DeSS so that outfits are weatherproof, warm, and soft when no individual garment or part provides all these features. Of note, skin and the rest of the body contains major evolved controls for thermoregulation, so that clothing can be considered an extension on top of the skin of the complex feedback controls involved in the exquisitely tight control of central temperature characteristic of healthy humans.  Adding layers in this way is an important consequence of layered architectures.

Though clothing layering is usually purely passive, the outer layer provides a barrier function to wind and rain, the mid layer provides a barrier to heat loss, and the inner layer provides a soft barrier between the possibly rough outer layers and the skin. It may seem strange to think of these as layers of passive control, but there is no other discipline that can integrate such passive mechanisms (which abound in engineering) into a full stack theory of active/passive/lossless control layers. 

A basic concept shown in the clothing example for understanding control architecture is ``barriers.'' We naturally think of active controllers as creating barriers in the state space of controller/plant feedback interconnection, and the theories of Lyapunov and barrier functions and robust control extensions are explicitly aimed to make this rigorous, useful, and scalable \cite{prajna2004safety,vidyasagar2002nonlinear,ames2016control,xu2015robustness}. What barriers in this sense allow for is showing that the set of possible controlled trajectories in state space robustly avoid ``bad'' regions. But if we want a more ``full stack'' theory of architecture where the higher levels and layers are typically active control, it will be necessary to include lower layer control that is passive or even lossless.  “Barriers” are  already familiar in studying passive controllers, but as post hoc analysis and less for design \cite{van2014port}.  We should probably think of active/passive/lossless as one example where there are both layers, e.g., in a car with active steering and braking, passive nonslip tires, and designed to be as lossless as can be in drag and friction, and levels, e.g., in a car active control is implemented in passive components plus power supplies, and physics tells us everything is microscopically lossless, which can be made rigorous using control theory \cite{sandberg2010lossless}. 

\hfill \emph{(continued on next page)}
\end{sidebar}
\begin{sidebar}{\continuesidebar}

\subsection{Clothing as a Layered Control Architecture}
Even the simplified proximal levels, layers, stages of dressing, and DeSS described here are minimal essentials to creating functional outfits, and nothing simpler will work in a harsh environment.  In particular, random piles of garments are vanishingly unlikely to make an outfit. Concretely, consider a small 30 garment wardrobe with 10 each of garments for shell$/$warm$/$soft layers. Layering allows potentially $10^3$ diverse but functional outfits, which is a much larger $n^3$ outfits versus $3n$ garments. But there is exponentially more $2^{30} ~= 1e9$ possible piles of garments and the piles/outfits ratio of $2^n/n^3$ obviously grows exponentially with $n$ garments in each layer.

One near universal in architectures is that they select functional but extremely thin and sparse subsets within the set of all possible “piles.” These thin, sparse subsets are even more extreme in the levels and layers below the garments level.  Baking is another familiar example with visible levels of ingredients and layers such as cake, frosting, crusts, filling, etc… The levels and stages of baking are explicit in a recipe, but the supply chains that provide the ingredients are typically hidden behind convenient consumer interfaces.  Random piles of ingredients and random stages of baking are extremely unlikely to produce anything even edible.

There are myriad tradeoffs and laws throughout the layers, levels, and stages that constrain what is possible, most obviously in the physical constraints on lower level materials and the high level users of the clothing architecture, but also on all the stages of supply chains. But many constraints are evolved or designed as part of the architecture, such as the fabric$\backslash$garment levels and outfit/garment layers, which were presumably not part of the earliest clothing using animal skins, even though all must obey physical laws.  These added constraints in higher layers and levels are “constraints that deconstrain” (\cite{gerhart2007theory}) in that they are essential to creating the DeSS that is the very goal of architecture.  The result is that a limited repertoire of fibers can create enormously diverse garments which are only functional due to the constraints imposed by the universal architecture used by designers, manufacturers, and users.  Baking has completely different details but is architecturally essentially the same.

This clothing architecture in harsh environments might be greatly simplified on others. Outfits in some tropical settings have one or even no layers, and garments can have a fabric made of plastics with low level polymers but no threads or fibers.  And so on.  So, diversity between architectures is as universal as the diversity that any one architecture enables, and once the centrality of this diversity is recognized, both diversities motivate an integrated theory to design and upgrade all important architectures.  But this is new and confusing even among experts, which we also hope to change.

\end{sidebar}

\noindent\textbf{Layers.} Layers are complementary to levels, and conceptually describe a functional decomposition of the overall behavior of a system.  Layered control architectures typically decompose across complexity and spatiotemporal scales, with more complex functionality implemented in higher global layers at a slower frequency, and more rigid/structured functionality implemented in lower local layers at a higher frequency, see for example Fig.~\ref{fig:three_layers}. 
Layers are the main architectural mechanism for taming complexity by breaking down a complex overall task into tractable subtasks (see \nameref{sec:layering}), and that enable DeSS \cite{Nakahira2021, Nakahira2019b} by matching the spatiotemporal resolution of each layer with a corresponding control subtask (see \nameref{sec:sweet-spots}). 

\noindent\textbf{Laws.} Almost universally, we observe that hardware components have speed-accuracy tradeoffs (SATs), which impose a law on the low level hardware that can then lead to high level laws or constraints on optimal controllers.  In neuroscience, vision is slower and more accurate than reflexes and proprioception. In immunology, adaptive immune responses take several days longer to mount than innate immune responses, but adaptive responses are more specific to the disease-causing pathogen. In computers, different storage components (e.g. registers, cache, RAM, disk) have extremely different speed, size, and cost.  Typically there are low level hardware laws from physics that can directly impact higher levels, as well as entirely new ones that arise at higher layers that have no parallel in physics and are associated with names like Turing, Shannon, and Bode.  Developing an integrated theory of laws across layers and levels is essential to a theory of architecture.

\noindent\textbf{Diversity-enabled Sweet Spots.} 
In engineering, complex system functionality requires diverse hardware, and most hardware is involved in diverse functions. If built out of homogeneous components, the SATs imposed by lower levels would make robust control impossible. However, these SATs allow for extreme diversity in the hardware, which can be leveraged with the right architectures to provide diverse functionality. Highly diverse hardware-level components (which are constrained by SATs) enable performance sweet-spots that largely overcome the severe hardware-level SATs of individual components. In computers, such sweet spots include virtual memory management systems. In neuroscience, extreme diversity in axon sizes, receptors, and neurotransmitters is abundant \cite{nakahira2015hard,Nakahira2021}, but largely hidden. By itself, diversity of components only enables sweet spots of function; to achieve these functional sweet spots requires specific architectures to maximize the utility of diverse components, which we call DeSS.  We proposed a quantitative theory of DeSS by viewing \nameref{sec:sweet-spots} and provided examples of these concepts at play.

\begin{sidebar}{Lego as a Layered Control Architecture}
\section[Lego as a Layered Control Architecture]{}\phantomsection
\label{sidebar:language}
Lego is a simple, convenient, and literally toy system that illustrates many essentials of architecture, uses conventional digital control, but has transparent processes for the supply chain to (dis)assemble toys \cite{csete2002reverse}. Consider a familiar scenario where a child is repeatedly assembling, operating, and disassembling Lego robots to build a new one, and further focus on the building of one robot from a box full of old partial robots and isolated basic parts. There are roughly 4000 diverse standard Lego parts which are produced by a manufacturing supply chain that is hidden (virtualized) from the child.   There are an infinite variety of possible robots, which are nevertheless a vanishingly small subset of all nonfunctional Lego assemblies. 

Focusing on building one robot, the minimal levels would be parts\textbackslash\!~robots consisting of the lower level parts that then make up an assembled robot, though additional levels could include various functional subassemblies.  The simplest stages would be dissassembly|parts|assembly which form a bowtie with a large but relatively thin knot of parts compared to the infinite variety of robots and assemblies as inputs and outputs.  This depends on a universal snap protocol to make both disassembly and assembly easy. Building a Lego toy is a minimal example of the classic thin knot consisting of a set of parts plus the protocols specifying how the parts can be assembled.  The most basic Lego has just one snap protocol and thousands of parts in its “knot.” Most architectures have many more of both, which are still tiny compared to the variety of systems with a shared architecture.  

This proximal part of the Lego supply chain would also have a control hourglass where a child builder would take instructions and convert them into step-by-step assembly via the snap protocol.  The thin middle waist layer would include the universal snap, here controlled repeatedly to control overall assembly.  The top layer would be the infinite possible instructions to assemble working robots, and the bottom layer would be the huge variety of supply chain steps that these instructions would control, and the robots and sub-assemblies this produces.

A universal feature this illustrates but that can be source of confusion is that the snap protocol is necessarily central to both the bowtie assembly knot and the control hourglass waist.  In the bowtie knot it is the physical mechanism that holds parts together and allows robots and their parts to be easily (dis)assembled.  This bowtie alone would be useless however, without an additional hourglass {\it control} of the snap process to in each step of (dis)assembly of a robot. The bowtie has essentially infinite diversity in the input of old robots or partial assemblies and the output of new robots.  The hourglass also has infinite diversity in the top layer of instructions and the low layer of physical assembly steps, with a thin mid-layer waist that performs snap by snap (dis)assembly according to instructions.  

\hfill \emph{(continued on next page)}
\end{sidebar}
\begin{sidebar}{\continuesidebar}
\subsection{Lego as a Layered Control Architecture}

The Lego snap is a sweet spot in the space of alternative connection and control protocols\cite{csete2002reverse}. 
One alternative is smooth bricks with no snap, which would be easier to assemble but would not be able to make robots. Another would be adding glue which would make the robot more robust to trauma but make reuse difficult.  The snap protocol is highly efficient, reusable, and robust, but fragile to finely targeted attacks such as removing imperceptibly thin and small bits of plastic just at the interface so that the snap would not hold.  The process and the built robot however would be largely robust to similar removals away from the snaps, except in the computers controlling the robot.  This extreme "robust yet fragile" feature is ubiquitous in real architectures \cite{CarlsonPNAS}, with one aspect captured in Bode's Integral Formula. 

The snap also makes it easy to manufacture new Lego parts that work with existing parts and architecture.  The knot and waist utilizing the snap protocol form the “core” of the architecture “crux” for the control of assembly, which here is done by a child infinitely more complex than any Lego robot.  This process could in principle be replaced by special purpose assembly machines not greatly more complex than the robots it builds but attempts to build a truly self-replication universal Lego robot or machine have proven challenging.  

In addition to the parts\textbackslash~\!robot levels the functioning robot has, minimally, layers of computer/(sense\&actuate)/plant where the “plant” here would be the uncontrolled raw robot.  This control system is distinct from the one doing assembly, and the computer would have sublayers of software/hardware, making this a toy version of a standard digital control system.  Note that the software would typically be vastly more complex than the rest of the robot, and computer hardware would introduce vastly more levels including microscopic components like transistors.  A robot toy without sensors, actuators, and computers would be infinitely simpler with only minimal functionality but would still have some important architectural features.  The complexity of the design process for a new robot would also be dominated by the control software, which would then be easily added in the assembly process.  This assumes the also complex computer hardware is designed and manufactured separately and arrives as a completed brick component.  The design and manufacture of the computer hardware would be vastly more complex than most robots using it as a component. 

These minimal starting points illustrate the most essential universal architectural features beyond stages, levels, and layers, including Diversity-enabled Sweet Spots (DeSS) and virtualization in both assembly and control.  There is obviously large diversity in the parts and huge diversity in the possible toys, but the integrated functionality of a built robot with a digital controller illustrates how diverse parts enable this functionality but require the specific architectural layers, level, and stages to realize this functionality.  An essential element of the DeSS is the use of virtualization, in both the assembly and control of a Lego robot.  The simplest is how the snap protocol is largely hidden in the assembled robot.  It has not disappeared completely, as disassembly would reveal, but it is completely hidden in normal operation.  This virtualizes both the parts and the assembly process so that the real-time control of the robot can ignore them.  The control layers of SW/HW/(sense\&act)/plant also has virtualization by every layer.  The sense\&act layer virtualizes the plant into an input-output system amenable to control, and the computer hardware virtualizes these details so that control can become a highly virtualized software design problem.  These hardware layers severely constrain what control is possible, but if well designed then great facilitate both the design and implementation of sophisticated control.  This creates a diversity-enable sweet spot (DeSS) where the resulting system has the flexibility and evolvability of software but the speed and accuracy in the sensing and actuation hardware of the robot plant.  Creating a DeSS is the most essential reasons why an architecture and virtualization is used at all, when no individual component alone could make up the robot or its control.

\hfill \emph{(continued on next page)}
\end{sidebar}
\begin{sidebar}{\continuesidebar}
\subsection{Lego as a Layered Control Architecture}

The big advantage of Legos as a case study is that for simple assemblies the process is transparent and doable by children, yet it illustrates many essential and universal features of architecture more generally.  It also illustrates that when active control is added via specialized parts for sensing, computing, and actuation, the complexity explodes, so that essentially all design challenges are dominated by control and software, and then the remaining physical parts only enable that control.  Together all these architectural universals create a highly virtualized system with a DeSS far beyond what any level or layer could provide by itself. This kind of efficiency, robustness, and evolvability tradeoffs addressed by virtualization dominate the design of most architectures in biology and technology, and necessarily lead successful architectures to adopt some nearly universal features.  Most are minimally present in toy Lego robots, and even more are present in the myriad cruxes in bacteria of replication, transcription, translation, metabolism, transport, and signal transduction, where all the cores have extremely conserved protocols for billions of years and even mostly conserved molecular machines.  

With an explicit inclusion of control, the gap between the complexity of Lego and free-living bacteria is enormous, where the latter makes not only all the parts and does self-replication, but controls allostasis and homeostasis in ways that typically don’t arise in robots with external supplies of parts and energy.  Nevertheless, they share striking universals from levels, stages, and layers to cruxes of bowtie stages and hourglass controls, to knot, waist, and core protocols, to virtualization and DeSS.  While there is no comparable universal terminology, we are proposing one here that is aimed to be maximally consistent with those specialized domains that do explicitly consider architecture.  Bacteria and Legos surprisingly illustrate the most essential universals, but a large variety of other less familiar domains could as well.  Particularly for experts in many domains of bio, med, neuro, and tech systems, there are equally rich if less accessible examples of universal architectures.  

Bacteria, however, are the original from which all else has evolved, and remain arguably the most perfect.  Their robustness and evolvability are due to the universal architectures that they share with all lineages from them, but their fragilities to hijacking are also devastatingly universal.

\end{sidebar}
\subsection{Bowties, hourglasses, virtualization, and abstraction}
Fortunately, some features of LCAs are very familiar, particularly universal \emph{bowties} and \emph{hourglasses} that appear in complex highly evolved systems at every scale and context. In both bowties and hourglasses two outer deconstrained stages and layers, with very diverse components that are evolvable and even swappable, are linked in the middle via a narrow, highly constrained knot/waist with little diversity or evolvability.  We call this \emph{constraints that deconstrain} (as in \cite{gerhart2007theory}).  The terminology of bowties and hourglasses is not standard and can be confusing, but the distinction between them is useful and important.  For a biologically motivated case-study, see \nameref{vig:bacteria}. 
Both the bowtie and hourglass enable virtualization via universal shared interfaces, like OSes, ATP, wall plugs, this text, ribosomes and translation, HTML, TCP/IP, HDMI, membrane potentials, faucets, dashboards, etc.

\noindent\textbf{Bowtie.} Diversity is the aspect of architectures that is most familiar and easiest to discuss in the stages making up supply chains.  Diverse proteins are a produced by highly conserved translation 
"knot" protocols with amino acid inputs and controlled by a transcription hourglass.  In metabolism, diverse carbon sources and molecules are linked via a thin ``knot'' of a few metabolic carriers and precursors.  Diverse electric power sources and user appliances are linked in a bowtie via standard knot protocols (e.g., 110v 60Hz) in power grids.  These examples all involve the flow of materials and energy thru various stages and, with respect to diversity, have a bowtie shape, with diverse sources and products at the edges and highly conserved and less diverse ``knots'' in the middle. This enables independent and thus rapid evolution on both ends of the bowtie.

\noindent\textbf{Hourglass.} An hourglass is used to describe the shape of layered communication and computing systems required to control bowties.  Diverse app software runs on diverse hardware in an hourglass linked via less diverse ``waist'' operating systems (OS) in computers and their networks.  Humans have diverse skills and memes and diverse tools, linked in an hourglass by shared languages and a poorly understood brain OS. Genes, apps, memes, words, technologies, and tools are highly modular and swappable, massively accelerating evolvability beyond what is possible with only the slow accumulation of small innovations.

\noindent\textbf{Virtualization.} Hourglasses rely on virtualization to enable the diversity both above and below the hourglass ``waist.''  For example, operating systems in computers act as a \emph{protocol that virtualizes} the wildly diverse hardware and computer networks in modern computing systems, which in turn has lead to the incredible progress and diversity of software and data.  In decision and control systems, low-level unstable dynamics are virtualized by the feedback control layer, allowing the planning layer to use simple, reduced order, and stable models for trajectory generation.  Indeed, a commonly used model for trajectory generation in robotics across a wide variety of platforms (e.g., quadrupeds, quadrotors, mobile robots) is the Dubins' car or unicycle model---we expounded on this particular example of virtualization in robotics in previous sections (see Fig.~\ref{fig:multi-rate} in \nameref{sidebar:multirate} and Ex.~\ref{ex:mpcdubins}).  Here, the reference trajectory serves as the protocol between diverse planning and control layers, wherein each can be constructed using a diversity of algorithms, abstractions (see below), hardware, and software.

\noindent\textbf{Abstraction:} Whereas the implementation of layered control architectures is enabled by bowties, hourglasses, and virtualization, the \emph{design} of layered architectures would be impossible without \emph{abstractions}.  
For example, when writing computer software engineers abstract OS/HW as memory and compute, often ignoring for example, device level drivers and timing constraints.  Note however that as software approaches the limits of what the underlying hardware can implement, these abstractions may no longer be valid; hence the need for, for example, real-time programming languages for embedded systems that directly access hardware resources. In decision and control systems, abstractions abound.  At the feedback control layer, the plant and controller are abstracted as mathematical operators operating on continuous- or discrete-time signals. At the trajectory planning layer, the potentially complex low-level closed-loop control system is abstracted using a simple dynamics model, e.g., a unicycle.  This abstraction is valid thanks to the virtualization enabled by the feedback control layer below, but also breaks down if the planned trajectories extend beyond the tracking capabilities of the closed-loop system, again showing that abstractions are useful only within operating ranges that virtualization can be reliably enforced.  Because virtualization greatly enables the use of effective abstractions, these two distinct concepts are often confused.

\begin{sidebar}{Bowties and Hourglasses in Bacterial Metabolism}
\setcounter{sequation}{0}
\renewcommand{\thesequation}{S\arabic{sequation}}
\setcounter{stable}{0}
\renewcommand{\thestable}{S\arabic{stable}}
\setcounter{sfigure}{0}
\renewcommand{\thesfigure}{S\arabic{sfigure}}
\section[Bowties and Hourglasses in Bacterial Metabolism]{}\phantomsection
   \label{vig:bacteria}

   \sdbarfig{\label{fig:bacteria-bowtie}
    \includegraphics[width=\columnwidth]{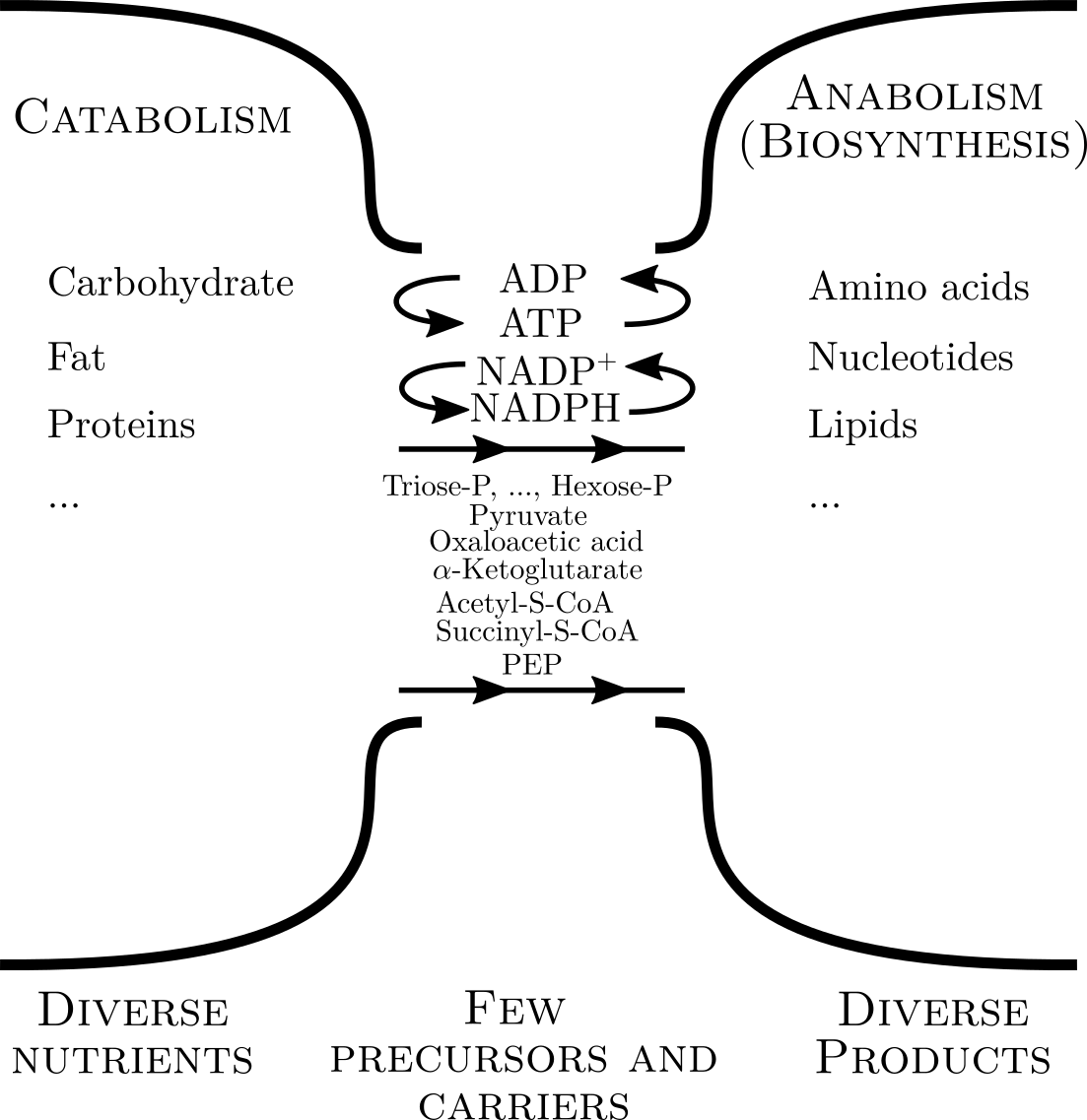}}{\footnotesize Bowtie in the low layer of metabolism stoichiometry in bacterial cells.}

An example that naturally embodies the structural features of layered architectures is the organization in the bacterial cell. A bacterial cell performs diverse types of complex functions on many timescales, from digesting nutrients and synthesizing macromolecules to adapting to environmental disturbances to cell cycling and decision making, to long term evolution. This wide range of functions is fundamentally 
enabled by the layered architecture with bowties and hourglasses (see Figs.~\ref{fig:bacteria-bowtie} and~\ref{fig:bacteria-hourglass}).

The lowest \textit{layer} or plant consists of metabolites connected by reactions and summarized in a stoichiometry matrix.  The overall organization has a \textit{bowtie} with very diverse input and output stages and a thin low diversity knot stage of precursors and carriers (ATP, NADH, ...) that are then cofactors throughout. Catabolic pathways convert input nutrients to knots and then anabolic pathways make output products. 
%
We can also crudely view a bacterial cell as having two layers, with a low layer of metabolism, and a high layer of gene expression and then add structural features within a layer and between layers in the bacterial cell.
Namely, each layer has a bowtie shape with a small knot (of carriers and precursors for metabolism) that connects diverse inputs and outputs on both sides.
The high layer, viewed as a controller of the low layer, has an hourglass shape, with a thin universal waist (OS-like of transcription and translation) controlling and connecting  diverse high layer genes to diverse low layer actuation by proteins. Bowties and hourglasses are further universal features of complex architectures but are hidden and cryptic in normal operation, enabling efficiency and flexibility but potentially hiding fragilities.

A bacterial cell's metabolism layer obtains energy and materials from nutrients in the environment using enzymes to catalyze reactions.
The stoichiometry captures the structure but not rates of these reactions, and has a bowtie shape, see Fig.~\ref{fig:bacteria-bowtie}. Diverse nutrient molecules are digested in the catabolism stage, and diverse macromolecules are synthesized in the anabolism stage, but the intermediate ``knot'' is a very thin stage of a few precursors and carriers.
The ATP/ADP pair is the carrier for energy from catabolism to use for anabolism.

As metabolism happens on a fast timescale that is intrinsically unstable, regulation of these rapid reactions is needed locally because of delays in diffusion.
Local enzymatic regulations through binding reactions serve as local actuators to be further controlled by higher layers. With these local enzymatic regulations stably maintaining a steady state of the cell's metabolism, this establishes a supply chain of molecules for energy, redox potential, and molecular building blocks used to perform tasks at a higher layer.
This higher layer then can perform dynamics that take this supply chain as given and focus on goals with a virtualized molecular supply chain.
For example, gene expression is one such layer. Here building blocks such as nucleic acids and amino acids are used to build up large molecules such as RNAs and proteins. 
The dynamics of gene expression can then focus on which RNAs and proteins are produced when and where, without worrying about the the supply chain of building blocks or energy for synthesis.
This ``digital layer'' is in contrast to those focusing on energy, redox potential, and molecular concentrations in the metabolism layer.

\hfill\emph{(continued on next page)}

 \end{sidebar}

   \begin{sidebar}{\continuesidebar}
   \renewcommand{\thesfigure}{S\arabic{sfigure}}

\subsection{Bowties and hourglasses in bacterial metabolism}

While the higher gene expression layer is fundamentally supplied and enabled by the lower metabolism layer, this higher layer regulates the lower metabolism layer on a slower time-scale. 
While metabolic reactions tend to happen faster than seconds, gene expressions tend to happen in tens of minutes.
The gene regulatory network can make more complex decisions and change the enzyme compositions precisely to actuate and coordinate at a global scale on the metabolism layer.
For example, while rapid fluctuations in the ATP concentration need to be stabilized by local enzymatic feedback, a shift in nutrient source requires the coordination on the gene expression level to stop expressing enzymes for old nutrients and start expressing enzymes to digest new nutrients.

\sdbarfig{\label{fig:bacteria-hourglass}
    \includegraphics[width=\columnwidth]{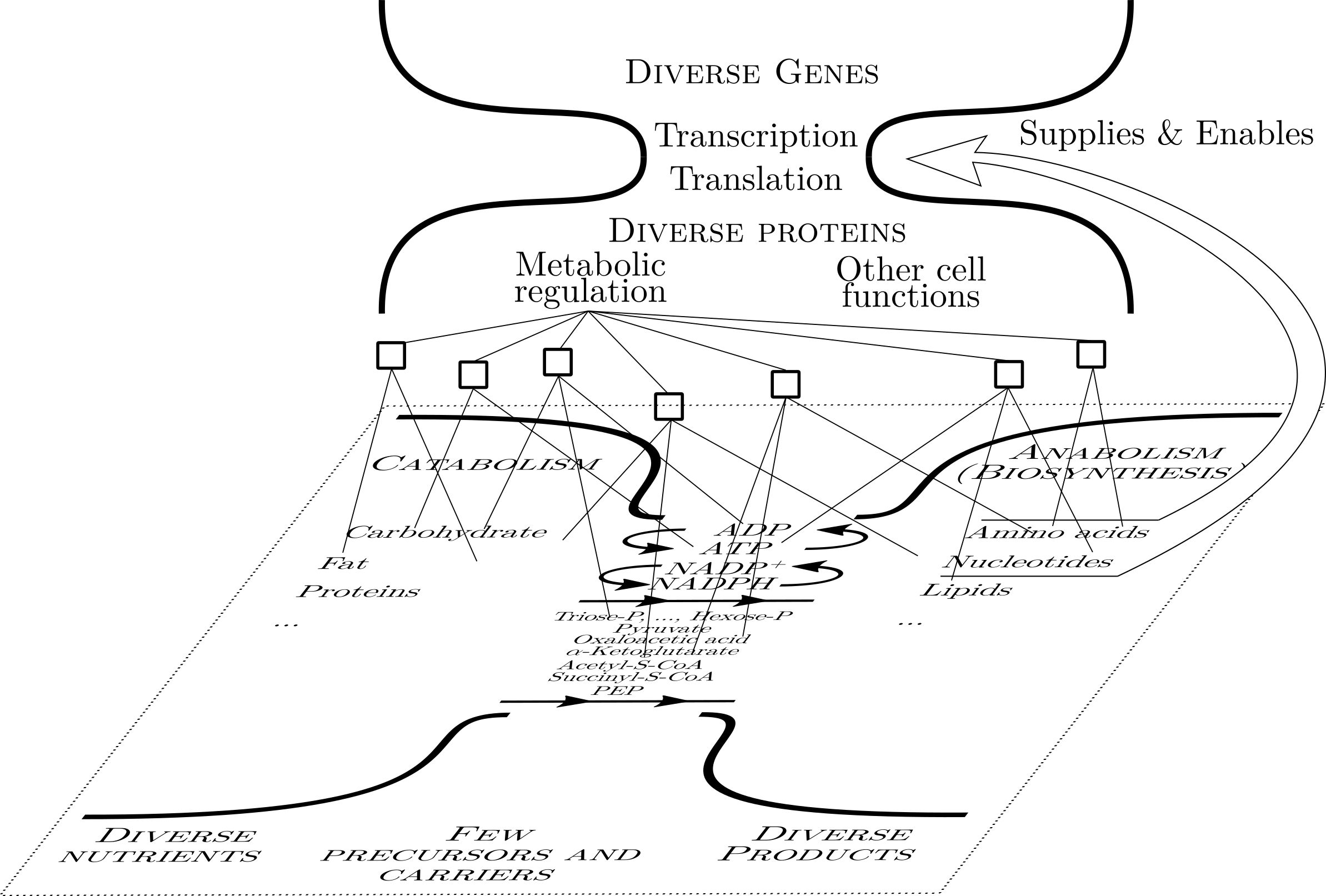}}{\footnotesize Hourglass in gene expression layer to control the metabolism stoichiometry layer. Squares represent enzymes that locally regulate some metabolic reactions. Graphics illustrating the structures in layered architectures of bacterial cell. \textbf{(a)} The low layer is the stoichiometry of bacterial metabolism, with a bowtie shape. \textbf{(b)} The metabolic reactions are locally regulated by enzymatic binding reactions such as allostery (squares), which is in turn regulated by the high layer of gene expression. The gene expression layer viewed as a controller for the low stoichiometry layer has an hourglass shape, connecting diverse genes with diverse proteins via a thin waist of transcription and translation machinery. The gene expression layer regulates the low metabolism stoichiometry layer, but also takes supply from and is enabled by the metabolism stoichiometry}   

In order to implement the gene expression layer's complex and diverse control of the metabolism layer, the cell organizes the gene expression controller in an hourglass shape, see Fig.~\ref{fig:bacteria-hourglass}. Diverse signals in the form of combinatorial gene activation are mapped to diverse actions in expressed enzymes and other regulatory proteins via a thin waist that is the universal protocol of transcription-translation machinery.
This hourglass structure is essential for the gene expression layer's control actions to scale up and facilitate diversity, namely coping with diverse and complex disturbances and performing diverse and complex actions on and via enzymes.
\end{sidebar}

\section{Conclusions}
\label{sec:conclusion}
We introduced a lexicon for key concepts in layered (control) architectures---levels, layers, stages, laws, DeSS, hourglasses, bowties, virtualization, abstraction---and instantiate them in familiar and diverse examples such as clothing, bacteria, GNC, robotics, and human sensorimotor control. These concepts are mostly familiar, but are referred to using different terms across domains: thus one primary goal of this manuscript was to establish a common language to unify the study of architecture.  Furthermore, for certain concepts, we also proposed quantitative frameworks for the analysis and synthesis of LCAs, grounded in robotics and sensorimotor applications.  

We are very much aware that this paper poses more questions than it answers, and is likely to confuse (and perhaps even anger) applied and theoretical researchers alike.  Nevertheless, we believe that underneath the cumbersome jargon and mathematical notation needed to convey our message there is a viable path towards a quantitative and universal theory of layered control architectures that the controls community is particularly well suited to pursue.  With that in mind, we hope that if the reader leaves this paper with but one core message, it is that complex systems are composed of diverse levels and layers, and that their analysis and design falls squarely within the skill set and expertise of the controls community.  Indeed, the impact in both theory and application of nascent versions of these concepts has already been astounding both within and outside of our community, and we are excited about the potential future impact that this emerging field of study will have.

\section{ACKNOWLEDGMENT}
The work of N. Matni was supported in part by NSF awards CPS-2038873, EECS-2231349, SLES-2331880, and NSF CAREER award ECCS-2045834.  The work of A. D. Ames is supported in part by NSF award CPS-1932091. 

\section{Author Information}
\begin{IEEEbiography}{Nikolai Matni}{\,}(Senior Member, IEEE) is an Assistant Professor in the Department of Electrical and Systems Engineering at the University of Pennsylvania, where he is also a member of the Department of Computer and Information Sciences (by courtesy), the GRASP Lab, the PRECISE Center, and the Applied Mathematics and Computational Science graduate group.  He has held positions as a Visiting Faculty Researcher at Google Brain Robotics, NYC, as a postdoctoral scholar in EECS at UC Berkeley, and as a postdoctoral scholar in the Computing and Mathematical Sciences at Caltech. He received his Ph.D. in Control and Dynamical Systems from Caltech in June 2016. He also holds a B.A.Sc. and M.A.Sc. in Electrical Engineering from the University of British Columbia, Vancouver, Canada. His research interests broadly encompass the use of learning, optimization, and control in the design and analysis of autonomous systems.  Nikolai is a recipient of the AFOSR YIP Award (2024), NSF CAREER Award (2021), a Google Research Scholar Award (2021), the 2021 IEEE CSS George S. Axelby Award, and the 2013 IEEE CDC Best Student Paper Award.  He is also a co-author on papers that have won the 2022 IEEE CDC Best Student Paper Award and the 2017 IEEE ACC Best Student Paper Award.
\end{IEEEbiography}

\begin{IEEEbiography}{Aaron D. Ames}{\,} (Fellow, IEEE) received the B.S. degree in mechanical engineering and the B.A.
degree in mathematics from the University of
St. Thomas, Saint Paul, MN, USA in 2001, and the
M.A. degree in mathematics and the Ph.D. degree in
electrical engineering and computer sciences from
the University of California, Berkeley, CA, USA,
in 2006. From 2006 to 2008, he served as a Post-Doctoral
Scholar in control and dynamical systems with
the California Institute of Technology (Caltech),
Pasadena, CA, USA. In 2008, he began his faculty career at Texas A\&M
University, College Station, TX, USA. He was an Associate Professor with
the Woodruff School of Mechanical Engineering and the School of Electrical
and Computer Engineering, Georgia Institute of Technology, Atlanta, GA,
USA. Since 2017, he has been a Bren Professor of Mechanical and Civil
Engineering and Control and Dynamical Systems at Caltech. His research
interests include the areas of robotics, nonlinear, safety-critical control, and
hybrid systems, with a special focus on applications to dynamic robots—both
formally and through experimental validation.
Dr. Ames was a recipient of the 2005 Leon O. Chua Award for Achievement
in Nonlinear Science and the 2006 Bernard Friedman Memorial Prize in
Applied Mathematics from the University of California, Berkeley. He received
the NSF CAREER award in 2010, the 2015 Donald P. Eckman Award, and
the 2019 IEEE CSS Antonio Ruberti Young Researcher Prize.
\end{IEEEbiography}

\begin{IEEEbiography}{John C. Doyle}{\,}received the B.S. and M.S. degrees in
electrical engineering from the Massachusetts Institute
of Technology, Cambridge, MA, USA, in 1977, and
the Ph.D. degree in mathematics from UC Berkeley,
Berkeley, CA, USA, in 1984.
He is currently the Jean-Lou Chameau Professor of Control and Dynamical Systems, Electrical
Engineer, and Bio-Engineering, Caltech, Pasadena,
CA, USA. His research interests include mathematical foundations for complex networks with applications in biology, technology, medicine, ecology,
neuroscience, and multiscale physics that integrates theory from control,
computation, communication, optimization, statistics (e.g., machine learning).
Dr. Doyle received the 1990 IEEE Baker Prize (for all IEEE publications),
also listed in the world top 10 most important papers in mathematics
1981-1993, the IEEE Automatic Control Transactions Award (3x 1988,
1989, 2021), the 1994 AACC American Control Conference Schuck Award, the
 ACM Sigcomm 2004 Paper Prize and 2016 Test of Time Award, and inclusion
in Best Writing on Mathematics 2010.  Individual awards include 1977
IEEE Power Hickernell, 1983 AACC Eckman, 1984 UC Berkeley Friedman,
1984 IEEE Centennial Outstanding Young Engineer (a one-time award for
IEEE 100th anniversary), 2004 IEEE Control Systems Field Award, and world records and championships in various sports.  
\end{IEEEbiography}

\bibliographystyle{unsrtnat}
\bibliography{arch_references,ames,references, comms, decentralized}
\endarticle

\begin{appendix}
    \paragraph{Dual to problem~\eqref{eq:minimax-ls}}
    \label{sec-app:dual}
    Recall that the minimax least-squares problem is given by:
    \begin{equation}\label{app:minimax-ls}
    \begin{array}{rl}
    \mathrm{minimize}_{x} & \max\{\|A_1 x - b_1\|_2^2, \|A_2 x - b_2\|_2^2\},
    \end{array}
\end{equation}

In order to derive an interesting dual problem, we consider the following equivalent problem
\begin{equation}\label{app:minimax-ls-equiv}
    \begin{array}{rl}
    \mathrm{minimize}_{x_1,x_2,t} & t \\
    \text{subject to} & \|A_1 x_1 - b_1\|_2^2 \leq t, \\
    & \|A_2 x_2 - b_2\|_2^2 \leq t, \\
    & x_1 = x_2.
    \end{array}
\end{equation}

The corresponding Lagrangian is given by
\begin{multline}
L(t, x_1, x_2, \lambda_1, \lambda_2, \nu) = t + \lambda_1(\|A_1 x_1 - b_1\|_2^2 - t) + \\ \lambda_2(\|A_2 x_2 - b_2\|_2^2 - t) + 2\nu^T(x_1-x_2),
\end{multline}
for $\lambda_1,\lambda_2>0.$

We rewrite the Lagrangian in the more suggestive form
\begin{multline}
L(t, x_1, x_2, \lambda_1, \lambda_2, \nu) = t(1-\lambda_1-\lambda_2) \\+ \lambda_1\|A_1 x_1 - b_1\|_2^2 + 2\nu^Tx_1 \\ + \lambda_2\|A_2 x_2 - b_2\|_2^2 - 2\nu^Tx_2.
\end{multline}
Recalling that the dual function is defined as
$$
g(\lambda_1, \lambda_2, \nu) = \inf_{t,x_1,x_2} L(t, x_1, x_2, \lambda_1, \lambda_2, \nu),
$$
and that the dual problem is given by
$$
\mathrm{maximize}_{\lambda_1,\lambda_2>0, \nu} \ g(\lambda_1, \lambda_2, \nu),$$
we immediately conclude that $g(\lambda_1,\lambda_2,\nu)$ is bounded below if and only if $\lambda_1 + \lambda_2 = 1$, $\nu \perp \mathrm{ker}(A_1) \iff \nu = A_1^T \mu_1$ for unconstrained $\mu_1$, and $\nu \perp \mathrm{ker}(A_2) \iff \nu = A_2^T \mu_2$ for unconstrained $\mu_2$.  Under these conditions, the minimizers are $x_i^* = -A_i^\dag\left(\frac{\mu_i}{2\lambda_i}+b_i\right)$, which after simplification and collecting like terms, results in the dual problem~\eqref{eq:minimax-ls-dual}.

\paragraph{Proof of Theorem~\ref{thm:sigma}} Under the assumptions of the Theorem, the set of feasible $\mu_1$ and $\mu_2$ in the dual problem~\eqref{eq:minimax-ls-dual} is given by $\mu_1 = \mu_2 = E_k \alpha$, for any $\alpha \in \R^k$. It follows from the KKT conditions of the dual problem~\eqref{eq:minimax-ls-dual} that the optimal solution is $(\mu_1^\star, \mu_2^\star, \lambda_1^\star, \lambda_2^\star) = (E_kE_k^T(b_1-b_2)/4, E_kE_k^T(b_1-b_2)/4, 1/2, 1/2)$, and that the corresponding $\sigma$-sweet-spot of the problem satisfies
    $$
    \sigma = \frac{1}{4}\|E_k^T(b_1-b_2)\|_2^2.
    $$

\paragraph{Experimental details for Fig.~\ref{fig:bicriterion} }
We set $m=5$ and $n=10$, and draw the entries of $A_1$ and $A_2$ i.i.d. according to a standard normal distribution (duplicating shared rows across $A_1$ and $A_2$ as necessary).  This ensures linear independence of rows with probability 1.  We draw the entries of $b_1$ i.i.d. according to a standard normal, and set $b_2 = b_1 + \Delta$, where $\Delta \sim{} \mathcal{N}(0,100 I_m).$

\end{appendix}

\end{document}